\documentclass[a4paper]{article}
\usepackage{authblk}
\usepackage{geometry}
\usepackage{graphicx}
\usepackage{graphbox}
\usepackage{caption}
\usepackage{subcaption}
\usepackage{amsmath,amssymb,amsthm}
\usepackage{amssymb}
\usepackage{bm}
\usepackage{commath,mathtools}
\usepackage{hyperref}
\usepackage{doi}
\usepackage[square,numbers]{natbib}
\usepackage[title,titletoc]{appendix}
\hypersetup{ pdfborder = 0 0 0, colorlinks=false,   linkcolor=black,citecolor=black, filecolor=black, urlcolor=black, }
\usepackage{enumitem}
\usepackage{rotating}
\usepackage{xcolor}
\usepackage{tikz}
\usetikzlibrary{math,calc,arrows.meta}

\date{}

\usepackage[capitalize]{cleveref} % Load cleveref last!

\newtheorem{lemma}{Lemma}
\newtheorem{remark}{Remark}

\newcounter{estno}
\newtheorem{estimate}[estno]{Error estimate}
\crefname{estimate}{Error estimate}{Error estimates}

\title{Quadrature error estimates for layer potentials evaluated near curved surfaces in three dimensions} 
   
\author{Ludvig af Klinteberg}
\author{Chiara Sorgentone}
\author{Anna-Karin Tornberg}

\affil{Department of Mathematics, KTH Royal Institute of Technology, Stockholm, Sweden}
 
\begin{document}
\maketitle
%\tableofcontents

%\section*{Notation}
\renewcommand{\v}[1]{\bm{#1}}
\newcommand{\pbar}{\bar{p}}
\newcommand{\pars}[1]{\left( #1 \right)}%
\newcommand{\braces}[1]{\left\{ #1 \right\}}%
\newcommand{\stdint}{{\mathcal I}}
\newcommand{\conj}[1]{\overline{#1}}
\newcommand{\trig}{\operatorname{F}}
\newcommand{\poly}{\operatorname{P}}
\newcommand{\taylor}{\operatorname{T}}
\newcommand{\reals}{\mathbb R}
\newcommand{\C}{\mathbb C}
\newcommand{\zplus}{{\mathbb Z}^+}
\newcommand{\Surf}{S}
\newcommand{\opint}{\operatorname{I}}
\newcommand{\opintTh}{\opint[\Theta_p](\v x)}
\newcommand{\opquadsym}{\operatorname{Q}}
\newcommand{\opremsym}{\operatorname{E}}
\newcommand{\opquad}{\opquadsym_{n}}
\newcommand{\oprem}{\opremsym_{n}}
\newcommand{\opremTh}{\oprem [\Theta_p](\v x)}
\newcommand{\opquads}{\opquadsym_{s,n_s}}
\newcommand{\oprems}{\opremsym_{s,n_s}}
\newcommand{\opquadt}{\opquadsym_{t,n_t}}
\newcommand{\opremt}{\opremsym_{t,n_t}}
\newcommand{\opremst}{\opremsym^2_{n_s,n_t}}
\newcommand{\Res}{\operatorname{Res}}
\renewcommand{\Re}{\operatorname{Re}}
\renewcommand{\Im}{\operatorname{Im}}
\newcommand{\imbrac}[1]{\Im \left[ #1 \right]}
\newcommand{\rebrac}[1]{\Re\left[#1\right]}
\newcommand{\est}{\operatorname{est}}
\newcommand{\gtilde}{{\tilde\gamma}}

%\begin{itemize}
%\item $\gamma(t)$ or $\gamma(s,t)$, parametrization into $\mathbb R^3$
%\item $\stdint = [-1,1]$, standard interval
%\item $\Gamma = \gamma(\stdint)$, curve in $\mathbb R^3$
%\item $Q = \gamma(\stdint^2)$, quadrilateral face in $\mathbb R^3$
%\item $K(\v x, \v y) = k(\v x, \v y) \norm{\v x - \v y}^{-2p}$, prototype kernel ($k$ smooth)
%\item $R^2(t) = \norm{\v\gamma(t) - \v x}^2$, squared distance function
%\item $t_0$, root s.t. $R^2(t_0)=0$
%\item Integral: $\opint[f]=\int f(t) \dif t$
%\item Quadrature: $\opquad[f]=\sum_i f(t_i) w_i$
%\item Error: $\oprem[f] = \opint[f] - \opquad[f]$
%\end{itemize}

\begin{abstract}
%  Why? What? How?
The quadrature error associated with a regular quadrature rule for
evaluation of a layer potential increases rapidly when the evaluation
point approaches the surface and the integral becomes nearly singular.
Error estimates are needed to determine when the accuracy is
insufficient and a more costly special quadrature method should be
utilized. 

The final result of this paper are such quadrature error estimates
for the composite Gauss-Legendre rule and the global trapezoidal rule,
when applied to evaluate layer potentials defined over smooth curved
surfaces in $\reals^3$. The estimates have no unknown coefficients and
can be efficiently evaluated given the discretization of the surface,
invoking a local one-dimensional root-finding procedure. 
They are derived starting with integrals over curves,
using complex analysis involving contour integrals, residue calculus
and branch cuts. By complexifying the parameter plane, the theory can
be used to derive estimates also for curves in $\reals^3$. These
results are then used in the derivation of the estimates for integrals
over surfaces.  
In this procedure, we also obtain error estimates for layer potentials
evaluated over curves in $\reals^2$. Such estimates combined with a local
root-finding procedure for their evaluation were earlier derived for
the composite Gauss-Legendre rule for layer potentials written in
complex form \cite{AfKlinteberg2018}. This is here extended to provide quadrature error
estimates for both complex and real formulations of layer potentials,
both for the Gauss-Legendre and the trapezoidal rule. 

Numerical examples are given to illustrate the performance of the
quadrature error estimates. The estimates for integration over curves
are in many cases remarkably precise, and the estimates for curved
surfaces in $\reals^3$ are also sufficiently precise, with
sufficiently low computational cost,  to be practically
useful.

\end{abstract}

%\begin{keywords}
%  Layer potential, Close evaluation, Quadrature, Nearly singular, Error estimate
%\end{keywords}
%
%\begin{AMS}
%  65R20, 65D30, 65D32, 65G99
%\end{AMS}
  
\section{Introduction}

Accurate evaluation of layer potentials is crucial when solving
partial differential equations using boundary integral methods. When an
evaluation point is close to the boundary,
%the integral defining such a layer potential is nearly singular.
the integral defining such a layer potential {has a rapidly varying
integrand. }
A regular quadrature method will then yield large errors,
and a specialized quadrature method must be used to keep errors low. 
There is a variety of specialized quadrature methods, but the
increased accuracy that they can provide comes at an additional
computational cost. It is therefore desireable to have error estimates
for the regular quadrature method that can be used to determine when
the accuracy will be insufficient and a special quadrature method must be applied.

In this paper, we study the errors incurred when using two
quadrature methods that are commonly applied to evaluate layer
potentials: the panel based Gauss-Legendre quadrature rule and the
global trapezoidal rule. 
The simplest example of a layer potential in 3D is the harmonic single layer
potential
\begin{equation}
S^{3D}_H[\sigma](\v x)=\int_{\Surf}\frac{\sigma(\v y)}{\norm{\v y - \v x}} \, \, \dif S(\v y),
\label{eq:3D_harmonic_sgl}
\end{equation}
but we will consider the more generic form 
\begin{align}
  u(\v x) = \int_{\Surf} \frac{k \pars{\v x, \v y}\sigma(\v y)
  }{\norm{\v y - \v x}^{2p}} \dif S(\v y) ,  \quad p=1/2,3/2,\ldots
\label{eq:layerpot_S}
\end{align}
where $k \pars{\v x, \v y}$ and $\sigma(\v y)$ are assumed to be
  smooth, {and the evaluation point $\v x \in \reals^3$ can be arbitrarily close to, but not
on, the surface $S$. The surface $\Surf$  in $\reals^3$, is supposed to be smooth over each separate
section (panel or other) that a quadrature rule will be applied
to.}

{Layer potentals in 2D can also be written in the generic form
(\ref{eq:layerpot_S}).} Here, $S$ is now instead a curve
in $\reals^2$, and $p=1,2,3,\ldots$.
%%and we will derive error estimates also for this case. 
An example is the harmonic double
layer potential, in 2D given by (with $\v n_y$ the normal vector at
$\v y \in \Surf$),
\[
D^{2D}_H[\sigma](\v x)=\int_{\Surf}\sigma(\v y) \frac{{\v n_y} \cdot (\v y - \v x)}{\norm{\v y - \v x}^2} \, \, \dif S(\v y).
\]

{Now, let $\v y^*$ be the closest point to $\v x$ on
$\Surf$. The closer $\v x$ is to $\Surf$, the more preaked the
integrand in (\ref{eq:layerpot_S}) will become around $\v y^*$ due to the factor $\norm{\v y -
  \v x}^{-2p}$. Analytically, the integral is well defined, but numerically it
will be difficult to approximate well. Following the terminology of
Elliott et al., \cite{Donaldson1972,Elliott2008,Elliott2011,Elliott2015}, we will henceforth call this a \emph{nearly singular
  integral}.  }

To start with the basics, 
consider a simple example of a 1D-integral 
\[
\int_{-1}^{1} \frac{1}{t^2+b^2} \, dt, 
\]
which will be nearly singular when $b>0$ is small. We can e.g. apply
an $n$-point Gauss-Legendre quadrature rule to approximate this
integral. The error will be large when $b$ is small, but decrease
exponentially with increasing values of $b$.  The classical error
estimate, available in e.g. Abramowitz and Stegun
\cite[eq. 25.4.30]{Abramowitz1972} or in the DLMF
\cite[\S3.5(v)]{NIST:DLMF}, includes the $2n^{th}$ derivative of the
integrand, and will largely over estimate the error
\cite[sec. 3.1.1]{AfKlinteberg2016quad}.

The above integral can be written in the following form (with $a=0$, $p=1$)), 
\begin{equation}
\int_{-1}^{1} \frac{1}{\left( (t-a)^2+b^2 \right)^p} \, dt = \int_{-1}^{1}
\frac{1}{\left(t-z_0 \right)^p \left(t-\bar{z}_0 \right)^p} \, dt, \quad z_0=a+ib, \quad a \in
\reals, b > 0, p\in \zplus.
\label{eq:cartesian_int_basic}
\end{equation}
Here, it is clear that with $b$ small, the integrand has poles in the
complex plane close to the integration interval along the real axis.
Donaldson and Elliott \cite{Donaldson1972} introduced a theory that
defines the quadrature error as a contour integral in the complex
plane over the integrand multiplied with a so-called remainder
function, that depends on the quadrature rule.  Using residue calculus
for this meromorphic integrand, Elliott et al. \cite{Elliott2008}
derived an error estimate for the error in the approximation of
(\ref{eq:cartesian_int_basic}) for $p=1$ with an $n$-point
Gauss-Legendre quadrature rule. Later, af Klinteberg and Tornberg
\cite{AfKlinteberg2016quad} derived error estimates for a general
positive integer $p$. In \cite{AfKlinteberg2016quad}, results were also
derived for the trapezoidal rule (hence with a different remainder
function), considering integration over the unit circle rather than an
open interval.

{Error estimates were also derived for the $n$-point Gauss-Legendre
quadrature rule by af Klinteberg and Tornberg in 
\cite{AfKlinteberg2016quad} for integrals in the form}
\begin{equation}
\int_{-1}^{1} \frac{1}{\left(t-z_0 \right)^p} \, dt, \quad z_0=a+ib, \quad a \in
\reals, b \ne 0, p\in \zplus.
\label{eq:complex_int}
\end{equation}
%ZZZZZZZZZZ
%Layer potentials in 2D can be rewritten using complex variables. Error
%estimates for the complex form of the harmonic double
%layer potential were derived in \cite{AfKlinteberg2016quad}, and more accurately
%generalized to panel-based Gauss Legendre quadrature for general smooth
%curves $\Gamma$ in \cite{AfKlinteberg2018}. Here we will provide error estimates on a
%form convenient to use for the form in (\ref{eq:layerpot_S}), with $p$
%integer. 
%ZZZZZZZZZ

In 2D, it is convenient to rewrite layer potentials using complex
variables. A typical form of integrals to evaluate over one segment of
the curve is
\begin{equation}
\int_{-1}^{1} \frac{f(t) \, \gamma'(t)}{(\gamma(t)-z_0)^p} \, dt, 
\label{eq:complex_int_curved}
\end{equation}
where $\gamma(t) \in \C$ is a parameterization of the curve segment,
which is assumed to be smooth.  In \cite{AfKlinteberg2018}, af Klinteberg and
Tornberg derived the quadrature error estimate for the Gauss-Legendre
method as applied to (\ref{eq:complex_int_curved}) for $p \in
\zplus$.
{Generalization of the error estimates for ``flat panels''
(\ref{eq:complex_int}) to curved segments
(\ref{eq:complex_int_curved}) introduces a geometry factor. Evaluation
of the estimates requires the knowledge of $t_0 \in \C$ such that
$\gamma(t_0)=z_0$. A numerical procedure is introduced to compute
$t_0$, given the Gauss-Legendre points used to discretize the
panel.}

{The harmonic double layer potential can be written in the form
(\ref{eq:complex_int_curved}) with $p=1$. In the combined field
Helmholtz potential, Hankel functions are present, and additional
steps are needed in the derivation of the error estimates. In both
cases, the resulting estimates are remarkably
accurate, when combined with the numerical
procedure to compute $t_0$ \cite{AfKlinteberg2018}. Using the same techniques, error
estimates for Stokes layer potentials are derived in
\cite{Palsson2019}, again with excellent results.}

{Evaluations of integrals in the form (\ref{eq:complex_int_curved}) with larger $p$ are required to
obtain expansion coefficients in the ``Quadrature by Expansion''
method (QBX) \cite{klockner2013}. 
In \cite{AfKlinteberg2018}, the derived error estimates were used to control the
coefficient error in the expansions, in the framework of an adaptive QBX
method applied to evaluate the harmonic double layer potential and the
combined field Helmholtz potential. This way, automatic parameter
selection in order to fulfil a desired accuracy was facilitated. }

%that was developed in \cite{AfKlinteberg2018}, to facilitate
%automatic parameter selection in order to fulfil a desired accuracy.

%%The resulting estimates become remarkably accurate.

%The derivations for the harmonic double layer potential and for the
%combined field Helmholtz potential can be found in
%\cite{AfKlinteberg2018}, and the estimates for Stokes equations
%derived using the same techniques can be found in \cite{Palsson2019}.

%The derivations for Laplace and Helmholtz equations can be found in
%\cite{AfKlinteberg2018}, and the estimates for Stokes equations
%derived using the same techniques can be found in \cite{Palsson2019}.

{The integral in (\ref{eq:cartesian_int_basic}) (or (\ref{eq:complex_int}) if using a
formulation in complex variables) is the simplest
prototype integral that can be related to an integral over a segment
of a curve. Considering instead a patch of a 3D surface, the simplest 
%Increasing the dimension, the simplest
two-dimensional integral to
consider is }
\begin{equation}
\int_{-1}^{1} \int_{-1}^{1} \frac{ds \, dt}{\left( (t-a)^2+
    s^2+b^2 \right)^p} 
\quad a \in \reals, b > 0.
\label{eq:int2d_basic}
\end{equation}
Elliott et al. \cite{Elliott2011} studied the approximation of this integral with an $n_s
\times n_t$ tensor product Gauss-Legendre rule. They derived error
estimates for the two cases with $p=1$ and $p=1/2$, the latter for $a=0$. 
This error analysis was further extended in  \cite{Elliott2015}, including a
higher order error term that was previously neglected. 
Considering such tensor product Gauss-Legendre rules, 
af Klinteberg and Tornberg \cite{AfKlinteberg2016quad}, derived
quadrature error estimates for QBX coefficients evaluated over flat 2D
patches. A similar estimate was also derived for the case of a spheroidal
surface, discretized with a tensor product rule with the
trapezoidal rule in the (periodic) azimuthal angle and the Gauss-Legendre rule in
the polar angle. 

In \cite{Morse2020}, Morse et al. derives an error estimate for the evaluation of
the double layer potential over a general surface in 3D discretized with
quadrilateral patches and a tensor product Clenshaw-Curtis quadrature rule. The
estimate however contains high derivatives of the Green's function which makes
it difficult to evaluate and hence, as the author acknowledges, cannot
really be applied. In this paper we aim to provide error estimates
without unknown coefficients that can rapidly be evaluated. 
They can then be directly applied to determine e.g. when a regular
quadrature rule is insufficient or how large upsampling is needed.
%%and therefore of direct practical use. XXX

%In 2D, major developments were needed to go from estimating the
%quadrature error from one flat panel with unit density, which is one
%way to interpret the integral in (\ref{eq:complex_int}), to actually
%estimating the quadrature errors of layer potentials over curved
%segments for the Gauss-Legendre quadrature rule.
%In this paper, we start by extending these results and evaluation
%procedure to the global trapezoidal rule. 
%Major developments both of theory and procedure are then required to advance from
%studying (\ref{eq:int2d_basic}) to layer potentials of the form
%(\ref{eq:layerpot_S}) over smooth curved surfaces in $\reals^3$. To do so, we will
%first write the layer potentials with real valued kernels, and
%consider the integral over curved segments in both $\reals^2$ and
%$\reals^3$. This will constitute the generalization of the base
%integral in (\ref{eq:cartesian_int_basic}) to generic layer potentials
%over curved segments, with error estimates both for the global
%trapezoidal rule and composite Gauss-Legendre quadrature. These 
%results can then directly be used for estimation of quadrature errors for 2D
%layer potentials. The results for curves in $\reals^3$ will be used in the
%following derivation of quadrature error estimates for the evaluation
%of layer potentials over surfaces in $\reals^3$. 

\section{Contributions and outline}
\label{sec:outline}

In this paper, we derive estimates for the numerical errors that
results when applying quadrature rules to nearly singular
integrals. Specifically, we consider
Gauss-Legendre (panel based) and trapezoidal (global) approximations
for evaluation of integrals of the type
\begin{align}
  u(\v x) = \int_{\Surf} \frac{k \pars{\v x, \v y}\sigma(\v y)
  }{\norm{\v y - \v x}^{2p}} \dif S(\v y) ,
\label{eq:layerpot_S2}
\end{align}
where $\v x$ can be close to, but not on, $S$.
We assume the functions $k \pars{\v x, \v y}$ and $\sigma(\v y)$ as
well as $S$ to be smooth and derive error estimates for $2p \in \zplus$
for the following cases:
\begin{enumerate}
\item $S$ is a curve in $\reals^2$ or
  $\reals^3$ that we denote $\Gamma$, where $\Gamma=\gamma(E)$, $E
  \subset \reals$. 

In this case, we can write the layer potential in (\ref{eq:layerpot_S2}) in the equivalent form
\begin{align}
  u(\v x) = \int_E \frac{k \pars{\v x, \v \gamma(t)}\sigma(\v\gamma(t)) }{\norm{\v\gamma(t)-\v x}^{2p}} \norm{\v\gamma'(t)} \dif t
  = \int_E \frac{f(t) \dif t}{\norm{\v\gamma(t)-\v x}^{2p}},
\label{eq:base_integral}
\end{align}
where we in the last step have collected all the components that are
assumed to be smooth into the function $f(t)$, which has an implicit
dependence on $\v x$.

\item $S$ is a two-dimensional surface in $\reals^3$, parameterized by
  $\v \gamma: E \rightarrow \reals^3$, $E=\left\{ E_1 \times E_2 \right\}
  \subset \reals^2$.
Now, the prototype layer potential (\ref{eq:layerpot_S2}) takes the form
\begin{align}
  u(\v x) = \iint_E 
  \frac{k \pars{\v x, \v \gamma(s,t)}\sigma(\v\gamma(s,t)) }
  {\norm{\v\gamma(s,t)-\v x}^{2p}} 
  \norm{\dpd{\v\gamma}{s} \times \dpd{\v\gamma}{t}}
  \dif t \dif s
  = \iint_E \frac{f(s,t) \dif t \dif s}{\norm{\v\gamma(s,t)-\v x}^{2p}}.
\label{eq:3Dpot}
\end{align}
Here we have again collected all the smooth components into the function $f(s,t)$, which depends implicitly on $\v x$.

\end{enumerate}
{It is possible to derive error estimates also for kernels with a logarithmic singularity \cite{Elliott2008}, relevant for 2D problems. However we are interested in error estimates for layer potentials in 3D, where it is not so natural to consider logarithmic singularities, so we will limit this study to case 1 and 2 above.
}

Considering the approximation of (\ref{eq:base_integral}), the
trapezoidal rule will always be applied to a closed curve, for which
it is spectrally accurate for smooth integrands. For the
Gauss-Legendre quadrature rule, we will consider the discretization of
one open segment of the curve. Any curve can be divided into several
such segments (panels), and the total quadrature error can be obtained
by adding the contributions from all panels.
%The trapezoidal rule is a global quadrature rule that for a curve will be applied
%to the full closed curve. The Gauss-Legendre quadrature rule on the
%other hand is panel based, and it is hence natural to divide the curve
%into segments and consider one segment at a time. 
Similarly, for surfaces in $\reals^3$, we will assume that we have a
parameterization for either the global surface or a quadrilateral
panel of the surface, and will apply a tensor product quadrature rule
based on the trapezoidal rule and Gauss-Legendre quadrature,
respectively.

Considering layer potentials and derivatives thereof, $p$ will
naturally be a positive integer for curves in the plane, and a
positive half-integer for surfaces in $\reals^3$, even if we do not
have to restrict the estimates to these cases. Our strategy in
deriving the error estimates for surfaces will be based on deriving an
estimate of the error in one direction first, then integrating it in
the other, as will be discussed in section
\ref{sec:surface_estimates}.  For this, we need error estimates for
the numerical integration over curves in $\reals^3$, which is included
in case 1 above. 

The outline of this paper is as follows: In section
\ref{sec:basictheory}, we briefly introduce the theory of Donaldson
and Elliott \cite{Donaldson1972},  to estimate quadrature errors using contour
integrals in the complex plane. When the integration is over a planar
curve, layer potentials are conveniently rewritten in complex form. In
\cref{sec:CurvesComplex}, we summarize earlier results for the Gauss-Legendre
quadrature rule for such complex values kernels as based on 
\cite{AfKlinteberg2016quad} and \cite{AfKlinteberg2018}.
Using the underlying derivations, it is straightforward to
arrive at corresponding results for the trapezoidal rule. {Since such
results have not been previously available, we
derive them in this section.}  Error estimates for solving the Laplace equation in 2D
are compared to actual errors for both discretizations, displaying a
remarkable precision.  

\Cref{sec:qerr1Dcurves} treats case 1 above, with the main theoretical
results in \cref{sec:trapezoidal-rule,sec:gauss-legendre-rule}. It is
also discussed how to evaluate these estimates in practice and the
excellent predictive accuracy of the estimates is illustrated with
numerical examples.

In \cref{sec:surface_estimates}, the results from
\cref{sec:qerr1Dcurves} for integration over curves in $\reals^3$ are
used when we extend the analysis to surfaces (case 2).  Also here,
numerical results are used to compare our final estimates to actual
measured errors. It is shown that while the estimates are not as
precise as those for integration over curves, they still have a good
predictive power and can be used to determine at which point the
accuracy of the regular quadrature becomes unable to meet a specified
error tolerance.

\section{Formulas for quadrature errors using complex analysis}
\label{sec:basictheory}

Let us introduce the base interval $E$, which for the Gauss-Legendre
quadrature will be  $[-1, 1]$ and for the trapezodial rule $[0, 2\pi]$.
Consider an integral over such a base interval
\begin{equation}
I[g]=\int_{E} g(t) dt.
\end{equation}
Applying an $n$-point quadrature rule to approximate this definite integral, we get
\begin{equation}
Q_n[g]=\sum_{\ell=1}^n g(t_{\ell}) w_{\ell},
\end{equation}
where the quadrature nodes $t_{\ell}$ and corresponding weights
$w_{\ell}$ depend on the quadrature rule. 

The decay of the error 
\begin{equation}
E_n[g]=I[g]-Q_n[g]
\end{equation}
as a function of $n$ will depend on the function $g$. Classical error
bounds for the Gauss-Legendre quadrature rule involve higher
derivatives of $g$ with increasing $n$.  However, such error bounds do
not work well for integrals such as (\ref{eq:base_integral}) when the
evaluation point $\v x$ is close to $\Gamma$. See the discussion in
\cite{AfKlinteberg2016quad}, that illustrates how a classical error
estimate will for some cases even predict a growth in error with increasing $n$ when the
actual error decays with $n$.

We can consider the integral over $E$ as an integral over a part of
the real line in the complex plane. When an integral is nearly
singular, that means that the complex continuation of the integrand
will have a very small region around $E$ where it is analytic. 
In these cases, much better error estimates can be achieved by using
the theory of Donaldson and Elliott \cite{Donaldson1972}, based on contour integrals in the
complex plane. 

Following their lead, we can write
\begin{equation}
Q_n[g]=\sum_{\ell=1}^n g(x_{\ell}) w_{\ell} = \frac{1}{2\pi i}
\int_{C} g(z) q_n(z) dz
\end{equation}
where $C$ contains the integration interval $E$ and where the complex
continuation of $g$ is analytic on and inside $C$. 
The integration interval is $[-1, 1]$ for Gauss-Legendre, and is simple to enclose.  
For the trapezoidal rule, the contour can be chosen as the rectangle
$[0,2\pi] \pm ia, a>0$. The sides of the rectangle cancel, leaving only the top and
bottom lines.
The function $q_n(z)$ is specific to each quadrature rule, as will be
discussed further below. 

We can furthermore write
\begin{equation}
g(t)= \frac{1}{2\pi i} \int_{C} \frac{g(z)}{z-t} dz,
\end{equation}
and hence
\begin{equation}
I[g]=\int_{C} g(z) m(z) dz,
\end{equation}
where 
\begin{equation}
m(z)=\int_{E} \frac{dt}{z-t} .
\end{equation}

From this, we can define $k_n(z)=m(z)-q_n(z)$ such that:
\begin{equation}
E_n[g] =I[g]-Q_n[g]= \frac{1}{2\pi i} \int_{C} g(z) k_n(z) dz,
\label{eq:En_w_kn}
\end{equation}
where $k_n(z)$ depends on the quadrature rule. 

There is no closed expression for $k_n(z)$ for the $n$ point
Gauss-Legendre rule. In the limit as $n \rightarrow \infty$ it can
however be shown to satisfy  \cite{Donaldson1972}, 
\begin{equation}
k_n(z) \simeq \frac{c_n}{(z+\sqrt{z^2-1})^{2n+1}},
\label{eq:knGL}
\end{equation}
where the constant 
\begin{equation}
c_n = \frac{2 \pi (\Gamma(n+1))^2 }{ \Gamma(n+1/2) \Gamma(n+3/2) }   
\simeq 2 \pi,
\end{equation}
with $\Gamma(.)$ the gamma function. Note that we have used $\Gamma$ without an
argument to denote a curve in $\reals^2$ and $\reals^3$. This should
be clear from the context such that it causes no confusion.  
{In eq. (\ref{eq:knGL})} and for the remainder of this paper, $\sqrt{z^2-1}$ is defined as 
$\sqrt{z+1}\sqrt{z-1}$ with $-\pi < \arg(z \pm 1) \le \pi$ \cite{Elliott2008}.

For the trapezoidal rule with $n$ points, we have \cite{Trefethen2014}
\begin{equation}
k_n(z) = 2 \pi i \left\{
\begin{array}{cc}
\frac{-1}{e^{-inz}-1} & \Im z > 0, \\
\frac{1}{e^{inz}-1} & \Im z < 0.
\end{array}
\right.
\label{eq:knTz}
\end{equation}

This theory by Donaldson and Elliot has been used on integrals of the
simple form (\ref{eq:cartesian_int_basic}). In Elliot et
al. \cite{Elliott2008}, integrals with $p=1$, but including also a
nominator $t^k$, $k=0,1,\ldots$ were considered.  In
\cite{AfKlinteberg2016quad}, the current authors studied integrals with positive integer
values of $p$. The results for the errors in approximation with an
$n$-point Gauss Legendre rule for $k=0$ in \cite{Elliott2008} and
$p=1$ in \cite{AfKlinteberg2016quad} coincide. 
In both these works, the fact that the integrand is meromorphic was
used as the estimates were derived starting from the contour integral
in (\ref{eq:En_w_kn}).

Elliott et al  \cite{Elliott2008} also studied the same integral with
a non-integer $p$, $0<p<1$. In this case, the integrand is no longer
meromorphic and one needs to work with branch cuts when analyzing the
contour integral in (\ref{eq:En_w_kn}).

When considering estimates for integrals with positive integer
values of $p$, $p>1$,  derivatives of $k_n$ are needed. 
For the Gauss-Legendre $n$-point quadrature we can estimate
\cite{AfKlinteberg2016quad}
\begin{align}
  k_n^{(q)}(z) \simeq
  \pars{-\frac{2n+1}{\sqrt{z^2-1}}}^{q}
  \frac{2\pi}{(z + \sqrt{z^2-1})^{2n+1}}
  .
  \label{eq:knq_gl_asy}
\end{align}
For the trapezoidal rule, an asymptotic form of (\ref{eq:knTz})
for $n\to\infty$ was derived \cite{AfKlinteberg2016quad}
\begin{align}
  k_n(z) \simeq 2\pi i
  \begin{cases}
    -e^{inz} & \Im z > 0,  \\
    e^{-inz} & \Im z < 0,
  \end{cases}
  \label{eq:kn_trapz_per_asy}
\end{align}
with derivatives
\begin{align}
  \dod[q]{k_n(z)}{z} \simeq 2\pi i
  \begin{cases}
    -(in)^q e^{inz} & \Im z > 0,  \\
    (-in)^q e^{-inz} & \Im z < 0 .
  \end{cases}
  \label{eq:knq_trapz_per_asy}
\end{align}
We are in most cases only interested in the magnitude of the error, in
which case it is useful to write
\begin{align}
  \abs{ k_n^{(q)}(z) } \approx \frac{2\pi n^{q}}{e^{n|\Im z|}}  .
  \label{eq:knp_trapz_abs}
\end{align}
\Cref{eq:kn_trapz_per_asy,eq:knq_trapz_per_asy,eq:knp_trapz_abs} are
good approximations to \eqref{eq:knTz} as long as
$e^{n|\Im z|} \gg 1$. As we will see, this is also a requirement for the error to be
less than $\mathcal O(1)$, so they are useful in most practical
applications.

\section{Quadrature errors near planar curves with kernels in complex form}
\label{sec:CurvesComplex}

The first aim of this paper is to derive error estimates for the
Gauss-Legendre quadrature and trapezoidal rule as applied to
(\ref{eq:base_integral}). We need to do so in a form such that they are
applicable both for planar and spatial curves $\gamma(t)$.

Before doing so, we however want to discuss some closely related
results for planar curves. In this case, layer potentials can be
rewritten in complex form, see e.g. \cref{app:DoubleLayerEst}
for the harmonic double layer potential.
We will summarize the results obtained in
\cite{AfKlinteberg2018} for layer potentials in 2D in
complex form. These error estimates are for the Gauss-Legendre
quadrature rule, and for completeness we also derive the corresponding
error estimates for the trapezoidal rule. 

For layer potentials, and derivatives thereof, expressed in complex
variables, the generic form can in analogy with
(\ref{eq:base_integral}) be written as
\begin{equation}
I[\psi_p](z_0)=\int_{E} \psi_p(t,z_0) \, dt, 
\label{eq:complex_int_psip}
\end{equation}
where 
\begin{equation}
\psi_p(t,z_0)=\frac{g(t) \gamma'(t)
}{(\gamma(t)-z_0)^p}=\frac{f(t)}{(\gamma(t)-z_0)^p},
\label{eq:complex_psip}
\end{equation}
with $\gamma(t) \in \C$  a parameterization of the curve
(segment), $z_0 \in \C$ the evaluation point and $p$ a positive
integer. Note that in what follows, we will work with the first form in
(\ref{eq:complex_psip}). 

Using that $\psi_p(t,z_0)$ is a meromorphic function, with a pole at $z_0$
of order $p$, the following estimate can be derived for the quadrature
error $E_n[\psi_p](z_0)$ as defined in  (\ref{eq:En_w_kn}) \cite{AfKlinteberg2018},
\begin{equation}
%|E_n[\psi_p](z_0)| \approx \frac{1}{(p-1)!} \frac{|g(t_0)|}{ |\gamma'(t_0)|^{p-1}} |k_n^{(p-1)}(t_0)|
|E_n[\psi_p](z_0)| \approx \frac{1}{(p-1)!} \abs{ \frac{g(t_0)}{ (\gamma'(t_0))^{p-1}} k_n^{(p-1)}(t_0)}
\label{eqn:En-complex-w-kn}
\end{equation} 
where 
%\begin{equation}
%G(t_0)=\frac{1}{\gamma'(t_0)}
%\end{equation}
%and 
{$t_0$ is the point in $\C$ closest to $E$ such that $\gamma(t_0)=z_0$.} Hence, for an integral
with $p=1$, $\gamma'(t)$ does not appear in the estimate. 
For the Gauss-Legendre $n$-point quadrature we use the estimate (\ref{eq:knq_gl_asy})
for $k_n^{(p-1)}(t_0)$ to obtain the following result.
%cite{AfKlinteberg2016quad}
%\begin{align}
%  k_n^{(q)}(z) \simeq
%  \pars{-\frac{2n+1}{\sqrt{z^2-1}}}^{q}
%  \frac{2\pi}{(z + \sqrt{z^2-1})^{2n+1}}
%  .
%  \label{eq:knq_asy0}
%\end{align}
%Combining this, we have the following result.
\begin{estimate}
The error in approximating the integral (\ref{eq:complex_int_psip})
with the $n$-point Gauss-Legendre quadrature rule can in the limit $n
\rightarrow \infty$ be estimated by
\begin{equation}
%|E_n[\psi_p](z_0)| \approx \frac{2\pi}{(p-1)!} 
 % \left|\frac{2n+1}{\sqrt{t_0^2-1}}\right|^{p-1}
%\frac{|g(t_0)|}{  |\gamma'(t_0)|^{p-1}}
%    \frac{1}{\rho(t_0)^{2n+1}}.
|E_n[\psi_p](z_0)| \approx \frac{2\pi}{(p-1)!} 
\abs{  \left(\frac{2n+1}{\sqrt{t_0^2-1}}\right)^{p-1}
\frac{g(t_0)}{  (\gamma'(t_0))^{p-1}}    \frac{1}{\rho(t_0)^{2n+1}} }.
%%  \frac{1}{|t_0 + \sqrt{t_0^2-1}|^{2n+1}}
\label{eq:errest-complexGL}
\end{equation} 
Here, $p$ is a positive integer, $\gamma(t) \in \C$  is a
parameterization of the curve,  {$t_0 $ is the point in $\C$ closest to $E$ such that $\gamma(t_0)=z_0$} and $\rho(t) = t +
\sqrt{t+1}\sqrt{t-1}$.
\label{est-compl-GL}
\end{estimate}
 This result can be found in equation (68) of \cite{AfKlinteberg2018}
if we adjust for the definition in equation (48) in
\cite{AfKlinteberg2018} as compared to (\ref{eq:complex_int_psip})
above. The above result is a generalization of \cite{AfKlinteberg2016quad} Thm. 1 to
curved panels. 

This result is an asymptotic result for $n \rightarrow \infty$, but
it is remarkably accurate for $n$ point Gauss-Legendre
quadrature already for moderate values of $n$. 
A larger $n$ is however needed for larger values of $p$. A rule of
thumb from \cite{AfKlinteberg2018} is that we need $n>2p$ to have a
good precision in the estimates.

This estimate can be used in practice, since there are no unknown
coefficients. Given a $z_0$, one does however need $t_0$ such that
$\gamma(t_0)=z_0$. We denote $t_0$ the pre-image of $z_0$. 
In \cite{AfKlinteberg2018} , a numerical procedure
was used to determine an accurate approximation of $t_0$, as 
discussed in the next section.  The above result is an estimate, and not
a bound. A bound on the error was given in \cite{AfKlinteberg2018},
which is theoretically of value, but it does overestimate the error by
a large factor.

The corresponding error estimate for the trapezoidal rule has not been
derived before, but is straightforward to do with all the components
that we have available. The error estimate
(\ref{eqn:En-complex-w-kn}) introduced above still holds. To derive it
one needs to use a different contour to enclose the integration
interval, as commented on in \cref{sec:basictheory}. Compare also to the
forthcoming discussion in \cref{sec:qerr1Dcurves}.
Combining (\ref{eqn:En-complex-w-kn}) with the estimate
(\ref{eq:knp_trapz_abs}) for $k_n^{(p-1)}(t_0)$ for the trapezoidal rule, 
we obtain the following result.
\begin{estimate}
The error in approximating the integral (\ref{eq:complex_int_psip})
with the $n$-point trapezoidal rule can in the limit $n
\rightarrow \infty$ be estimated by 
\begin{equation}
%|E_n[\psi_p](z_0)| \approx 2\pi \frac{ n^{p-1}}{(p-1)!} \frac{|g(t_0)|}{
%  |\gamma'(t_0)|^{p-1}} e^{-n |\Im t_0|}
|E_n[\psi_p](z_0)| \approx 2\pi \frac{ n^{p-1}}{(p-1)!} \abs{ \frac{g(t_0)}{
  (\gamma'(t_0))^{p-1}} } e^{-n |\Im t_0|}.
\label{eq:errest-complexTZ}
\end{equation}
Here, $p$ is a positive integer, $\gamma(t) \in \C$  is a
parameterization of the curve and  {$t_0 $ is the point in $\C$ closest to $E$ such that $\gamma(t_0)=z_0$}.
\label{est-compl-TZ}
\end{estimate}
%\ref{est-compl-GL}
%\ref{est-compl-TZ}

In \cite{Barnett2014}, Barnett studied the error in the evaluation of
the harmonic double layer potential with the trapezoidal rule. This
would correspond to $p=1$ and a specific choice of $f$ in the estimate
above. He proved that there exist constants $C$ and $n_0$ such that
the error is bounded by $C e^{-n |\Im t_0|}$ for all $n>n_0$ (Theorem
2.3 in \cite{Barnett2014}).

Hence, from both our estimate and this bound, we see the exponential
decay of the error with $n$, but also that it is the distance of the
pre-image $t_0$ to the real line that determines the decay rate. 
%\subsection{Root finding in the complex plane}
%\label{sec:mapping}
{
\begin{remark}
Note that given a $z_0 \in \C$, there is in general more than one $t_0 \in \C$ such that
$\gamma(t_0)=z_0$. In our error estimates, we only include the contribution
from the $t_0$ closest to $E$. This is motivated by the fact that the
error decays rapidly with the distance from $E$.
\end{remark}}

\subsection{Examples for planar curves with kernels in
  complex form}
\label{sec:results_complex}

In order to evaluate the estimates given in {\em Error estimate
  \ref{est-compl-GL}} and {\em Error estimate \ref{est-compl-TZ}}
above, we need to know $t_0$ to be able to evaluate $\gamma'(t_0)$ and
$g(t_0)$.  To obtain the pre-image $t_0$, we need to solve
$\gamma(t_0)=z_0$.  We however frequently do not have analytical
expressions neither for $\gamma$ nor $g$.

%The function $g$ typically contains a layer density obtained by
%numerically solving an integral equation, and $\gamma$ might also only
%be discretely defined.

In \cite{AfKlinteberg2018}, a polynomial $\poly_n[\gamma](t) \in \C$ of
degree $n-1$ is defined as an approximation of $\gamma(t)$ with the
Legendre polynomials as a basis.  Using the polynomial
$\poly_n[\gamma](t)$ as an approximation of the analytic continuation
of $\gamma(t)$, we can now find an accurate approximation of $t_0$ by solving
$\poly_n[\gamma](t_0)=z_0$.  This can be done efficiently and robustly
using Newton's method. For details regarding this procedure and the related
evaluation of $\gamma'(t_0)$ and $g(t_0)$, see the discussion in
\cite{AfKlinteberg2018}.

For the Gauss-Legendre rule, $n$ is the number of points on one panel,
i.e. along one segment on the curve, where as for the trapezoidal
rule, $n$ is the number of points used to discretize the full
curve. Using a global approximation of $\gamma(t)$ based on
e.g. trigonometric polynomials hence adds an unnecessarily large extra
cost. Here, we instead use a local $5$th order Taylor expansion to
approximate the curve in the root finding process. This will in section
\ref{sec:curve-rootfinding} be discussed in a
more general setting for root finding that can be used both in $\reals^2$ and $\reals^3$.

\begin{figure}[htbp]
  \centering
  \begin{minipage}{0.45\textwidth}  
    \includegraphics[width=\textwidth]{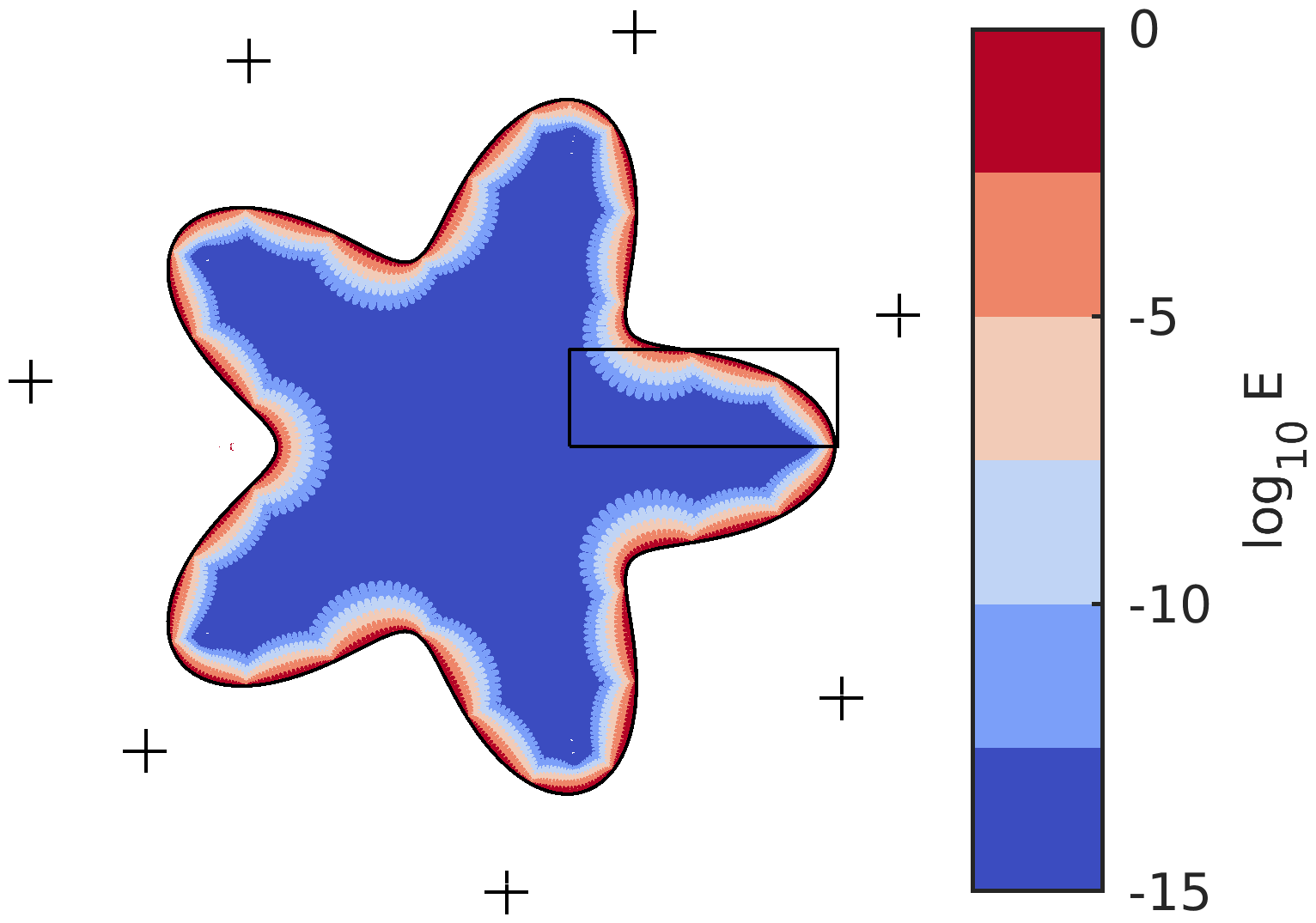}
  \end{minipage}
  \begin{minipage}{0.54\textwidth}
    \includegraphics[width=\textwidth]{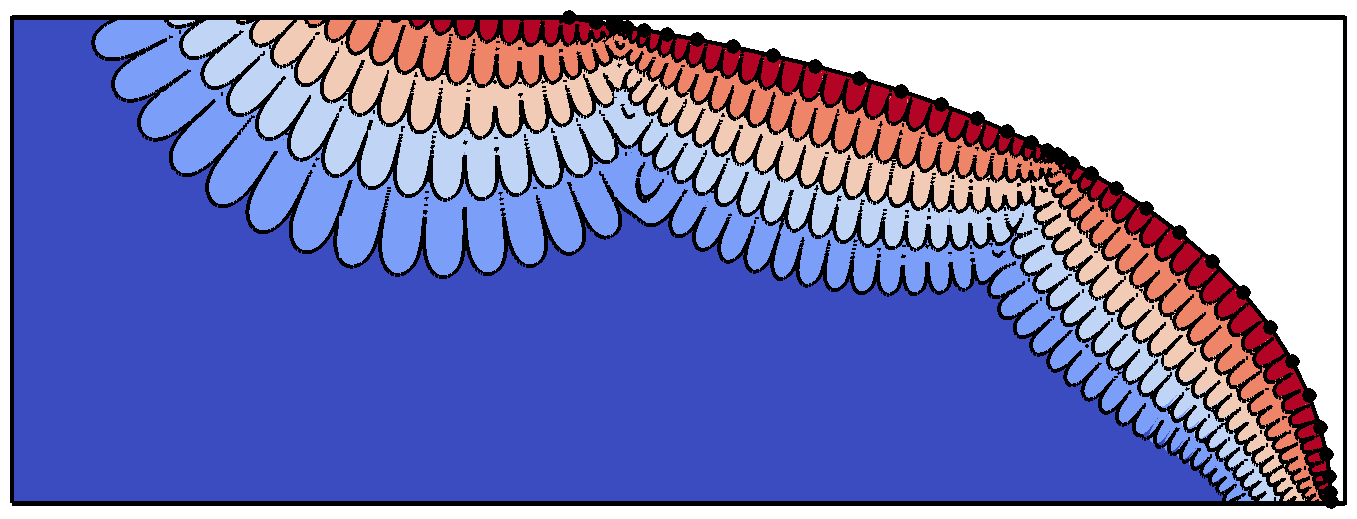} \\
    \includegraphics[width=\textwidth]{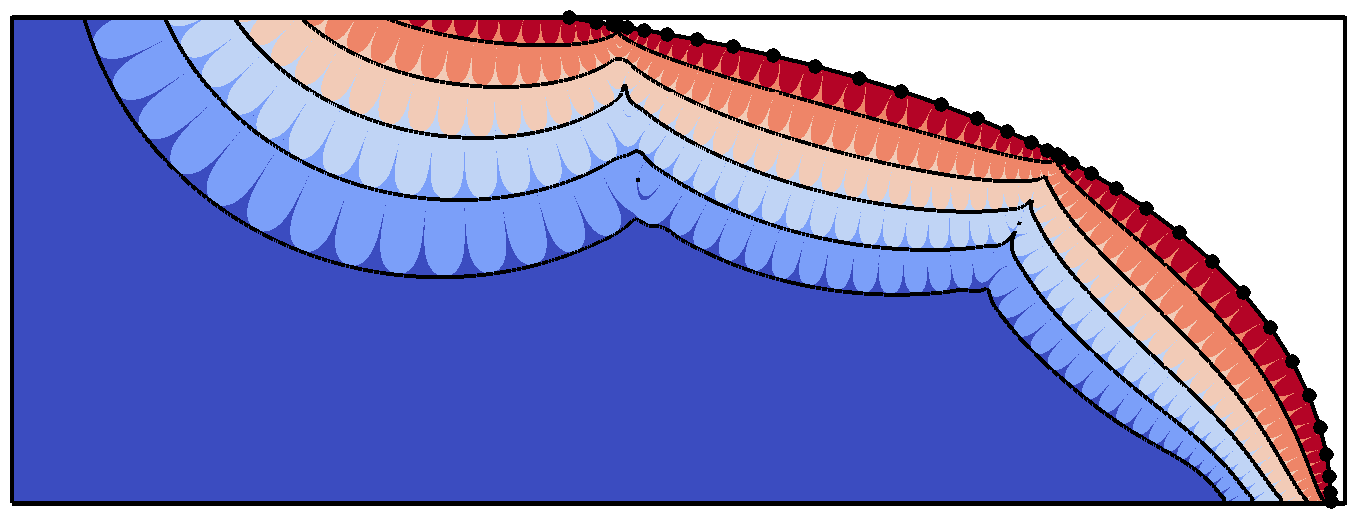}    
  \end{minipage}
  \caption{Figure~3 from \cite{AfKlinteberg2018}. The field shows the
    error when using a panel based 16-point Gauss-Legendre quadrature for
    evaluating the 2D Laplace double layer potential, which represents
    the solution to the interior boundary value problem constructed by
    using the field from the point sources marked with plus ($+$)
    signs as boundary data. 
The error estimates are plotted with black contours; the top right
plot displays estimate (\ref{eq:est-HDL-GL-im}) and the bottom right
plot estimate (\ref{eq:est-HDL-GL}).
    \\{\scriptsize Figure originally published in: L. af Klinteberg
      and A.-K. Tornberg. Adaptive Quadrature by Expansion for Layer
      Potential Evaluation in Two Dimensions. SIAM J. Sci. Comput.,
      40(3):A1225–A1249, 2018. Copyright \copyright 2018 Society for
      Industrial and Applied Mathematics.  Reprinted with permission.
      All rights reserved.}}
  \label{fig:aqbx_plot}
\end{figure}

We now present some numerical results for the harmonic double layer
potential. We solve the interior Dirichlet Laplace problem on a
starfish shaped domain depicted in figure \ref{fig:aqbx_plot}. The
boundary data is taken from the field obtained from point sources
whose locations are marked in the same picture. The solution is
obtained in two steps. First, we solve an integral equation to obtain a
layer density $\sigma$, defined on the boundary of the domain. Then,
at any point in the domain where we want to compute the solution, we
evaluate the harmonic double layer potential as given in 
(\ref{eq:doublelayerC}) in \cref{app:DoubleLayerEst}.
Since we know the exact solution by construction, we can measure the
pointwise numerical error. With this, we can compare our estimate
of the error with the actual error. 

As discussed in \cref{app:DoubleLayerEst}, the error estimate
(\ref{eqn:En-complex-w-kn}) for the kernel (\ref{eq:doublelayerC})
becomes simply $E_n\approx |\sigma(t_0) k_n(t_0)|$, if we ignore
taking the imaginary part that is in the kernel (we have $p=1$).  If we do include the
imaginary part, we get instead $E_n\approx |\Im \{\sigma(t_0) k_n(t_0)
\}|$.
Estimate (\ref{eq:errest-complexGL}) in {\em Error estimate
  \ref{est-compl-GL}} for the Gauss-Legendre rule and
(\ref{eq:errest-complexTZ}) in {\em Error estimate \ref{est-compl-TZ}}
for the trapezoidal rule, and the variants of taking the imaginary
parts hence yield four estimates
\begin{align}
E_n^{GL} & \approx 2 \pi |\sigma(t_0)/ \rho(t_0)^{2n+1}|, 
\label{eq:est-HDL-GL} \\
E_n^{GL,Im} & \approx 2 \pi |\Im \{\sigma(t_0)/ \rho(t_0)^{2n+1} \}|, 
\label{eq:est-HDL-GL-im} \\
E_n^{TZ} & \approx 2 \pi |\sigma(t_0) \,  e^{-n |\Im t_o|},
\label{eq:est-HDL-TZ} \\
E_n^{TZ,Im} & \approx 2 \pi |\Im \{\sigma(t_0) \, e^{-n |\Im t_o|} \}
              |,
\label{eq:est-HDL-TZ-im} 
\end{align}
where $\rho(t)=t+\sqrt{t+1}\sqrt{t-1}$.

In figure \ref{fig:aqbx_plot}, we present the results for a
Gauss-Legendre discretization. The figure is reproduced with
permission from \cite{AfKlinteberg2018}. The discretization is made
with a 16-point Gauss-Legendre rule using $27$ panels, and all details
can be found in \cite{AfKlinteberg2018}.  The scaling of the layer potential in
\cite{AfKlinteberg2018} removes the factor of $2\pi$ in the error
estimates  (\ref{eq:est-HDL-GL})-(\ref{eq:est-HDL-TZ-im}), but will
also rescale the integral equation such that the layer density
$\sigma$ gets a $2\pi$ factor larger magnitude, so the end result displayed in the figure is the same. 
The color fields in figure \ref{fig:aqbx_plot} show the measured numerical
error, and for comparison, the error estimates ($n=16$) are plotted on top with
black contours in the enlarged plots for a part of the domain. 
To evaluate the error estimates, contributions from the two
panels closest to the evaluation point have been added. 
The error estimates are remarkably accurate, given the simplifications
that have been made. Keeping the imaginary part in the error estimate
(\ref{eq:est-HDL-GL-im}), we can even capture the oscillations of the
error.

\begin{figure}[htbp]
  \centering
  % scripts/twodim_complex_test.m
  \includegraphics[width=0.49\textwidth]{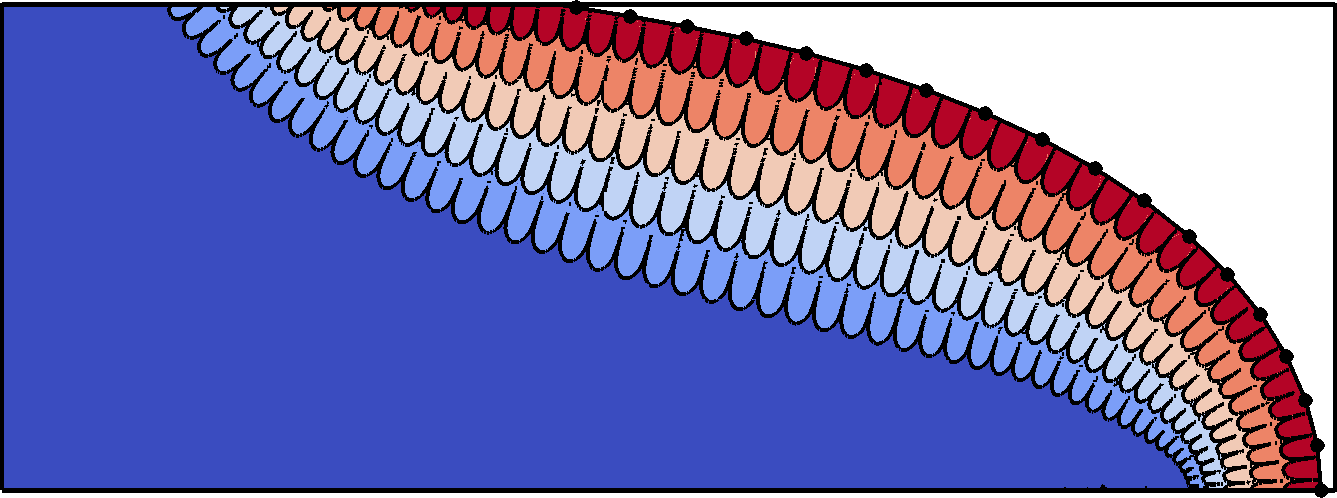}
  \includegraphics[width=0.49\textwidth]{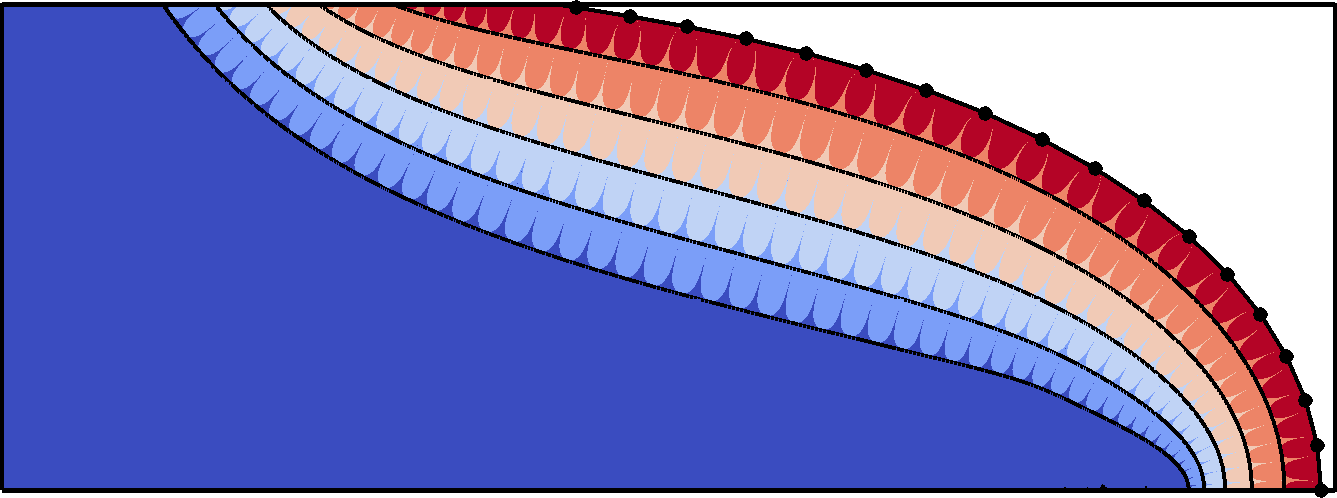}
  \caption{Plots corresponding to the two right plots in
    \cref{fig:aqbx_plot}, but for a discretization with trapezoidal
    rule with $n=250$. Again, the error estimates are plotted with black contours; the left
plot displays estimate (\ref{eq:est-HDL-TZ-im}) and the right
plot estimate (\ref{eq:est-HDL-TZ}).}
  \label{fig:trapz_complex}
\end{figure}

In figure \ref{fig:trapz_complex}, plots corresponding to the two
right plots in figure \ref{fig:aqbx_plot} are shown for a
discretization based on the trapezoidal rule. Here, the full curve is
discretized with $n=250$ points, and the estimates used are
(\ref{eq:est-HDL-TZ-im}) and (\ref{eq:est-HDL-TZ}). 
The error contours
look different for this approximation with uniformly spaced
discretizations points as compared to the panel based Gauss-Legendre
discretization, but again, the precision of the estimates is
excellent. 

Error estimates have also been derived for the Helmholtz and Stokes
equations for $n$-point Gauss-Legendre discretizations, and
derivations and corresponding plots can be found in
\cite{AfKlinteberg2018} and \cite{Palsson2019}, respectively. Using
what has been discussed above, it is straightforward to derive the
corresponding results for the trapezoidal rule. 

\section{Quadrature errors near one-dimensional curves}
\label{sec:qerr1Dcurves}

In this section we will derive error estimates for the numerical
evaluation of the layer potential (\ref{eq:base_integral}). 
The curve $\Gamma$ can be in $\reals^2$ or $\reals^3$, and we will
denote  $\Gamma=\gamma(E)$, $E \subset \reals$ for both cases. 
The form of layer potentials in $\reals^2$ and $\reals^3$ will be such
that $p$ is a positive integer in  $\reals^2$ and a positive
half-integer in  $\reals^3$. In our analysis, we will keep the two
cases of $\reals^2$ and $\reals^3$ together, and will derive error estimates for all $p$ such that $2p \in \zplus$. 

As was commented on in \cref{sec:outline}, we will consider
closed curves for the trapezoidal rule and open curves (segments) with
the Gauss-Legendre quadrature rule. As introduced in section
\ref{sec:basictheory}, the base interval $E$ is set to $[0, 2\pi]$ and
$[-1, 1]$, respectively. 

\subsection{General results}
We now introduce the squared distance function for a curve in $\mathbb
R^d$ ($d=2$ or $3$), given an evaluation point $\v x$, 
\begin{align}
  R^2(t,\v x) \coloneqq \norm{\v\gamma(t)-\v x}^{2} 
  = \sum_{i=1}^d (\gamma_i(t)-x_i)^2,
  \label{eq:R2_def}
\end{align}
such that we can write our integral of interest
(\ref{eq:base_integral}) in the form
\begin{align}
 %% I[\Theta_p](\v x) =
\opintTh = \int_E \Theta_p(t,\v x) \dif t, \qquad
\Theta_p(t,\v x)=\frac{f(t)}{\pars{R^2(t,\v x)}^p}.
\label{eq:cartesian_int}
\end{align}
{Now, if \eqref{eq:base_integral} is
computed using an $n$-point quadrature rule, then the error is given
by the contour integral 
\begin{align}
  \opremTh =
    \frac{1}{2\pi i} \int_{C}  
 \Theta_p(t,\v x) k_n(t)\dif t =
  \frac{1}{2\pi i} \int_{C}  
  \frac{f(t) k_n(t)\dif t}{\pars{R^2(t,\v x)}^p}.
%  =
 % \frac{1}{2\pi i}
 % \int_{C}  \frac{h(t) k_n(t)}
 % {(t-t_0)^p(t-\conj t_0)^p} \dif t,
 \label{eq:base_integral_contour}
\end{align}
Here, the function $k_n(z)$ is specific to the quadrature rule used,
and was given in equations (\ref{eq:knGL}) and (\ref{eq:knTz}),
respectively. $C$ is a contour containing the interval $E$, on and
within which $\Theta_p(t,\v x)$ is analytic. }
The region of analyticity of $\Theta_p(t,\v x)$ is bounded by its
singularities, which, under the assumption that $f$ is smooth, are
given by the roots of the squared distance function $R^2$. Since
$R^2(t,\v x)$ is real for real $t$, the roots will come in complex
conjugate pairs. Let $\braces{t_0, \conj t_0}$ be the pair closest to
$E$, such that
\begin{align}
  R^2(t_0,\v x) = R^2(\conj t_0, \v x) = 0.
\end{align}
We will refer to these points both as roots (of $R^2$) and
singularities (of the integrand). They are in most applications not
known a priori, but can be found numerically for a given target point
$\v x$ (how to do this is discussed in \cref{sec:curve-rootfinding}).

{We can deform the contour $C$  in (\ref{eq:base_integral_contour}) away from $E$, avoiding the
singularities  $t_0$ and $\conj t_0$, see \cref{fig:contourC}. We assume that
the integrand of \eqref{eq:base_integral_contour} vanishes
  faster than $|t|^{-1}$ as $|t|\to\infty$. This means that the
  contributions from those parts of $C$ that are well separated from
  the interval $E$ will tend to zero. If we let the contour tend to
  infinity, deforming it to avoid also other pairs of singularities
  further away from $E$, the error will be given by the sum of contributions
  from all the singularities. Considering the fast decay of the
  contribution from a singularity with the distance from $E$, we further
  assume that we can ignore all roots to $R^2$ except $t_0$ and $\conj t_0$.
For further discussion on multiple roots of $R^2$, see the
discussion in connection to the quadrature method introduced in \cite{AfKlinteberg2020line}.}
  
These assumptions let us approximate \eqref{eq:base_integral_contour}
using only the contributions from the pair of closest singularities, $\braces{t_0, \conj t_0}$. How these
contributions are computed depends on whether or not $p$ is an
integer, as we shall see in the following sections.

{The following derivation will be made for a given evaluation point
$\v x$, and we will temporarily drop the argument $\v x$ and replace $\opremTh \rightarrow \oprem$ for ease of notation.}

\subsubsection{Integer $p$}

If $p$ is integer, then $t_0$ and $\conj t_0$ are simply $p$th order
poles. Starting from (\ref{eq:base_integral_contour}) and letting the
contour $C$ go to infinity, we estimate the
quadrature error using only the residues (as based on the assumptions
made above), 
\begin{align}
  \oprem &\approx 
           -\sum_{w = \{t_0, \conj t_0\}}
           \Res \left[ 
           \frac{ f(t) k_n(t) }
           { \pars{R^2(t)}^p }, w \right] \\
         &= 
           -\sum_{w = \{t_0, \conj t_0\}}
           \frac{1}{(p-1)!} \lim_{t \to w}
           \dod[p-1]{}{t} \pars{ 
           f(t)
           \pars{\frac{t-w}{R^2(t)}}^p
           k_n(t)
           } 
           .
\end{align}
Following \cite{AfKlinteberg2018}, we simplify the derivative in the
above expression by only keeping the term with the highest derivative of
$k_n$. In addition, we define the \emph{geometry factor} $G$, which for a root $w$ of $R^2$ is defined as
\begin{align}
  G(w) = 
  \lim_{t \to w} \frac{t-w}{R^2(t)}
  =
  \left(2\pars{\v\gamma(w)-\v x}\cdot\v\gamma'(w) \right)^{-1}
  . \label{eq:geometry_factor}
\end{align}
This allows us to write 
\begin{align}
  \oprem 
  &\approx
    - \frac{1}{(p-1)!}
    \pars{
    f(t_0)G(t_0)^p k_n^{(p-1)}(t_0)
    +
    f(\conj t_0)G(\conj t_0)^p k_n^{(p-1)}(\conj t_0)
    }
  \label{eq:En_w_conj}\\
  &=
    - \frac{2}{(p-1)!}
    \rebrac{
    f(t_0)G(t_0)^p k_n^{(p-1)}(t_0)
    }
    .    
\label{eq:base_est_int_w_Re}
\end{align}
As the target point $\v x$ moves parallel to the curve, this estimate
oscillates in the same way as the error (see e.g. \cref{fig:field_err}). Capturing these oscillations with an estimate
can in some cases be hard, and is in any case of limited practical
use. We therefore use the triangle inequality on (\ref{eq:En_w_conj}) instead, to get a final, slightly
conservative, estimate for the absolute value of the error,
\begin{align}
  \abs{\oprem} 
  & \approx
   \frac{2}{(p-1)!}  \abs{f(t_0) G(t_0)^p  k_n^{(p-1)}(t_0) }
          = \frac{2} {(p-1)!}  \abs{f(t_0)} \abs{G(t_0)}^p         
           \abs{ k_n^{(p-1)}(t_0) } . 
           \label{eq:base_est_integer}
\end{align}
For a given quadrature rule with corresponding error function $k_n$,
this expression is straightforward to evaluate, and we will do so in
section  \ref{sec:trapezoidal-rule} for the
trapezoidal rule, and in section \ref{sec:gauss-legendre-rule} for the Gauss-Legende rule.

As was mentioned earlier, layer potentials in the plane can be
reformulated using complex variables. In \cref{app:DoubleLayerEst}, we
perform such a rewrite for the harmonic double layer potential in two
dimensions, and show that the estimate in (\ref{eq:base_est_integer})
matches the estimate in (\ref{eqn:En-complex-w-kn}) as applied to that
reformulated integral.

\subsubsection{Half-integer $p$}
\label{subsec:halfint}

We now consider the case when $p$ is a half-integer,
$p = \bar p + 1/2$, $\bar p\in \mathbb Z$. For this, we will be
following the approach of \cite{Elliott2008}.

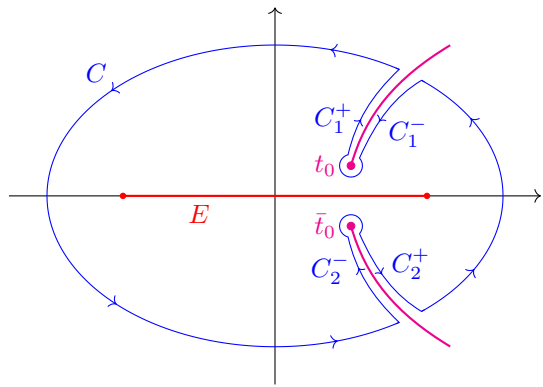
\begin{figure}
\begin{center}
\begin{tikzpicture}[scale=1]

% dati numerici
\tikzmath{
    % dati cerchietto primo quadrante
    \Cx = 1; \Cy = 0.4; \R = .15;
    %
    % angoli speciali cerchietto (in gradi)
    \ca1 = 42; % inizio buco
    \ca2 = 105; % fine buco
    %
    % angoli speciali ellisse (in gradi)
    \ea1 = 30; \succea1 = \ea1 + 1;
    \ea2 = 50; \succea2 = \ea2 + 1;
    \ea3 = 57; \succea3 = \ea3 + 1;
    \ea4 = 75; \succea4 = \ea4 + 1;
    \ea5 = 135; \succea5 = \ea5 + 1;
    \ea6 = 360-\ea5; \succea6 = \ea6 + 1;
    \ea7 = 360-\ea4; \succea7 = \ea7 + 1;
    \ea8 = 360-\ea3; \succea8 = \ea6 + 1;
    \ea9 = 360-\ea2; \succea9 = \ea7 + 1;
    \ea0 = 360-\ea1; \succea0 = \ea0 + 1;
}

% assi cartesiani
\draw[->] (-3.5,0) -- (3.5,0);
\draw[->] (0,-2.5) -- (0,2.5);

% ellisse guida
% \draw[dotted] (0,0) ellipse (3 and 2);

% definisci tutti i punti dell'ellisse
% (P1), (P2), ..., (P360) corrispondenti ai gradi
\foreach \a in {0, ..., 360} {
    \tikzmath {
        \x = 3*cos(\a);
        \y = 2*sin(\a);
    }
    \coordinate (P\a) at (\x,\y);
}

% cerchietti guida
%\draw[dotted] (\Cx,\Cy) circle (\R);
%\draw[dotted] (\Cx,-\Cy) circle (\R);

% definisci tutti i punti del cerchietto sopra
% (Q1), (Q2), ..., (Q360) corrispondenti ai gradi
\foreach \a in {0, ..., 360} {
    \tikzmath {
        \x = \Cx + \R*cos(\a);
        \y = \Cy + \R*sin(\a);
    }
    \coordinate (Q\a) at (\x,\y);
}

% definisci tutti i punti del cerchietto sotto
% (R1), (R2), ..., (R360) corrispondenti ai gradi
\foreach \a in {0, ..., 360} {
    \tikzmath {
        \x = \Cx + \R*cos(\a);
        \y = -\Cy - \R*sin(\a);
    }
    \coordinate (R\a) at (\x,\y);
}

% pallino magenta superiore
\draw[fill, magenta] (\Cx,\Cy) circle (0.05) node[left, xshift=-2pt] {$t_0$};

% curva magenta superiore
\coordinate (ext1) at ($(P53)+(.5,.4)$);

\coordinate (CC1) at 
($(\Cx,\Cy)!.3!(ext1)+(-.2,.2)$);

\coordinate (CC2) at 
($(\Cx,\Cy)!.6!(ext1)+(-.2,.2)$);

\draw[thick, magenta] 
(\Cx,\Cy) .. controls (CC1) and (CC2) .. (ext1);

% pallino magenta inferiore
\draw[fill, magenta] (\Cx,-\Cy) circle (0.05)
node[left, xshift=-2pt] {$\bar{t}_0$};

% curva magenta inferiore
\coordinate (ext2) at ($(P307)+(.5,-.4)$);

\coordinate (CC3) at 
($(\Cx,-\Cy)!.3!(ext2)+(-.2,-.2)$);

\coordinate (CC4) at 
($(\Cx,-\Cy)!.6!(ext2)+(-.2,-.2)$);

\draw[thick, magenta] 
(\Cx,-\Cy) .. controls (CC3) and (CC4) .. (ext2);

% INIZIO PERCORSO BLU

% primo trato di ellisse
\draw[blue] (P0)
\foreach \a in {1, ..., \ea1} {
    -- (P\a)
};

% prima freccetta
\draw[blue, ->] (P\ea1) -- (P\succea1);

% secondo tratto di ellisse
\draw[blue] (P\ea1)
\foreach \a in {\succea1, ..., \ea2} {
    -- (P\a)
};

% C_1^-
\coordinate (C1) at ($(P\ea2)!.3!(Q\ca1)+(-.1,.1)$);

\coordinate (C2) at ($(P\ea2)!.6!(Q\ca1)+(-.1,.1)$);

\tikzmath{
    \precea2 = \ea2 -1;
    \precca1 = \ca1 - 1;
}

\draw[blue] (P\precea2) -- (P\ea2) 
.. controls (C1) and (C2) .. (Q\ca1) -- (Q\precca1);

\draw[blue, -{>[scale=.8]}] 
($(C2)+(.05,0)$) -- +(-.02,-.03) 
node[right, yshift=-4] {$C_1^-$};

% cerchietto blu
\draw[blue] (Q\ca1)
\foreach \a in 
{\ca1, ..., 0, 360, 359, ..., \ca2} {
    -- (Q\a)
};

% C_1*ì
\coordinate (C3) at ($(Q\ca2)!.3!(P\ea3)+(-.1,.1)$);

\coordinate (C4) at ($(Q\ca2)!.6!(P\ea3)+(-.1,.1)$);

\tikzmath{
    \succea3 = \ea3 + 1;
    \succca2 = \ca2 + 1;
}

\draw[blue] (Q\succca2) -- (Q\ca2) 
.. controls (C3) and (C4) .. (P\ea3)
-- (P\succea3);

\draw[blue, -{>[scale=.8]}] 
($(C3)+(.05,0)$) -- +(.02,.03) 
node[left, xshift = .5, yshift=.4] {$C_1^+$};

% nuovo tratto di ellisse
\draw[blue] (P\ea3)
\foreach \a in {\ea3, ..., \ea4} {
    -- (P\a)
};

\draw[blue, ->] (P\ea4) -- (P\succea4);

\draw[blue] (P\ea4)
\foreach \a in {\ea4, ..., \ea5} {
    -- (P\a)
};

\draw[blue, ->] (P\ea5) -- (P\succea5) 
node [above left, xshift=2pt] {$C$};

\draw[blue] (P\ea5)
\foreach \a in {\ea5, ..., \ea6} {
    -- (P\a)
};

\draw[blue, ->] (P\ea6) -- (P\succea6);

\draw[blue] (P\ea6)
\foreach \a in {\ea6, ..., \ea7} {
    -- (P\a)
};

\draw[blue, ->] (P\ea7) -- (P\succea7);

% C_2^-
\coordinate (C5) at 
($(R\ca2)!.3!(P\ea8)+(-.1,-.1)$);

\coordinate (C6) at ($(R\ca2)!.6!(P\ea8)+(-.1,-.1)$);

\tikzmath{
    \precea8 = \ea8 - 1;
    \succca2 = \ca2 + 1;
}

\draw[blue] (R\succca2) -- (R\ca2)
.. controls (C5) and (C6) .. (P\ea8)
-- (P\precea8);

\draw[blue, -{>[scale=.8]}] 
($(C5)+(.05,0)$) -- +(-.02,+.03) 
node[left, xshift = .5, yshift=-.4] {$C_2^-$};

% cerchietto blu
\draw[blue] (R\ca1)
\foreach \a in 
{\ca1, ..., 0, 360, 359, ..., \ca2} {
    -- (R\a)
};

% C_2^-
\coordinate (C7) at ($(P\ea9)!.3!(R\ca1)+(-.1,-.1)$);

\coordinate (C8) at
($(P\ea9)!.6!(R\ca1)+(-.1,-.1)$);

\tikzmath{
    \succea9 = \ea9 + 1;
    \precca1 = \ca1 - 1;
}

\draw[blue] (P\succea9) -- (P\ea9) 
.. controls (C7) and (C8) .. (R\ca1)
-- (R\precca1);

\draw[blue, -{>[scale=.8]}] 
($(C8)+(.05,0)$) -- +(.02,-.03) 
node[right, yshift=4] {$C_2^+$};

\draw[blue] (P\ea7)
\foreach \a in {\ea7, ..., \ea8} {
    -- (P\a)
};

\draw[blue] (P\ea9)
\foreach \a in {\ea9, ..., \ea0} {
    -- (P\a)
};

\draw[blue, ->] (P\ea0) -- (P\succea0);

\draw[blue] (P\ea0)
\foreach \a in {\ea0, ..., 360} {
    -- (P\a)
};

% segmento rosso
\draw[thick, red] (-2,0)  -- (2,0);

% pallocchi rossi
\draw[fill, red] (-2,0) circle (1pt);
\draw[fill, red] (2,0) circle (1pt);

% etichetta rossa
\draw[red] (-1,0) node[below] {$E$};
\end{tikzpicture}
\end{center}
 \caption{{The contour $C$, the deformations $C_1$ and $C_2$ going around the singularities $t_0$ and $\bar{t_0}$.}}
 \label{fig:contourC}
\end{figure}

Consider again the integral \eqref{eq:base_integral_contour}, with the
contour $C$ {depicted in \cref{fig:contourC}}. The integrand now
has singularities of the form $(t-t_0)^{\bar p+1/2}$, with branch
points at the singularities. Since these singularities are not poles,
we can no longer use residue calculus. Instead, we let $C_1$ and $C_2$
be the deformations of $C$ going around $t_0$ and $\conj t_0$ respectively, following
the branch cuts going from the singularities. We now let $C$ go to infinity,
and again, based on the assumptions introduced below
(\ref{eq:base_integral_contour}), consider only the contributions from $C_1$ and $C_2$,
\begin{align}
\oprem \approx \frac{1}{2\pi i}\int_{C_1}  \frac{f(t) k_n(t)}{\norm{\v\gamma(t) - \v x}^{2p}} \dif t
  + \frac{1}{2\pi i}\int_{C_2}  \frac{f(t) k_n(t)}{\norm{\v\gamma(t) - \v x}^{2p}} \dif t
  =: E_1 + E_2.
\end{align}

Now consider the contribution from $C_1$. We multiply and divide the
integrand with $(t-t_0)^p$ and integrate by parts $\bar p$ times,
(ignoring endpoint contributions, since we are considering a section of
a closed contour)
\begin{align}
  E_1
  & = 
    \frac{1}{2\pi i} \int_{C_1}
    \frac{1}{(t-t_0)^p}
    \frac{(t-t_0)^pf(t)k_n(t)}{\norm{\v\gamma(t) - \v x}^{2p}}
    \dif t
  \\
  & = 
    \frac{1}{2\pi i}
    \frac{1}{\prod_{q=1}^{\bar p}(q-p)}
    \int_{C_1} 
    \frac{1}{\sqrt{t-t_0}}
    \dod[\bar p]{}{t}
    \pars{
    \frac{(t-t_0)^pf(t)k_n(t)}{\norm{\v\gamma(t) - \v x}^{2p}}
    }
    \dif t.
\label{eq:branch_partint}
\end{align}
We can simplify this further using \cref{lemma:main} in
\cref{app:lemmas},
\begin{align}
  \frac{1}{\prod_{q=1}^{\bar p}(q-p)} = \frac{\Gamma(1-p)}{\sqrt{\pi}} .
\end{align}
Note that the factor $(t-t_0)^pf(t)/\norm{\v\gamma(t) - \v x}^{2p}$ is
smooth on $C_1$, and we assume that it varies much slower than
$k_n$. Analogous to the integer case, we first simplify by only
differentiating $k_n$, and then we approximate the smooth part with
its value at $t_0$ (where $k_n$ is largest),
\begin{align}
  E_1
  & \approx \frac{1}{2\pi i}
    \frac{\Gamma(1-p)}{\sqrt{\pi}} 
    G(t_0)^p f(t_0)
    \int_{C_1} \frac{k_n^{(\bar p)}(t)}{\sqrt{t-t_0}} \dif t .
\end{align}
Denote the sides of $C_1$ by $C_1^+$ (going out) and $C_1^-$ (going
in). Defining the jump across the branch cut as
\begin{align}
  (t-t_0)\Big|_{C_1^+} = (t-t_0)\Big|_{C_1^-} e^{-2\pi i},
\end{align}
we have that
\begin{align}
  (t-t_0)^{-1/2} \Big|_{C_1^+} = -(t-t_0)^{-1/2} \Big|_{C_1^-},
\end{align}
which lets us write the $C_1$ contribution as
\begin{align}
  \int_{C_1} \frac{k_n^{(\bar p)}(t)}{\sqrt{t-t_0}} \dif t
  &=
    \int_{t_0}^\infty \frac{k_n^{(\bar p)}(t)}{\sqrt{t-t_0}  \Big|_{C_1^+}} \dif t
    + \int_{\infty}^{t_0} \frac{k_n^{(\bar p)}(t)}{\sqrt{t-t_0}  \Big|_{C_1^-}} \dif t
  =
    -2\int_{t_0}^\infty \frac{k_n^{(\bar p)}(t)}{\sqrt{t-t_0}} \dif t
.
\end{align}
Going back to using $p$ rather than $\bar p$, we define
\begin{align}
  J(t_0,n,p) = \int_{t_0}^\infty \frac{k_n^{(p-1/2)}(t)}{\sqrt{t-t_0}} \dif t,
  \label{eq:J_def}
\end{align}
where the integration from $t_0$ to $\infty$ is to follow the
branch cut. With this, we get
\begin{align}
  E_1
  & \approx 
    \frac{\Gamma(1-p)}{i\pi^{3/2}} 
    G(t_0)^p f(t_0)
    J(t_0,n,p) .
\end{align}
Repeating the calculations for $C_2$ with the pole at $\bar{t}_0$, we
find $E_2=\bar{E}_1$. 
Similar to the previous section, we use the triangle
inequality when adding up $E_1$ and $E_2$ to get
a slightly conservative estimate. In addition, we simplify using the relation
$\abs{\Gamma(1-p)}=\pi/\Gamma(p)$, which is a special case of Euler's
reflection formula for $p$ half-integer. 
Our final form for the absolute value of the error is then
\begin{align}
  \abs{\oprem} \approx \frac{2}{ \Gamma(p) \sqrt{\pi}} 
  \abs{J(t_0,n,p) G(t_0)^p f(t_0)} =
\frac{2}{ \Gamma(p) \sqrt{\pi}}   \abs{J(t_0,n,p)} \abs{G(t_0)}^p \abs{f(t_0)}.
  \label{eq:base_est_halfinteger}
\end{align}
In order for this expression to be useful, a closed-form estimate for
$J(t_0,n,p)$ is needed. Such an estimate can be derived by defining
a suitable branch cut with respect to the error function $k_n$, as we
shall see in the cases of the trapezoidal and Gauss-Legendre
quadrature rules.

\subsection{Trapezoidal rule}
\label{sec:trapezoidal-rule}

For the trapezoidal rule, we are considering the integral in
(\ref{eq:cartesian_int}), with the integration interval $E=[0,2\pi)$
and an integrand that is assumed to be periodic in $t$. The corresponding
error function is given in (\ref{eq:knTz}), with an asymptotic form
for $n\to\infty$ in  (\ref{eq:kn_trapz_per_asy}). The derivatives of
this function is given in (\ref{eq:knq_trapz_per_asy}) with a
somewhat simpler expression for the magnitude in (\ref{eq:knp_trapz_abs}). 

%\Cref{eq:kn_trapz_per_asy,eq:knq_trapz_per_asy,eq:knp_trapz_abs} are
%good approximations to \eqref{eq:knTz} as long as
%$e^{n|\Im z|} \gg 1$. This is also a requirement for the error to be
%less than $\mathcal O(1)$, so they useful in most practical
%applications.

\subsubsection{Trapezoidal rule with integer $p$}

For the trapezoidal rule with integer $p$, formulating an error estimate is just a matter of combining \eqref{eq:base_est_integer} and \eqref{eq:knp_trapz_abs}, giving \cref{est:trapz_p_int}.

\begin{estimate}[Trapezoidal rule with integer $p$]
  \label{est:trapz_p_int}
Consider the integral in (\ref{eq:cartesian_int}), where $\gamma(E)$ is the parameterization
of a smooth closed curve in $\mathbb R^2$ or $\mathbb R^3$. 
The integrand is assumed to be periodic in $t$ over the integration interval  $E=[0,2\pi)$.
The error in approximating the integral with the $n$-point trapezoidal
rule can in the limit $n \rightarrow \infty$ be estimated as
  \begin{align}
    \abs{\opremTh} \approx  
    \frac{4\pi n^{p-1}}{(p-1)!} \abs{f(t_0)} \abs{G(t_0)}^p e^{-n|\Im t_0|}.
    \label{eq:trapz_p_int}
  \end{align}
Here, $p$ is a positive integer, and the geometry
factor $G$ is defined in  (\ref{eq:geometry_factor}). 
The squared distance function is defined in (\ref{eq:R2_def}), and 
$\braces{t_0, \conj t_0}$ is the pair of complex conjugate roots of this
$ R^2(t,\v x)$ closest to the integration interval $E$.
\end{estimate}

\subsubsection{Trapezoidal rule with half-integer $p$}

For the trapezoidal rule with half-integer $p$, we must derive an
expression for $J(t_0,n,p)$, as defined in \eqref{eq:J_def}, with
derivatives of $k_n$ as given in \eqref{eq:knq_trapz_per_asy}. We can without loss of
generality assume that $\Im t_0>0$. Let the branch cut going from $t_0$ to infinity be
\begin{align}
  B(t_0)=\{t(s) \in \mathbb{C} : t(s)=t_0+is, \ \ 0 \leq s < \infty \},
\end{align}
and let this be the branch cut enclosed by the path $C_1$ in
\eqref{eq:branch_partint}. Note that along this cut
\begin{align}
  t'(s) &= i, \\
  t(s)-t_0 &= is,
\end{align}
such that
\begin{align}
  J(t_0,n,p)
  =  \sqrt{i} \int_{0}^\infty \frac{k_n^{(p-1/2)}(t_0+is)}{\sqrt{s}}\dif s 
  \approx  -2\pi i \sqrt{i} (in)^{p-1/2} e^{i n t_0}
    \underbrace{
    \int_{0}^\infty \frac{
    e^{is}
    }{\sqrt{s}}\dif s
    }_{\sqrt{\pi/n}}.
\end{align}
Considering only the absolute value,
\begin{align}
  \abs{J(z_0,n,p)} \approx 
  2 \pi^{3/2} n^{p-1} e^{-n|\Im t_0|},
\end{align}
finally yields \cref{est:trapz_p_halfint}.

\begin{estimate}[Trapezoidal rule with half-integer $p$]
  \label{est:trapz_p_halfint}
Consider the integral in (\ref{eq:cartesian_int}), where $\gamma(E)$ is the parameterization
of a smooth closed curve in $\mathbb R^2$ or $\mathbb R^3$. 
The integrand is assumed to be periodic in $t$ over the integration interval  $E=[0,2\pi)$.
The error in approximating the integral with the $n$-point trapezoidal
rule can in the limit $n \rightarrow \infty$ be estimated as 
 \begin{align}
  \abs{\opremTh} 
 \approx  
    \frac{4\pi n^{p-1} }{\Gamma(p)}  \abs{f(t_0)} \abs{G(t_0)}^p e^{-n|\Im t_0|}.
    \label{eq:trapz_p_halfint}
  \end{align}
Here, $p$ is a positive half-integer, $\Gamma(p)$ the gamma function, and the geometry
factor $G$ is defined in  (\ref{eq:geometry_factor}). 
The squared distance function is defined in (\ref{eq:R2_def}), and 
$\braces{t_0, \conj t_0}$ is the pair of complex conjugate roots of this
$ R^2(t,\v x)$ closest to the integration interval $E$. 
\end{estimate}
Interestingly, this is identical to the estimate
\eqref{eq:trapz_p_int} for integer $p$, if we generalize the factorial
to non-integer $p$ as $(p-1)!=\Gamma(p)$. This generalization can be
found in \cite{AfKlinteberg2016quad} for both the trapezoidal and
Gauss-Legendre rules, where it was noted that it works well for
half-integer $p$. What we have shown here (and will show in \cref{sec:gl-half-integer}) is why it works well.

\subsection{Gauss-Legendre rule}
\label{sec:gauss-legendre-rule}

For the Gauss-Legendre quadrature rule we consider the integral in 
(\ref{eq:cartesian_int}) over the base interval $E=[-1,1]$. 
The error function is not available in closed form, but can in the
limit $n\to\infty$ be shown to asymptotically satisfy the formula (\ref{eq:knGL})
\cite{Elliott2008}, here written as
\begin{align}
  k_n(z) \simeq \frac{2\pi}{\xi(z)^{2n+1}},
  \label{eq:kn_gl_asy}
\end{align}
where
 \begin{align}
  \xi(z) = z + \sqrt{z^2-1} .
  \label{eq:xi_def}
\end{align}
As was introduced below equation (\ref{eq:knGL}),
 $\sqrt{z^2-1}$ is defined as
$\sqrt{z+1}\sqrt{z-1}$ with $-\pi < \arg(z \pm 1) \le \pi$
\cite{Elliott2008}. Alternatively, we can write
$\xi(z)=z \pm \sqrt{z^2-1}$, with the sign defined such that
$|\xi| \ge 1$.
The approximation of the derivatives of $k_n(z)$ as introduced in
(\ref{eq:knq_gl_asy}) will include the same factor of $1/\xi(z)^{2n+1}$. 
%can furthermore be approximated as \cite{AfKlinteberg2016quad}
% \begin{align}
%  k_n^{(q)}(z) \simeq
%  \pars{-\frac{2n+1}{\sqrt{z^2-1}}}^{q}
%  \frac{2\pi}{\xi(z)^{2n+1}}.
%  \label{eq:knq_gl_asy}
%\end{align}

The main characteristic of the asymptotic Gauss-Legendre error
function~\eqref{eq:kn_gl_asy} is that its magnitude is constant on the
level sets of the function
\begin{align}
  \rho(z) = \abs{\xi(z)},
\end{align}
which we denote the \emph{Bernstein radius} of $z$. This follows from
the notion of a Bernstein ellipse, which is an ellipse with foci
$\pm 1$ where the semimajor and semiminor axes sum to $\rho>1$. It can
be constructed as the image of the circle $|\xi|=\rho$ under the
Joukowski transform
\begin{align}
  z(\xi) = \frac{\xi + \xi^{-1}}{2},
  \label{eq:joukowski}
\end{align}
which is the inverse of \eqref{eq:xi_def}.

\subsubsection{Gauss-Legendre rule with integer $p$}

For the Gauss-Legendre rule with integer $p$, we can combine
\eqref{eq:base_est_integer} and~\eqref{eq:knq_gl_asy} to get the
following error estimate. 
\begin{estimate}[Gauss-Legendre rule with integer $p$]
  \label{est:gl_p_int}
 Consider the integral in (\ref{eq:cartesian_int}) , where $\gamma(E)$ is the parameterization
of a smooth curve in $\mathbb R^2$ or $\mathbb R^3$, with $E=[-1,1]$.
The error in approximating the integral with the $n$-point
Gauss-Legendre rule can in the limit $n \rightarrow \infty$ be estimated as
  \begin{align}
    \abs{\opremTh} \approx
    \frac{ 4 \pi }{(p-1)!}
    \abs{\frac{2n+1}{\sqrt{t_0^2-1}}}^{p-1} \abs{f(t_0)} \abs{G(t_0)}^p
    \frac{1}{\rho(t_0)^{2n+1}}.
    \label{eq:gl_p_int}
  \end{align}
Here, $p$ is a positive integer, the geometry
factor $G$ is defined in  (\ref{eq:geometry_factor}) and $\rho(t) = |t + \sqrt{t+1}\sqrt{t-1}|$.
The squared distance function is defined in (\ref{eq:R2_def}), and 
$\braces{t_0, \conj t_0}$ is the pair of complex conjugate roots of this
$ R^2(t,\v x)$ closest to the integration interval $E$. 
\end{estimate}

\subsubsection{Gauss-Legendre rule with half-integer $p$}
\label{sec:gl-half-integer}

In order to evaluate the integral \eqref{eq:J_def} with the
Gauss-Legendre error function \eqref{eq:knq_gl_asy}, we will closely
follow the steps outlined by Elliott et al. \cite{Elliott2008}.  We
define a branch cut going from $t_0$ to infinity (again, assuming
without loss of generality that $\Im t_0>0$) using a scaled Joukowsky transform,
\begin{equation}
B(t_0)=\left\{
  t(s) \in \mathbb{C} : 
  t(s) = \frac{1}{2} \pars{\zeta(s) + \frac{1}{\zeta(s)}}, \ \  1 \le s < \infty
\right\},
\end{equation}
where $\zeta(s)$ parameterizes a radial line from the point
$\zeta_0=t_0+\sqrt{t_0^2-1}$ to infinity, 
\begin{align}
  \zeta(s) = \zeta_0 s
\end{align}
and $\rho(t_0)=|\zeta_0|$ as introduced above.
This parametrization of the branch cut satisfies $t(1)=t_0$, and leads
to the following useful relations,
\begin{align}
  t'(s) &= \frac{\zeta_0}{2}\pars{1-\frac{1}{\zeta_0^2s^2}}, \\
  t(s)-t_0 &=  \frac{\zeta_0(s-1)}{2} \pars{1 - \frac{1}{\zeta_0^2 s}}, \\
  \sqrt{t(s)^2-1} &= \frac{\zeta(s)}{2}\pars{1 - \frac{1}{\zeta(s)^2}} = s t'(s) .\\
\end{align}
Together with the relation $\xi(t(s)) = \zeta(s)$, this lets us write
\eqref{eq:knq_gl_asy} as
\begin{align}
  k_n^{(\bar p)}(t(s)) \approx 
  \pars{-
  \frac{2n+1}
  { st'(s) }
  }^{\bar p} \frac{2\pi}{\zeta(s)^{2n+1}}. 
\end{align}
With substitution of the above relations, reverting to use only
$p=\bar p +1/2$ we can write
\eqref{eq:J_def} as
\begin{align}
  \begin{split}
    &J(t_0,n,p) = \int_1^\infty \frac{k_n^{(p-1/2)}
      (t(s))}{\sqrt{t(s)-t_0}} t'(s)
    \dif s \\
    & \approx
    \pars{-(4n+2)}^{p-1/2} \frac{\pi\sqrt{2}}{\zeta_0^{2n+p}}
    \int_1^\infty \frac{1}{s^{2n+p+\frac 12}\sqrt{s-1}}
    \frac{1}{\sqrt{1 - \frac{1}{\zeta_0^2 s}}}
    \pars{1-\frac{1}{\zeta_0^2s^2}}^{3/2-p} \dif s .
  \end{split}
    % LK CHECKED OK 5 May
\end{align}
By assuming that most of the contribution to the integral comes from
the neighborhood of $s=1$, we make the simplifications
%(upper bounds, actually)
\begin{align}
  1 - \frac{1}{\zeta_0^2 s^2} \approx
  1 - \frac{1}{\zeta_0^2 s} \approx 
  1 - \frac{1}{\zeta_0^2}  =
  \frac{2}{\zeta_0} \sqrt{t_0^2-1}
  .
\end{align}
Then,
\begin{align}
  J(t_0,n,p) \approx 
  \frac{2 \pi i \sqrt{(2n+1)}}{\zeta_0^{2n+1}}
  \pars{ -\frac{2n+1}{\sqrt{t_0^2-1}} }^{p-1}  
  \int_1^\infty
  \frac{1}{s^{2n+p+\frac 12}\sqrt{s-1}}
  \dif s .      
\end{align}
This can be simplified using the result from \cite{Elliott2008} that for
$N$ large,
\begin{align}
  \int_1^\infty \frac{\dif s}{s^N\sqrt{s-1}} \approx
  \sqrt{\frac{\pi}{N}},
\end{align}
such that, taking the absolute value and reintroducing
$|\zeta_0|=\rho(t_0)$ we get
\begin{align}
  \abs{J(t_0,n,p)} \approx 
  \pi^{3/2}
  \sqrt{\frac{4(2n+1)}{2n+p+1/2}}
  \abs{ \frac{2n+1}{\sqrt{t_0^2-1}} }^{p-1}    
  \frac{1}{\rho(t_0)^{2n+1}}.
\end{align}
Assuming $n$ large and $p \ge 1/2$ of moderate values, we approximate
$(2n+1)/(2n+1/2+p) \approx 1$, such that
\begin{align}
  \abs{J(t_0,n,p)} \approx 
  2 \pi^{3/2}
  \abs{ \frac{2n+1}{\sqrt{t_0^2-1}} }^{p-1}    
  \frac{1}{\rho(t_0)^{2n+1}}.
\end{align}
Inserting this into \eqref{eq:base_est_halfinteger},
we get
\begin{estimate}[Gauss-Legendre rule with half-integer $p$]
  \label{est:gl_p_halfint}
 Consider the integral in (\ref{eq:cartesian_int}) , where $\gamma(E)$ is the parameterization
of a smooth curve in $\mathbb R^2$ or $\mathbb R^3$, with $E=[-1,1]$.
The error in approximating the integral with the $n$-point Gauss-Legendre
rule can in the limit $n \rightarrow \infty$ be estimated as
  \begin{align}
    \abs{\opremTh} \approx     
    \frac{4\pi}{ \Gamma(p)}   
    \abs{ \frac{2n+1}{\sqrt{t_0^2-1}} }^{p-1}   \abs{f(t_0)}  \abs{G(t_0)}^p
    \frac{1}{\rho(t_0)^{2n+1}}.
    \label{eq:gl_p_halfint}
  \end{align}
Here, $p$ is a positive half-integer, $\Gamma(p)$ the gamma function, the geometry
factor $G$ is defined in  (\ref{eq:geometry_factor}) and $\rho(t) = |t + \sqrt{t+1}\sqrt{t-1}|$.
The squared distance function is defined in (\ref{eq:R2_def}), and 
$\braces{t_0, \conj t_0}$ is the pair of complex conjugate roots of this
$ R^2(t,\v x)$ closest to the integration interval $E$. 
\end{estimate}
Analogously to the trapezoidal rule case, this estimate is identical
to the estimate \eqref{eq:gl_p_int} for integer $p$, with the
generalization $(p-1)!=\Gamma(p)$.

%\subsection{Results for planar curves}
\subsection{Examples for one-dimensional curves in $\reals^2$}
\label{sec:curve-examples}

The estimates in \cref{sec:trapezoidal-rule,sec:gauss-legendre-rule}
for the nearly singular quadrature error are derived using a number of
simplifications, in order to get closed-form expressions without any
unknown constants. Nevertheless, they have excellent predictive
accuracy. To demonstrate this, we consider the simple layer potential
\begin{align}
  u(\v x) = \int_0^{2\pi} 
  \frac{\norm{\v\gamma'(t)} \dif t}
  {\norm{\v\gamma(t)-\v x}^{2p}},
  \label{eq:curve-example-integral}
\end{align}
for $p=\frac{1}{2}, 1, \frac{3}{2}, 2$ and points $\v x$ near the planar curve defined by
\begin{align}
  \v\gamma(t) = \pars{1 + 0.1\cos(5t)}
  \begin{pmatrix}
    \cos(t) \\ \sin(t)
  \end{pmatrix}, \quad t \in [0,2\pi) .
  \label{eq:planar_potato}
\end{align}
We discretize the integration interval $[0, 2\pi)$ in two ways: using
the composite Gauss-Legendre method with 20 equisized panels and $n$
points per panel, and using the global trapezoidal rule with $n$
points along the curve. This gives us the quadrature value $\opquad$
for each $\v x$. We compute the reference value using an adaptive
quadrature routine (Matlab's \texttt{integral} with \texttt{AbsTol=0}
and \texttt{RelTol=0}). Then we can compute an accurate value of
$\oprem$, which we compare to
\cref{est:trapz_p_int,est:trapz_p_halfint,est:gl_p_int,est:gl_p_halfint}.
Note that in the case of composite Gauss-Legendre, $\oprem$ is
calculated as the sum of the contributions from the 3 panels nearest
to $\v x$.

\subsubsection{Error contours}
\label{sec:twodim-error-contours}

In order to compare our estimates to the actual quadrature error, we
need target points $\v x$ for which the corresponding root $t_0$ is
known. Before entering the discussion of how to compute $t_0$, we will
demonstrate our estimates for target points $\v x$ for which $t_0$ is
analytically known. We define such points through complexification of
the curve parametrization, such that they by construction are roots to
the squared distance function \eqref{eq:R2_def}. That is, we first set
$t_0\in\mathbb C$, and then construct the corresponding $\v x$ as
\begin{align}
  \v x(t_0) =
  \begin{pmatrix}
    \Re \omega(t_0)\\
    \Im \omega(t_0)
  \end{pmatrix}
  \quad \mbox{ where } \quad
  \omega(t) = \gamma_1(t) + i\gamma_2(t).
  \label{eq:complexification}
\end{align}
\begin{figure}[htbp]
  % matlab/scripts/twodim_est_tests.m
  \centering
  \begin{subfigure}{.4\textwidth}
    \centering
    \includegraphics[width=\textwidth]{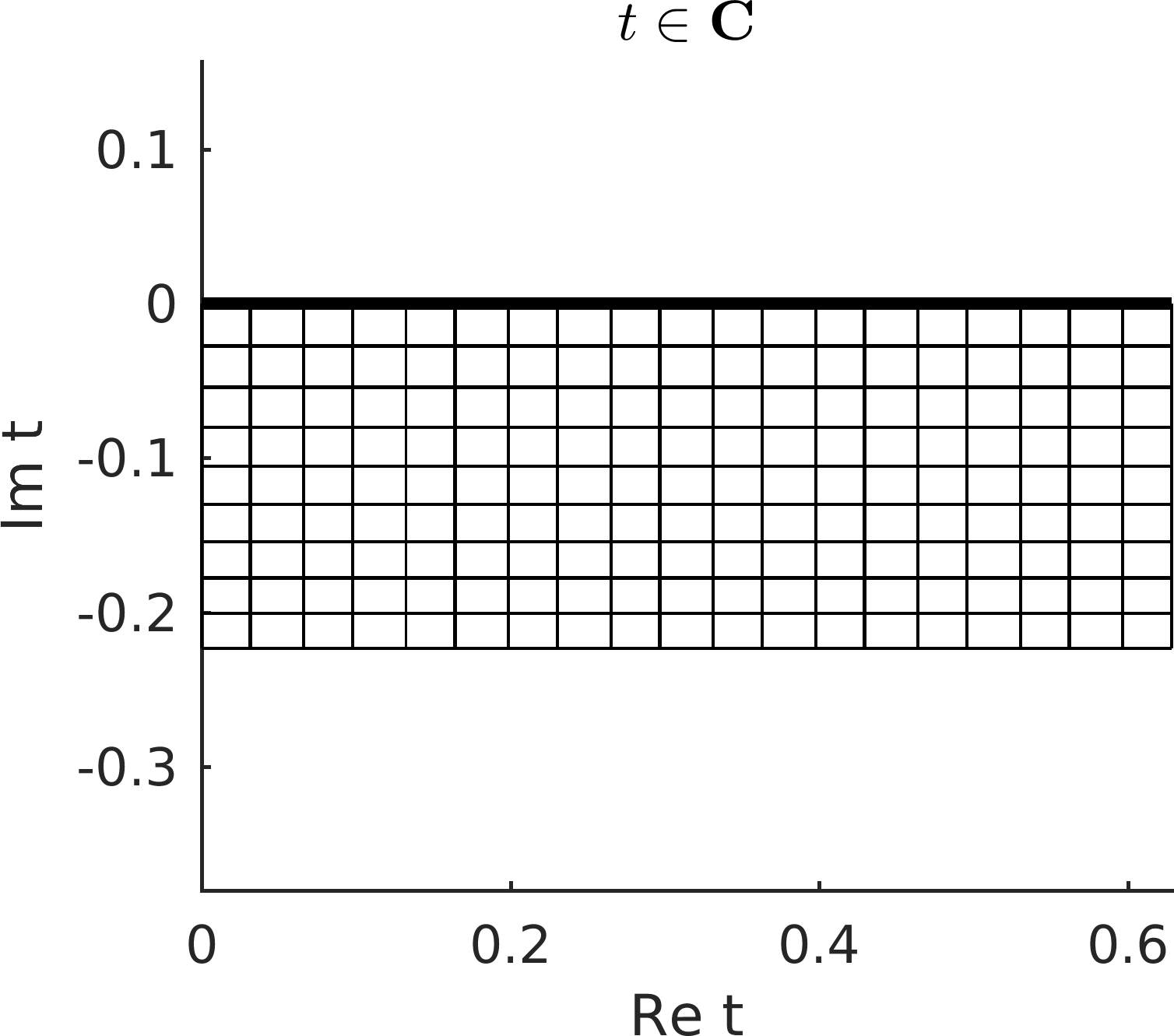}
    \caption{Grid in complex parametrization space.}
    \label{fig:field_grid_param}
  \end{subfigure}
  \hspace{.05\textwidth}
  \begin{subfigure}{.4\textwidth}
    \centering
    \includegraphics[width=\textwidth]{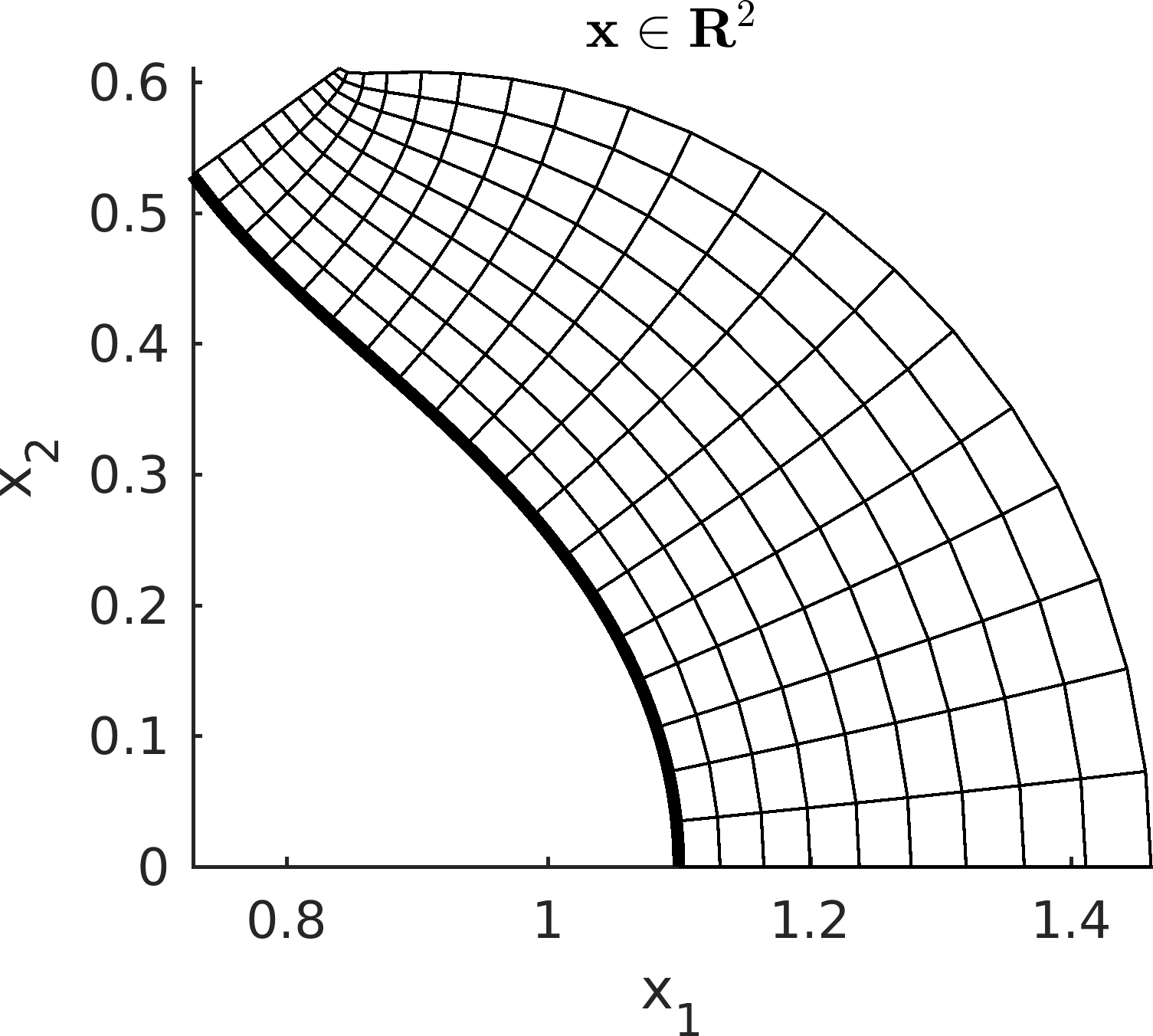}
    \caption{Image of grid in physical space.}
    \label{fig:field_grid_field}
  \end{subfigure}
  \caption{Correspondence between complex points near the real axis
    (thick line), and points near a segment of the curve parametrized
    by \eqref{eq:planar_potato}, under the complexification
    \eqref{eq:complexification}.}
  \label{fig:field_grid}
\end{figure}
This mapping is illustrated in \cref{fig:field_grid}. We construct a
grid of $200 \times 100$ points covering the region shown in that
figure, and consider the evaluation of the integral \eqref{eq:curve-example-integral} at
those target points, for $p=3/2$. We compute the integral with Gauss-Legendre
and 20 panels with 16 points each, and with the trapezoidal rule and
200 points on the entire curve. Then we estimate the quadrature error
using \eqref{eq:trapz_p_halfint} and
\eqref{eq:gl_p_halfint}. Comparing the level contours of the errors
and the estimates, see \cref{fig:gl_field_error}, it is clear that
they match well. The estimate is a smooth envelope of the
node-frequency oscillations in the quadrature error, and therefore
provides an estimated upper bound of the error. The enveloping is due
to our use of the triangle inequality when combining the errors from
the two singularities. A more precise estimate that includes the node
frequency oscillations can be obtained by skipping this step, compare
e.g. section \ref{sec:results_complex}, but it does risk
underestimating the error at some points if the oscillations do not match perfectly.
\begin{figure}[htbp]
  % matlab/scripts/twodim_est_tests.m
  \centering
  \begin{subfigure}{.49\textwidth}
    \centering
    \includegraphics[width=\textwidth]{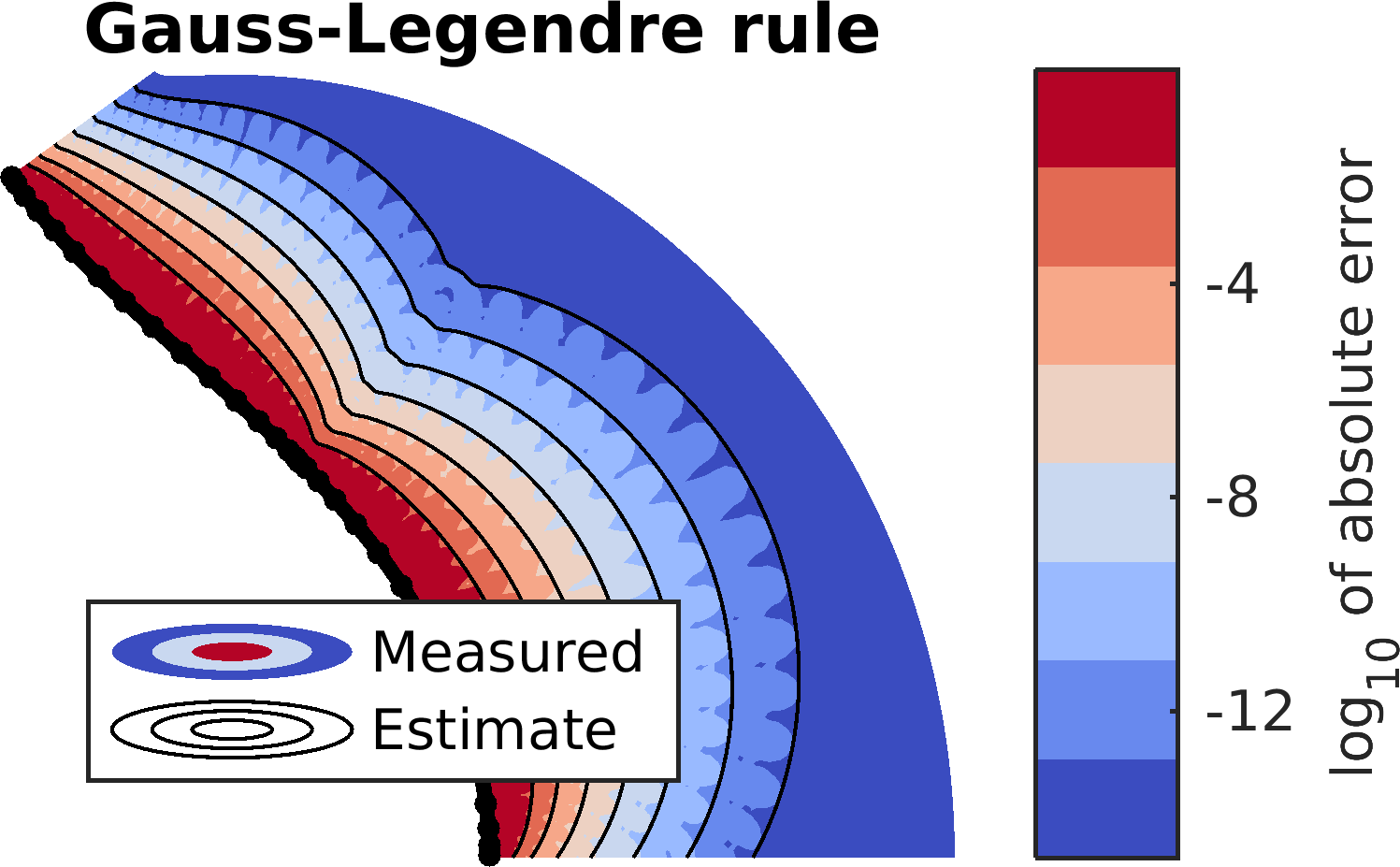}
    \caption{GL, 20 panels, $n=16$}
    \label{fig:gl_field_error}
  \end{subfigure}
  \begin{subfigure}{.49\textwidth}
    \centering
    \includegraphics[width=\textwidth]{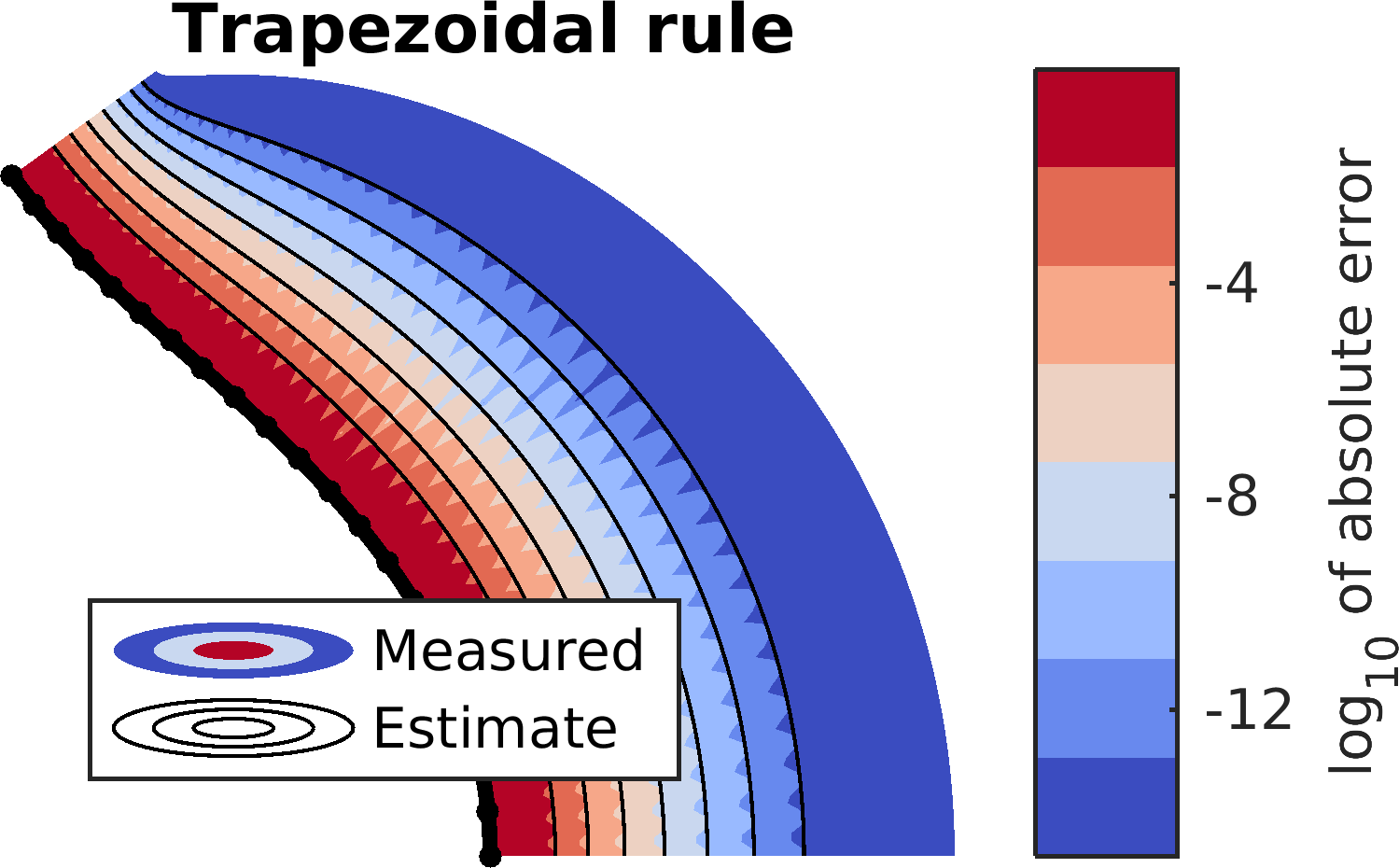}
    \caption{Trapezoidal, $n=200$}
    \label{fig:trapz_field_error}
  \end{subfigure}
  \caption{Quadrature error vs estimate for $p=3/2$.}
  \label{fig:field_err}
\end{figure}

\subsubsection{Convergence of errors and precision of error estimates}

\begin{figure}[htbp]
  % matlab/scripts/twodim_est_tests.m
  \centering
  \begin{subfigure}{.33\textwidth}
    \centering
    \includegraphics[height=\textwidth]{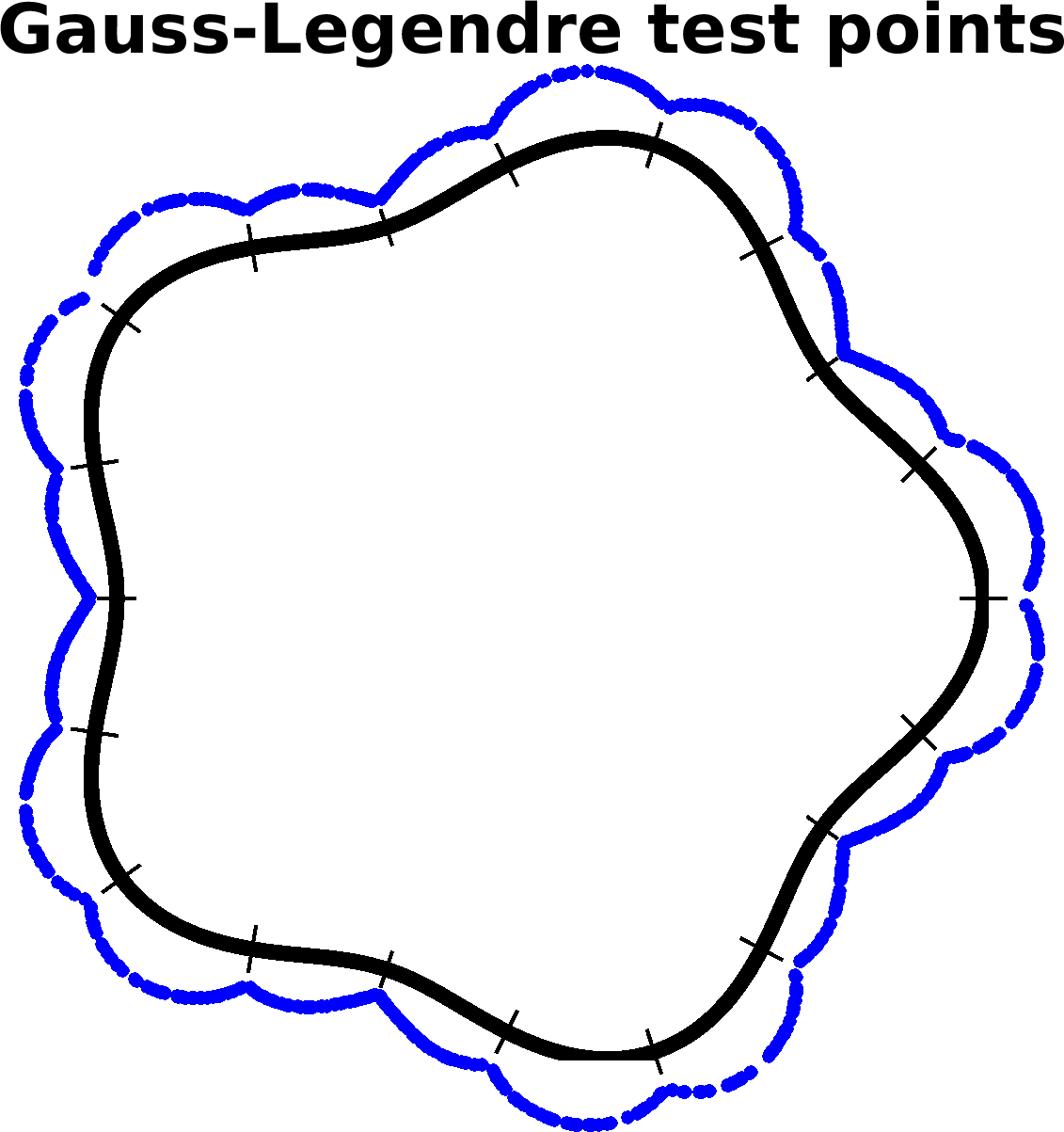}
    \caption{Target points with $\rho(t_0)=1.05$ in the
      parametrization of the nearest panel.}
    \label{fig:gl_convergence_geo}
  \end{subfigure}
  \hspace{0.1\textwidth}
  \begin{subfigure}{.33\textwidth}
    \centering
    \includegraphics[height=\textwidth]{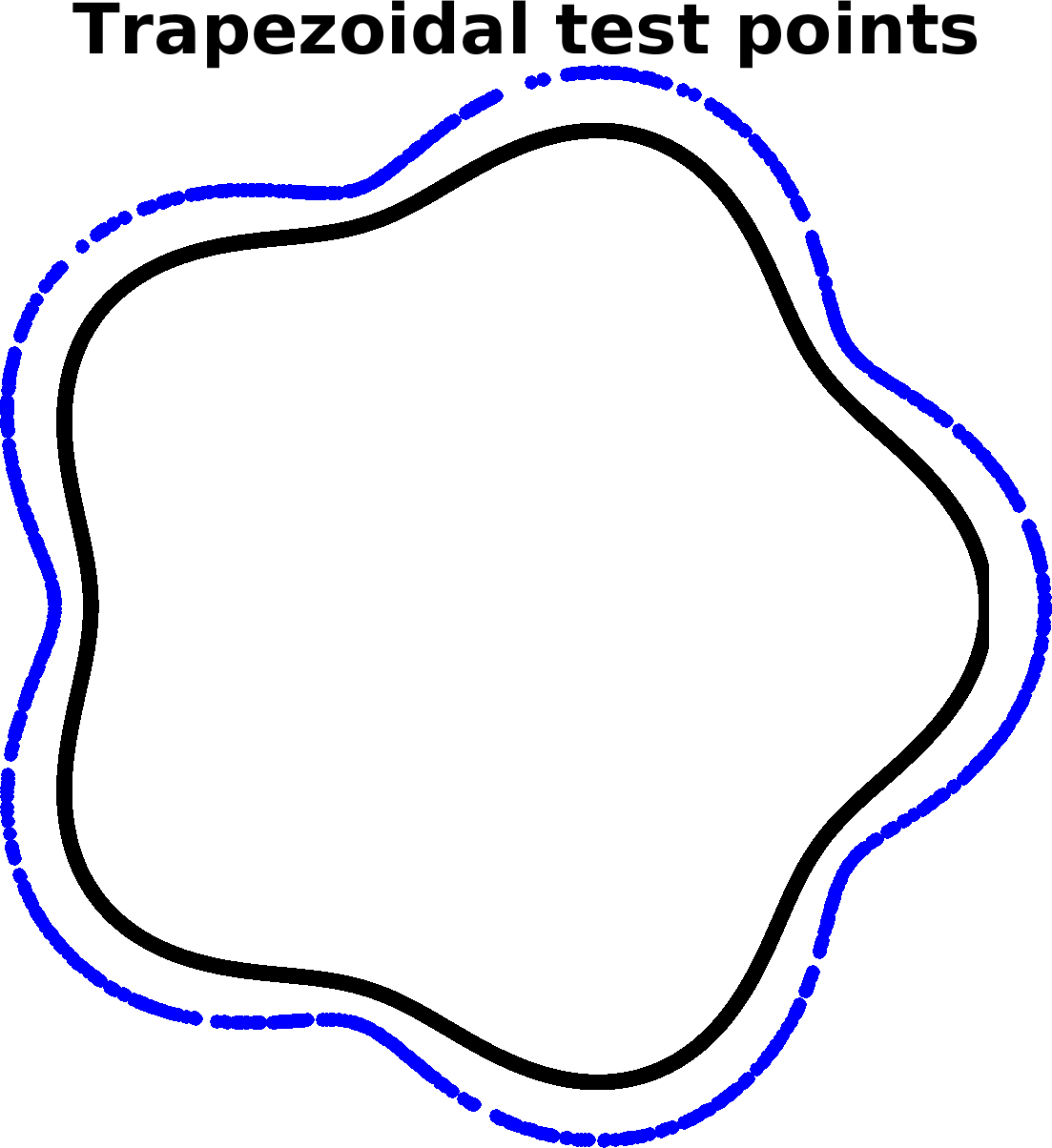}
    \caption{Target points with $|\Im t_0|=0.1$.}
    \label{fig:trapz_convergence_geo}
  \end{subfigure}
  \caption{Source geometry (black) and test points (blue) used when
    evaluating quadrature errors and error estimates for the
    Gauss-Legendre and trapezoidal rules, respectively. In each case
    the number of test points is 1000. For the Gauss-Legendre case,
    the subdivision into 20 panels is also shown.}
  \label{fig:convergence_geo}
\end{figure}

In \cref{fig:field_err} we see that the the error estimates can match
the errors well in space, for a given value of $n$ and $p$. In an attempt to show
how the error varies and how well this is captured by our estimates,
we now study the convergence with respect to $n$ for a number of
values of $p$. For this, we first attempt to construct a large number
of random points for which the magnitude of the error is of the same
magnitude, i.e. points that lie along one of the contours in
\cref{fig:field_err}. In the case of the trapezoidal rule, this is
straightforward: we simply generate random points with $|\Im t_0|$
fixed and $\Re t_0\in [0,2\pi)$. For the Gauss-Legendre rule, we need
points such that $\rho(t_0)$ has a fixed value in the parametrization
of the panel closest to the points. We create them for
each panel by applying the Joukowski transform \eqref{eq:joukowski} to
random points on a semicircle of radius $\rho$, and then keeping the
points such that $-1<\Re t_0<1$. In both cases, we determine $\v x$
from $t_0$ using the complexification
\eqref{eq:complexification}. \Cref{fig:convergence_geo} shows our test
points generated in this way.

\begin{figure}[htbp]
  \begin{subfigure}{\textwidth}
    % matlab/scripts/twodim_est_tests.m
    \centering
    \includegraphics[width=0.33\textwidth]{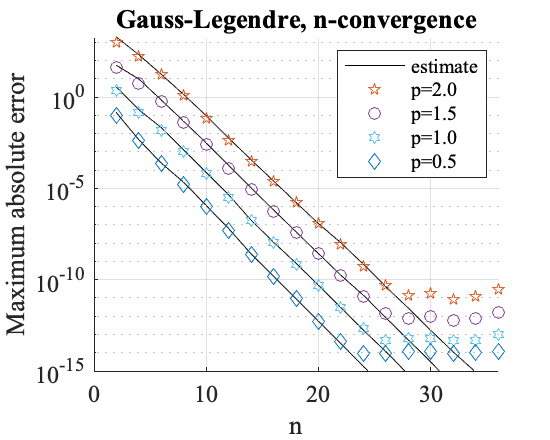}
    \includegraphics[width=0.66\textwidth]{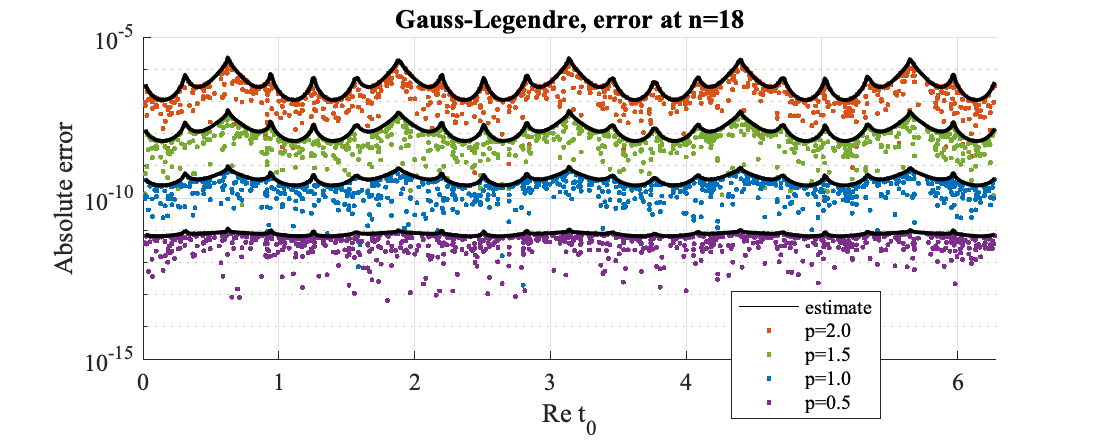}
    \caption{Gauss-Legendre convergence}
    \label{fig:gl_convergence}
  \end{subfigure}
  \begin{subfigure}{\textwidth}
    % matlab/scripts/twodim_est_tests.m
    \centering
    \includegraphics[width=0.33\textwidth]{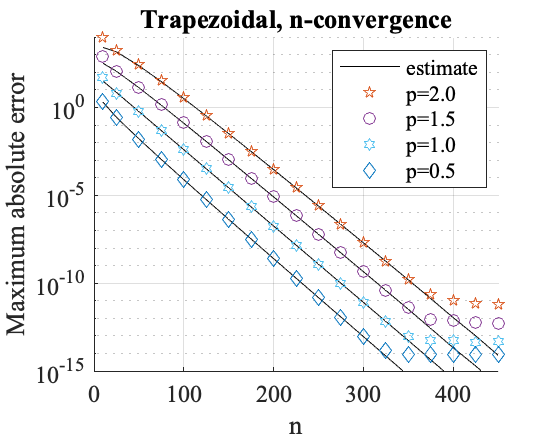}
    \includegraphics[width=0.66\textwidth]{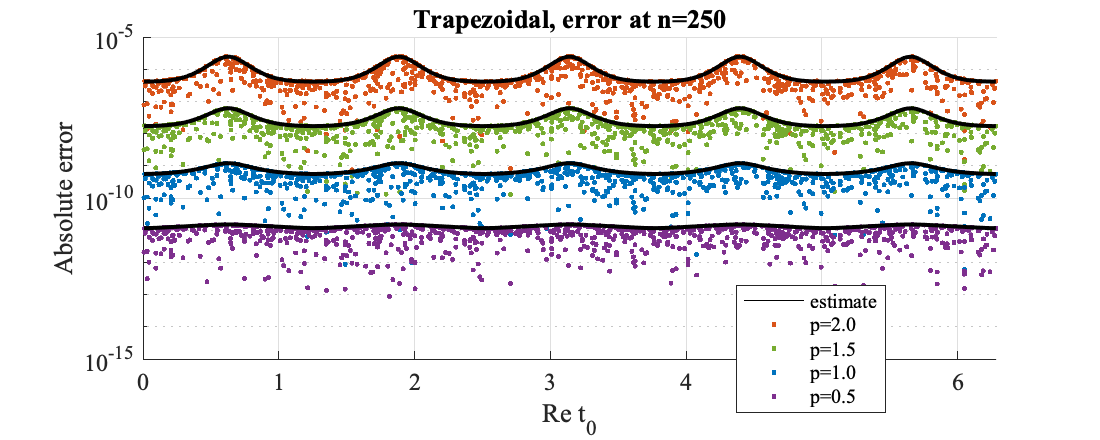}
    \caption{Trapezoidal convergence}
    \label{fig:trapz_convergence}
  \end{subfigure}
  \caption{Errors in approximation the layer potential in   (\ref{eq:curve-example-integral}) by
    the (a) Gauss-Legendre rule applied to each of $20$ panels, and (b) the trapezoidal rule. The
    errors are measured at the points shown in
    \cref{fig:convergence_geo}$a)$ and \cref{fig:convergence_geo}$b)$,
    respectively. 
For the $n$-convergence we report the value of the estimate at the point with the largest error.} 
  \label{fig:convergence}
\end{figure}

Consider the results in \cref{fig:convergence}. The plots to the right
show the errors as the sets of target points in
\cref{fig:convergence_geo} are traversed, for a number of $p$ and
fixed $n$. The estimates clearly provide a good approximate upper
bound of the error. The plots to the left shows how the maximum error
over all the target points converges towards zero as $n$ increases,
for a number of $p$. That value is compared to the estimate, which is
computed at the point with the maximum error. We see that the
estimates capture both the magnitude of of the errors and the rates of
convergence as $n$ increases.

\subsection{Root finding}
\label{sec:curve-rootfinding}

As we have seen, we can accurately predict the magnitudes of the
nearly singular quadrature errors for one-dimensional curves
discretized using the trapezoidal and Gauss-Legendre rules. However,
in order to do so for given target point $\v x$, we need to know the
location of the nearest complex root $t_0$ of the squared distance
function~\eqref{eq:R2_def}. Fortunately, finding $t_0$ numerically is
both fast and robust, using only the discrete quadrature nodes on the
curve. A method for the Gauss-Legendre case was introduced in~\cite{AfKlinteberg2018} for curves in $\mathbb R^2$, and further
developed in~\cite{AfKlinteberg2020line} for curves in both $\mathbb R^2$
and $\mathbb R^3$. We will here summarize these results, generalize
them for the trapezoidal rule, and then introduce simplifications that
will prove useful in the three-dimensional case.

In order to determine $t_0$ without explicit knowledge of the
parametrization $\v\gamma(t)$, we form an approximation to it, denoted
$\tilde{\v\gamma}(t)$. The most straightforward way of doing this, and
also the most accurate and expensive way, is to use the values at the
$n$ quadrature nodes $t_{\ell}$, $\ell=1,\dots,n$. For a
Gauss-Legendre panel $t_{\ell} \in [-1,1]$, for the trapezoidal rule
$t_{\ell} \in [0,2\pi)$. Then, an interpolant is created for each of
the $d$ components of $\v \gamma \in \reals^d$ using suitable
orthogonal basis functions. For Gauss-Legendre we use a basis
$\{P_j\}$ of orthogonal polynomials on $[-1,1]$ (e.g. Chebyshev or
Legendre),
\begin{align}
  \tilde\gamma_i(t) = \poly_n[\gamma_i](t) 
  = \sum_{j=0}^{n-1} c_j P_j(t),  \quad i=1,\ldots,d,
  \label{eq:gamma_poly}
\end{align}
while for trapezoidal we use a trigonometric polynomial,
\begin{align}
  \tilde\gamma_i(t) = \trig_n[\gamma_i](t)  
   = \sum_{k=-n/2}^{n/2-1} \hat{\gamma}_i(k) e^{ikt},
  \quad i=1,\ldots,d.
  \label{eq:gamma_trig}
\end{align}
Once we have $\tilde{\v\gamma}$, we can form an approximation to the
squared distance function in (\ref{eq:R2_def}), 
\begin{align}
  \tilde R^2(t) = \sum_{i=1}^d (\tilde\gamma_i(t)-x_i)^2 = 0.
  \label{eq:R2tilde}
\end{align}
The roots to this equation can be found to high accuracy using
Newton's method and a suitable initial guess, see discussion in
\cite{AfKlinteberg2020line}. The method typically converges rapidly,
and has a cost of $\mathcal O(n)$ per iteration, for the evaluation of
$\tilde{\v\gamma}$.  This is related to the procedure introduced
in~\cite{AfKlinteberg2018} for complex kernels as described in section
\ref{sec:CurvesComplex}. There, we create a complex-valued
approximation $\tilde\gamma(t)\approx\gamma(t)$, and given $z_0$ we
solve $\tilde\gamma(t_0)=z_0$ for the pre-image $t_0$. This procedure
can however not be generalized to three dimensions. For a planar
curve, the pre-image corresponds to one of the two roots of
$\tilde R^2(t)$.

As described, this is a reasonably efficient scheme for a
Gauss-Legendre panel, where $n$ rarely is more than $16$, but cannot
be considered efficient for the trapezoidal rule, where $n$ is the
number of points on the entire curve. This is especially true since
estimating the quadrature error is not in itself a necessary
computation, and should not incur a significant extra cost. However,
we only need to know $t_0$ with sufficient accuracy to estimate the
quadrature error to the correct order of magnitude. Thus, we can
consider ways of computing $t_0$ that are faster, but less accurate.

Let $\tilde{\v\gamma}$ approximated using \eqref{eq:gamma_poly} or
\eqref{eq:gamma_trig} be denoted the \emph{global approximation} (with
respect to the quadrature rule). Note however that it is only for the
trapezoidal rule that it is truly global, for Gauss-Legendre it
involves one full panel. 
Then, we denote the \emph{local
  approximation} to be the $q$th order Taylor expansion
\begin{align}
  \tilde\gamma_i(t) = \taylor_q[\gamma_i, t^*](t)
  = \sum_{j=0}^q \gamma_i^{(j)}(t^*)
  \frac{ \pars{t-t^*}^j }{j!},
  \label{eq:gamma_expa}
\end{align}
where $t^*$ is the value of the parametrization at the quadrature node that is closest to $\v x$, 
\begin{align}
\label{eq:t_star}
  t^* = \operatornamewithlimits{argmin}_{t \in \{t_1,\dots,t_n\}}
  \norm{\v\gamma(t) - \v x}
\end{align}
{and can be identified using e.g. a tree-based search algorithm.}

Solving \eqref{eq:R2tilde} using Newton's method, with the
approximation \eqref{eq:gamma_expa} and a moderate value of $q$,
allows us to compute $t_0$ rapidly and with sufficient accuracy for
error estimation (as we shall demonstrate). The prerequisite is that
we need to know all derivatives of $\v\gamma$ up to
$\v\gamma^{(q)}(t)$ at all quadrature nodes, which may not be
available. These can however be computed numerically at the time of
discretization, which is a one-time cost (as opposed to finding $t_0$
for all target points $\v x$, which we consider an on-the-fly
cost). The resulting root finding scheme is of course independent of
the quadrature, but is most useful for the trapezoidal rule where the
alternative of global approximation really involves the
discretization of the whole curve. 

\begin{figure}[htbp!]
  % computed using scripts/twodim_map_test.m
  \centering
  \begin{subfigure}{.32\textwidth}
    \centering
    \includegraphics[width=\textwidth]{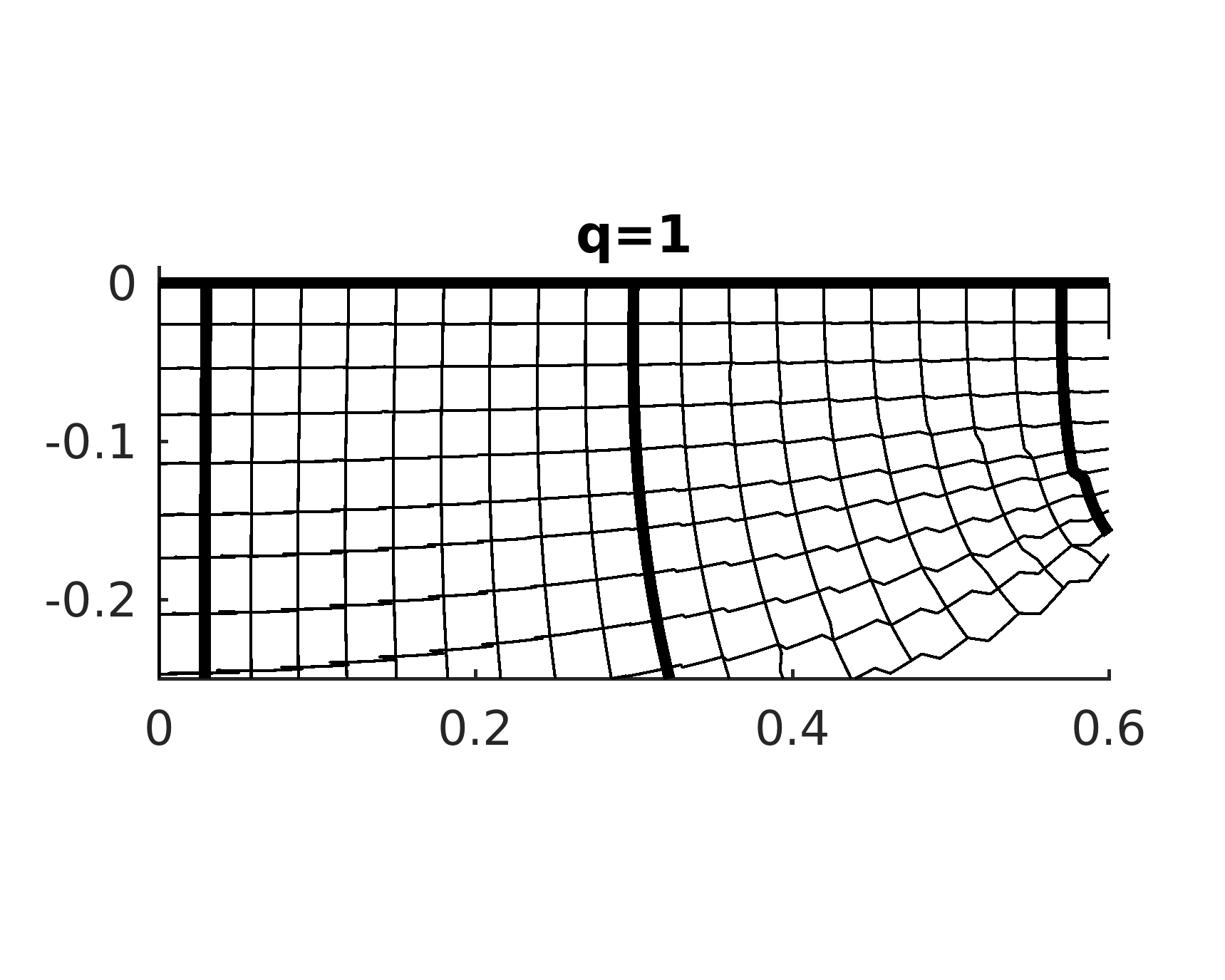}
    \caption{Inverse of $\taylor_1[\v\gamma]$.}
    \label{fig:taylor_map1}
  \end{subfigure}
  \begin{subfigure}{.32\textwidth}
    \centering
    \includegraphics[width=\textwidth]{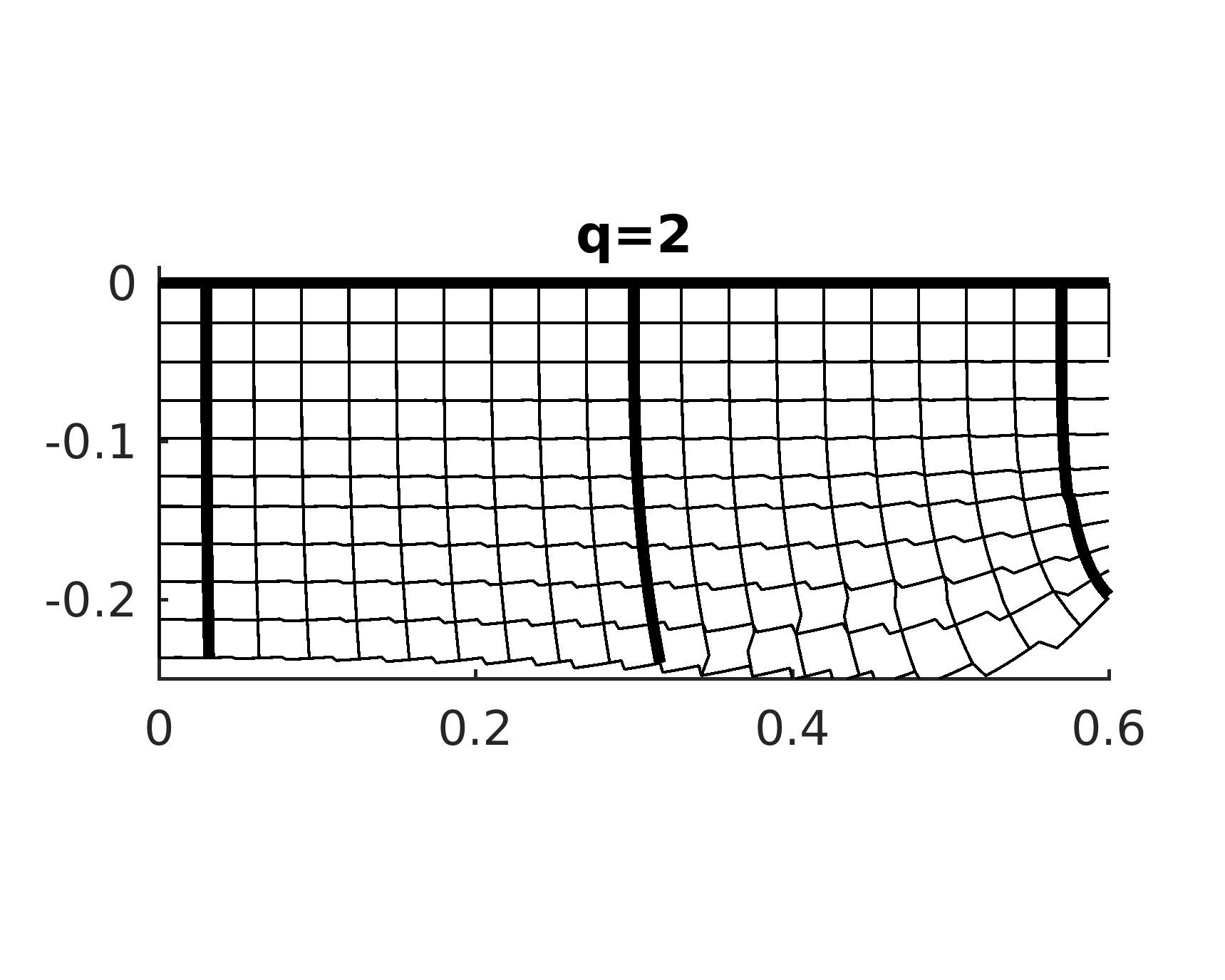}
    \caption{Inverse of $\taylor_2[\v\gamma]$.}
    \label{fig:taylor_map2}
  \end{subfigure}
  \begin{subfigure}{.32\textwidth}
    \centering
    \includegraphics[width=\textwidth]{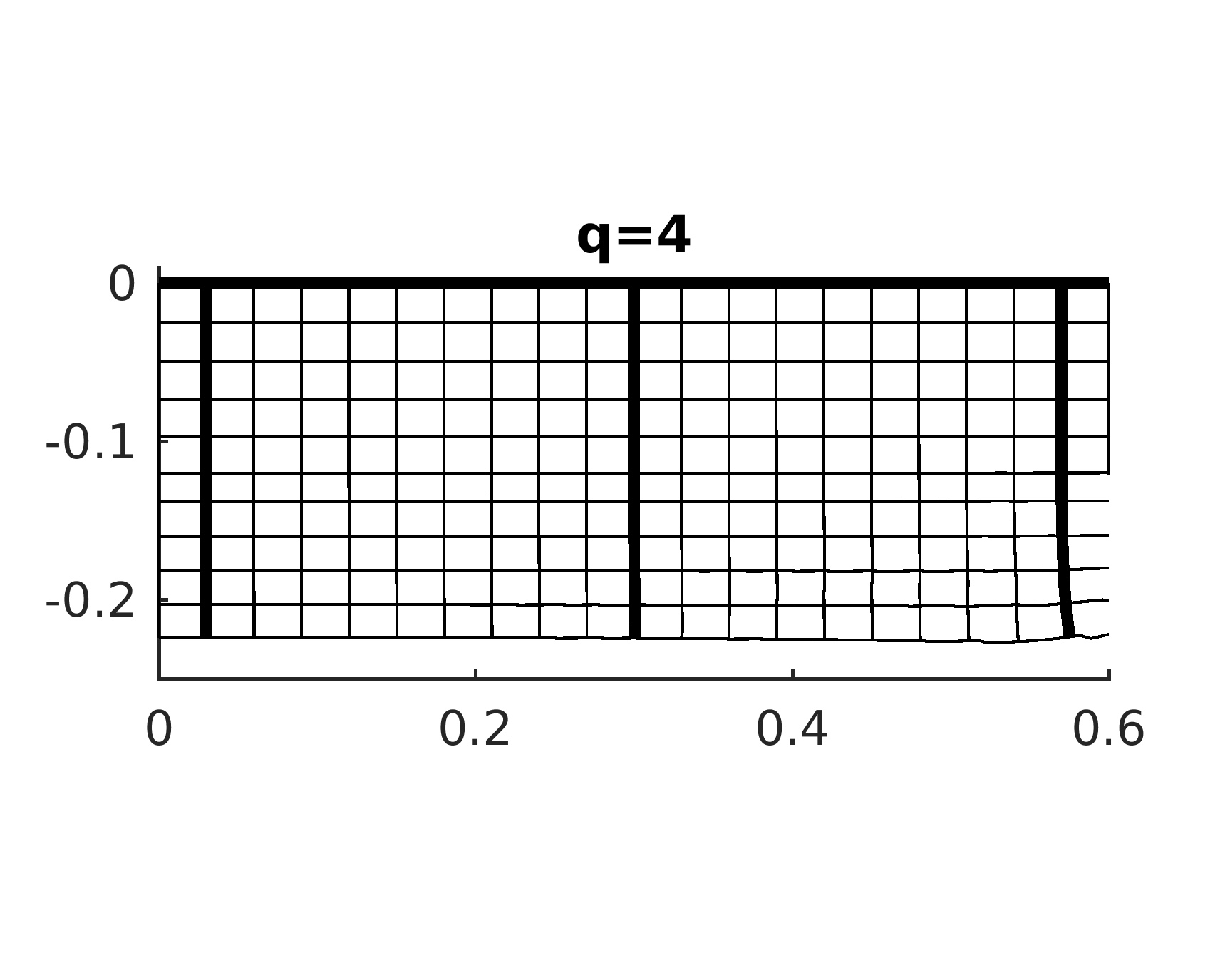}
    \caption{Inverse of $\taylor_4[\v\gamma]$.}
    \label{fig:taylor_map4}
  \end{subfigure}
  \\
  \vspace{1em}
  \begin{subfigure}{\textwidth}
    \centering
    \includegraphics[width=0.32\textwidth]{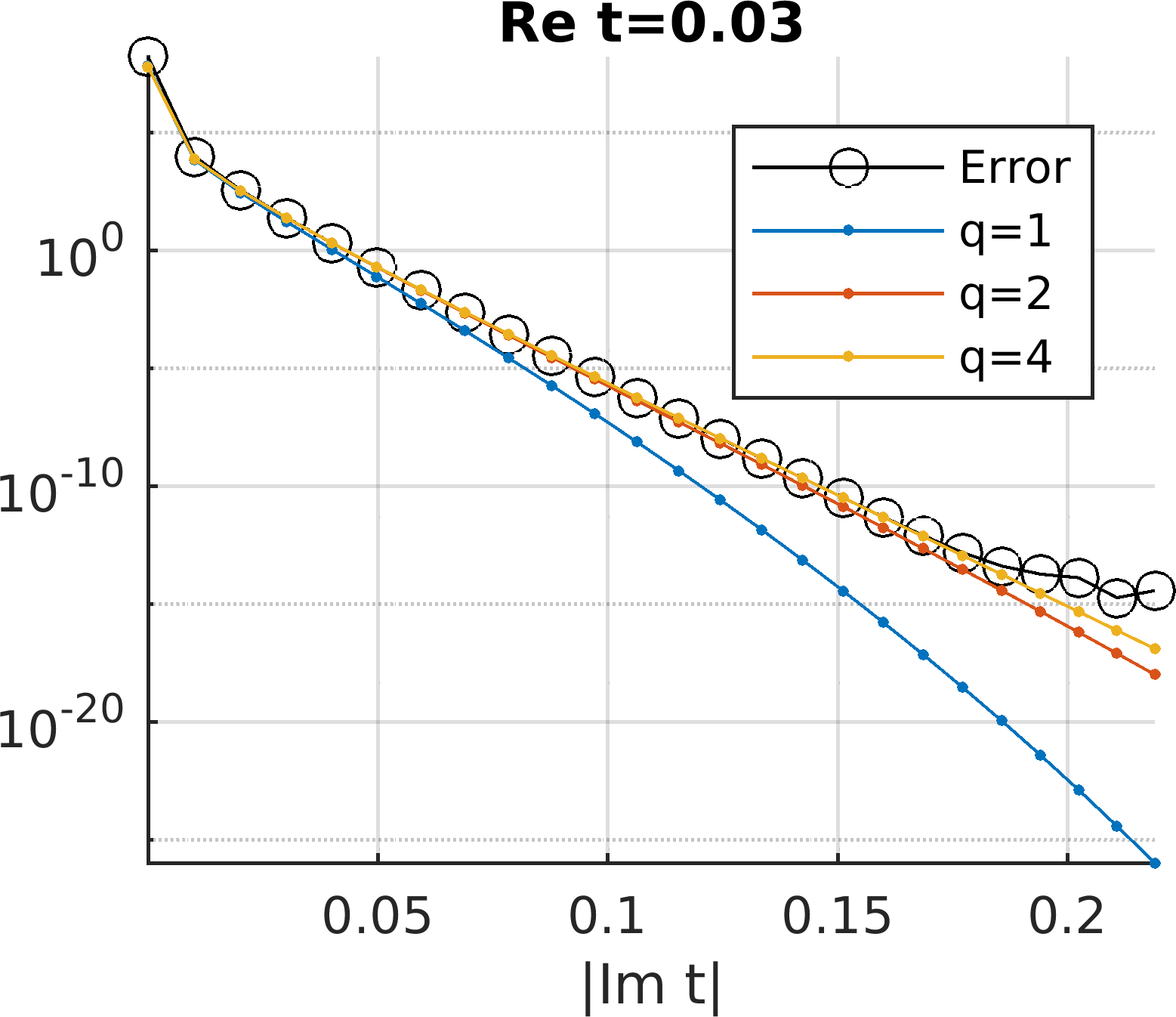}
    \includegraphics[width=0.32\textwidth]{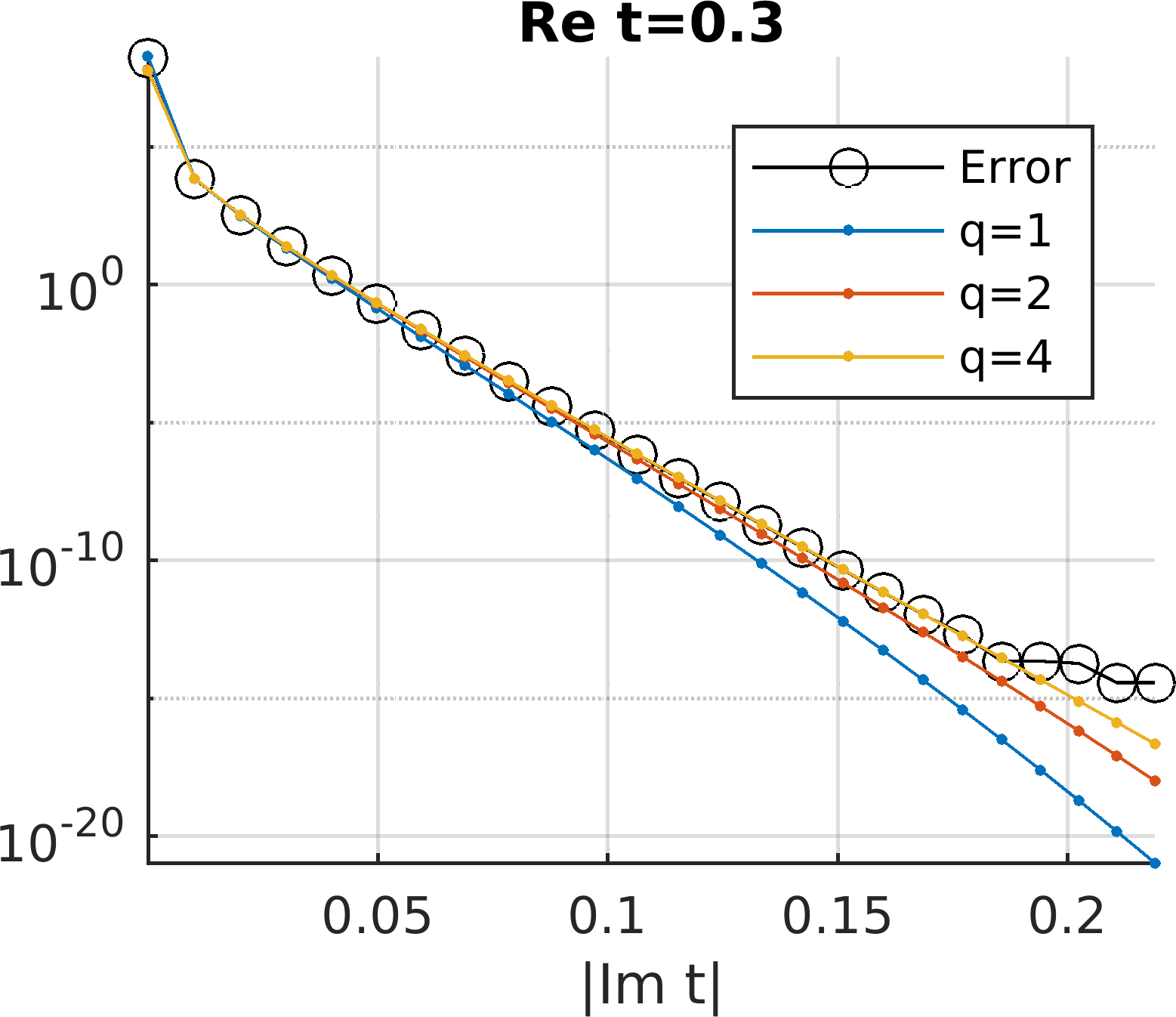}
    \includegraphics[width=0.32\textwidth]{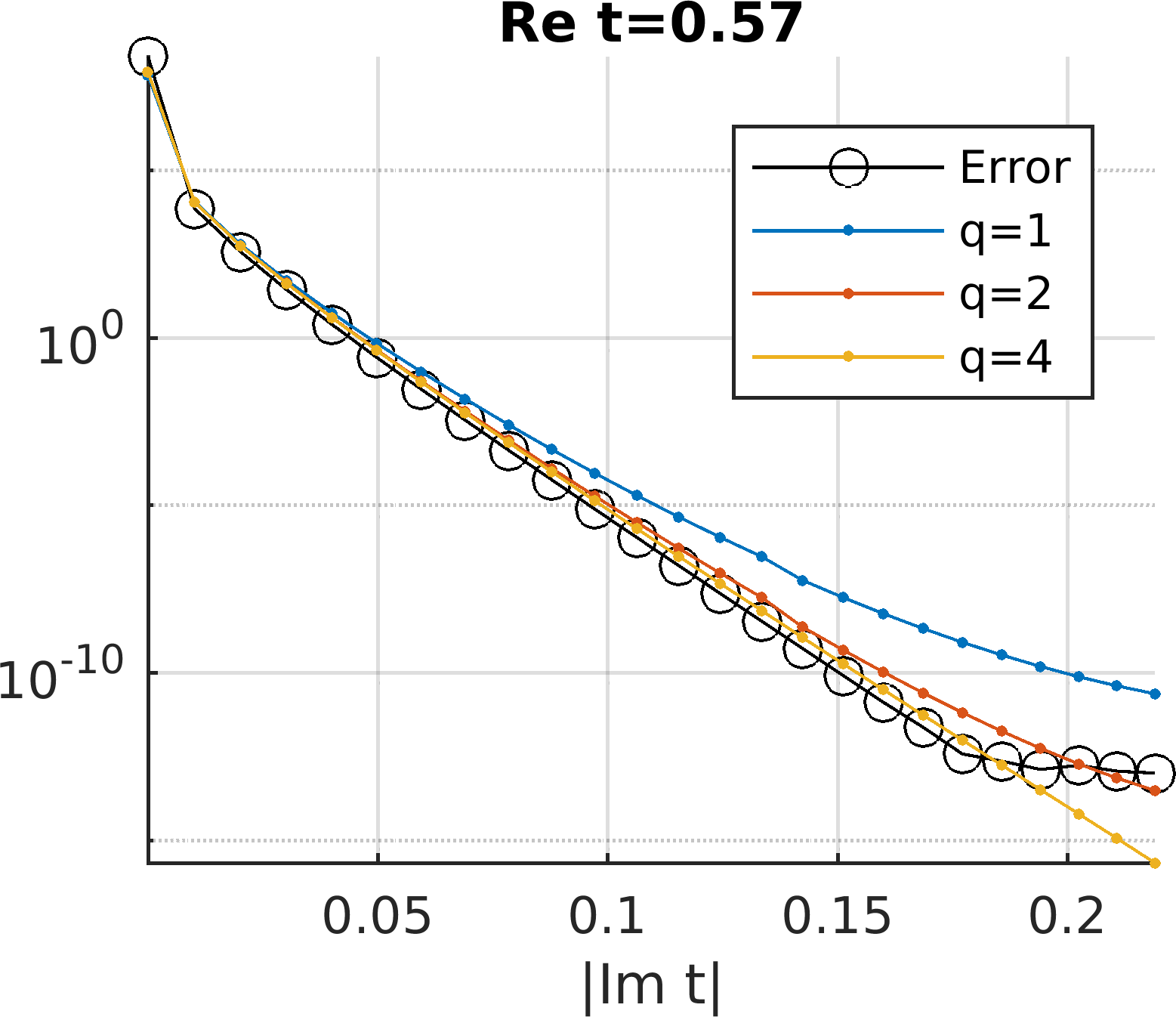}
    \caption{Errors vs estimates, computed using $q=\{1,2,4\}$, on
      lines extending out from the curve. From left to right, the
      errors are computed on lines corresponding to
      $\Re t=\{0.03, 0.3, 0.57\}$.}
    \label{fig:taylor_map_lines}
  \end{subfigure}  
  \caption{Effects of using a local Taylor expansion of order $q$ for
    finding roots based on the trapezoidal rule discretization, on the problem shown in
    \cref{fig:field_grid,fig:field_err}. \Cref{fig:taylor_map1,fig:taylor_map2,fig:taylor_map4}
    show the roots $t\in\mathbb C$ that are found when the method is
    applied to the points in the grid in \cref{fig:field_grid_field}
    (the true roots are those shown in
    \ref{fig:field_grid_param}). The ``steps'' occur where a different
    node on the curve is used as the expansion center.  The plots in
    \cref{fig:taylor_map_lines} focus on the thick black lines marked
    in the abovementioned plots, comparing the true error to the
    estimates computed using root finding of varying order. Clearly,
    $q \ge 2$ is sufficient for capturing the magnitude of the error,
    up to the point where roundoff starts to dominate.}
  \label{fig:backmaps_taylor}
\end{figure}

\Cref{fig:backmaps_taylor} revisits our trapezoidal rule example from
\cref{sec:twodim-error-contours}. This time, instead of using the
known roots to evaluate the estimates, we use roots that are computed
with the combination of Newton's method and Taylor expansion, for a
few different orders $q$. Essentially, what we are are doing is
computing the inverse of the complexification $\omega(t)$ in
\eqref{eq:complexification}. Not surprisingly, low order
approximations introduce a distortion in this map, that increases with
local curvature and distance from the curve. However, with only
moderate values of $q$ it is possible to compute the root sufficiently
well for the contours of the estimate to follow the error contours
closely, at least in this particular example.

%%\subsection{Computing $f(t_0)$ and $G(t_0)$}
\subsubsection{Evaluating quantities at the root}

So far we have derived the error estimates for one-dimensional curves,
and discussed how to find the root $t_0$ that corresponds to $\v x$.
However, in order to evaluate the estimates, we also need the values
$f(t_0)$ and $G(t_0)$. The geometry factor $G(t_0)$, defined in
\eqref{eq:geometry_factor}, can easily be computed using the
approximation $\tilde{\v\gamma}$, which we have already constructed in
order to find the root %
(it is in fact computed in the Newtons iterations).%
The value $f(t_0)$ does not come ``for free'' in the same way, and the
cheapest way to compute it is to simply bound it (approximately) using
the maximum value on the curve,
$f(t_0) \lesssim \norm{f}_{\infty(\Gamma)}$, or on a section of the
curve. Using the maximum value over the Gauss-Legendre panel
generally works well, see comparison in
\cite[Fig.~4]{AfKlinteberg2018}. Slightly more accurate, and slightly
more costly, is to construct an approximation $\tilde f$, in the same
way that we constructed $\tilde{\v\gamma}$, and then evaluate
$\tilde f(t_0)$. This is the method that we use in this paper.

\subsection{Examples for one-dimensional curves in $\reals^3$}
\label{sec:lines_in_3D}

\begin{figure}[htbp!]
  % scripts/threedim_est_tests.m
  \centering
  \includegraphics[width=0.50\textwidth]{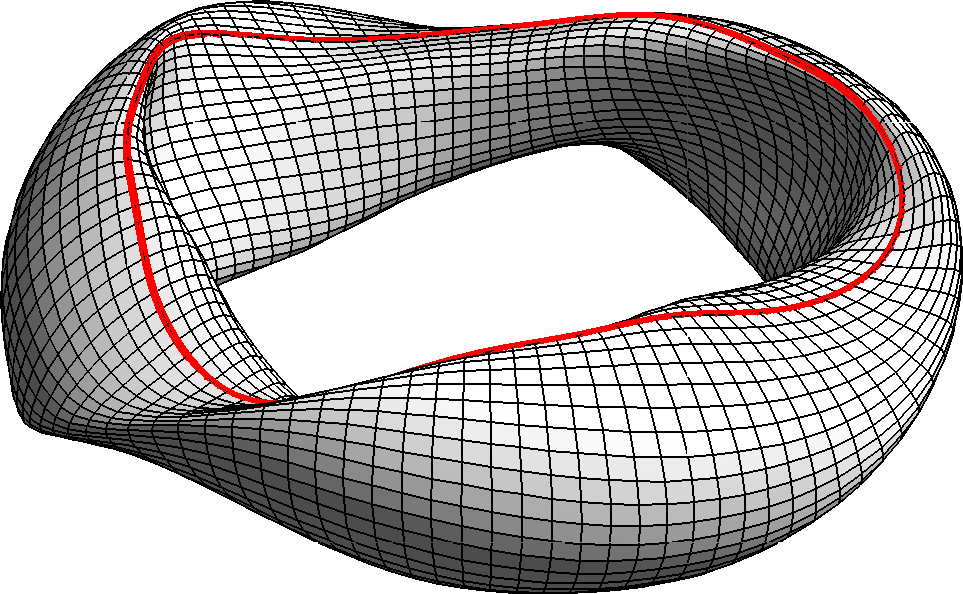}
  \caption{One-dimensional curve on the surface of the QAS3
    stellarator.}
  \label{fig:threedim_line_geo}
\end{figure}

The estimates that we have derived so far are for quadrature errors in
layer potentials near one-dimensional curves. The most obvious use for
these are in the context of boundary integral methods for planar
geometries, but there is in fact nothing limiting them to 2D problems,
as long as the source geometry is one-dimensional. In 3D,
one-dimensional source geometries appear in slender-body
approximations of fluid flow or electrical fields (see discussion in
\cite{AfKlinteberg2020line}). To demonstrate the application of our
estimates on a curve in 3D, we consider a curve (shown in
\cref{fig:threedim_line_geo}) defined on the surface known as the QAS3
stellarator \cite{Garabedian2002}. This surface was used as an example
for the integral equation solver developed in \cite{Malhotra2019}, and
we will use it for our surface error estimates in
\cref{sec:surface_estimates}.

\begin{figure}[htbp!]
  \centering
  % scripts/threedim_est_tests.m
  \begin{subfigure}{0.4\textwidth}
    \includegraphics[align=c,width=\textwidth]{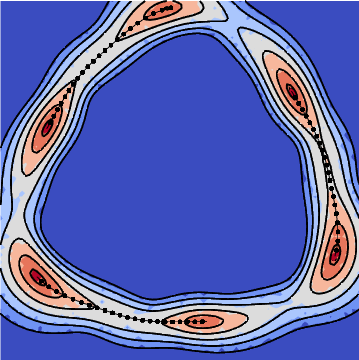}
    \caption{Trapezoidal rule.}
    \label{fig:threedim_line_field_trapz}
  \end{subfigure}
\hfill
\begin{subfigure}{0.4\textwidth}
  \includegraphics[align=c,width=\textwidth]{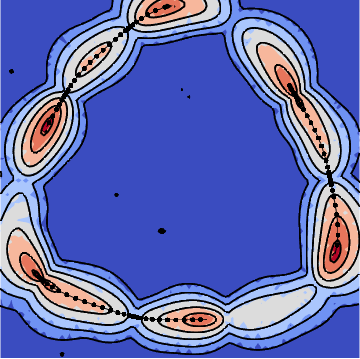}
  \caption{Gauss-Legendre rule.}
     \label{fig:threedim_line_field_gl}
\end{subfigure}
\hfill
\includegraphics[align=c,height=0.45\textwidth]{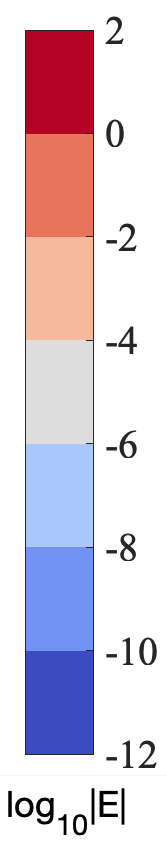}
\caption{Quadrature errors (filled contours) vs error estimates (black
  contours) for the 3D harmonic single layer potential evaluated from the
  source line marked in \cref{fig:threedim_line_geo}.
Errors are measured in the $xy$-plane for a fixed $z=1.16$.} 
  \label{fig:threedim_line_field}
\end{figure}

On the QAS3 stellerator, we define a curve $\Gamma$ by fixing the
poloidal\footnote{See \cite[Fig. 2]{Malhotra2019} for illustration of
  the toroidal/poloidal directions.} angle at $\phi=\pi/2$. The curve
is then parametrized in the toroidal angle $\theta \in [0,2\pi)$. We
set our layer potential to be the 3D harmonic single layer potential
$u(\v x)=S^{3D}_H[\sigma](\v x)$, given in \eqref{eq:3D_harmonic_sgl},
with the simple density $\sigma(\v y)=y_1 y_3$. We discretize $\Gamma$
in two ways: using the trapezoidal rule with $n=120$, and using the
composite Gauss-Legendre rule with 10 panels and $n=16$. Using these
discretizations, we evaluate $u(\v x)$ on the plane $z=1.16$ (the mean
$z$-coordinate of $\Gamma$). We compute the error using adaptive
quadrature, and estimate the error using \eqref{eq:trapz_p_halfint}
and \eqref{eq:gl_p_halfint} with $p=1/2$. The root $t_0$ is for each
target point found by approximating $\gamma$ using a 5th order Taylor
expansion in the case of the trapezoidal rule, and a 16th order
Legendre polynomial on each panel in the case of the Gauss-Legendre
rule. The results, shown in \cref{fig:threedim_line_field}, indicate
that our estimates work well also for one-dimensional curves in
3D. {The black spots in the blue area in \cref{fig:threedim_line_field_gl} is due to the fact that the root finding process has not converged to the correct root for these evaluation points. This is not of practical concern, since these locations are far from the curve. See also the discussion in connection to \cref{fig:stellerator_gl_cut_error}.}

\section{Quadrature errors near two-dimensional surfaces in $\reals^3$}
\label{sec:surface_estimates}

Let us now consider the three-dimensional case, for which our
prototype layer potential (\ref{eq:layerpot_S2}) takes the form
(\ref{eq:3Dpot}). Here, $S \subset \mathbb R^3$ is a two-dimensional
surface parametrized by $\v\gamma : E \to \mathbb R^3$,
$E = \{E_1 \times E_2\} \subset \mathbb R^2$.

Our goal is now to find a way of estimating the error committed when
\eqref{eq:3Dpot} is evaluated using an $n_s \times n_t$ tensor product
quadrature rule, based on the trapezoidal and/or Gauss Legendre quadrature
rule. The base intervals $E_1$ and $E_2$ are set according to which
rule is considered in the two directions. The use of the
Gauss-Legendre rule means that we are considering one panel that only
covers part of a full surface. To obtain the full error estimate, the
error contribution from different panels must be added. Due to the
localized nature of the errors, in practice, only the panels closest
to the target point $\v x$ need be considered. 

\subsection{Error estimates for surfaces}
If we introduce the convenience notation
$g_p=f/\norm{\v\gamma-\v x}^{2p}$, then we can write the tensor product
quadrature as
\begin{align}
  \opquads \opquadt g_p
  &= \pars{\opint_s - \oprems}\pars{\opint_t - \opremt} g_p.
\end{align}
The operators $I[g]$, $Q_n[g]$ and $E_n[g]$ where introduced in the beginning of section
\ref{sec:basictheory}. Here we use them with a subindex indicating if
they are applied in the $s$ or the $t$ direction. For ease of
notation, we have also skipped the brackets above, such that $I_s I_t g$ means an
integration of $g$ first in the $t$ and then in the $s$
direction. The full expression of the term $\opint_s \opremt g_p$ can be
found in (\ref{eq:IsEt}) below.

Neglecting the quadratic error term, and using that $\oprems\opint_t=\opint_t\oprems$, we can approximate the tensor product quadrature error as
\begin{align}
  \opremst g_p &\coloneqq 
             \left(  \opint_s\opint_t - \opquads \opquadt \right) \, g_p
             \approx
             \left( \opint_s \opremt + \opint_t\oprems \right) \, g_p
             .
             \label{eq:proderr_approx}
\end{align}
Elliott et al. \cite{Elliott2015} has shown for some basic integrals,
that the {\em remainder of the remainder} term that we here neglect can
have an important contribution. It is however a higher order
contribution, and this is only true when the quadrature error is
large, and we will proceed without it. 

\begin{figure}
% scripts/stellarator_line.m
  \centering
  \fbox{\includegraphics[width=0.65\textwidth]{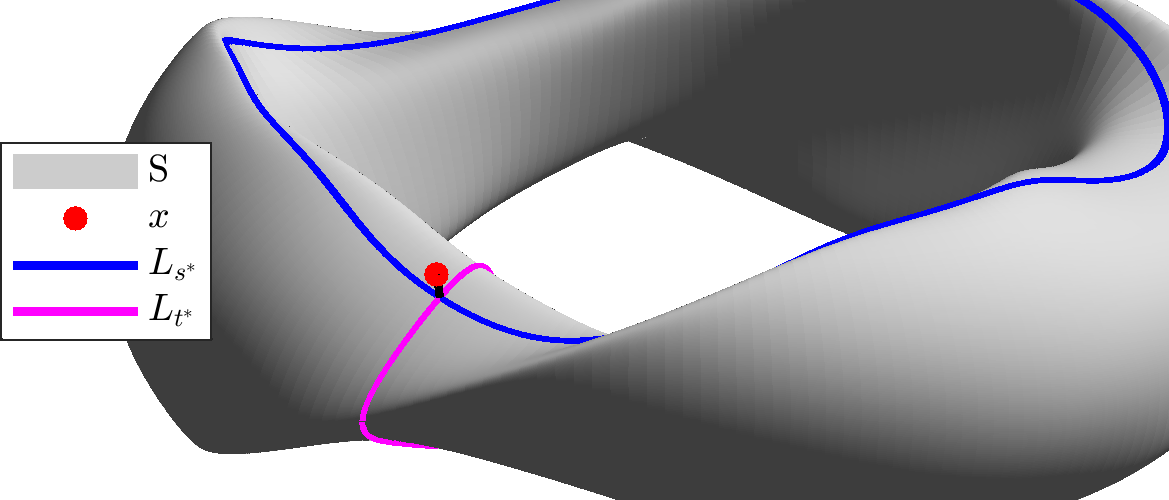}}
  \caption{Lines $L_{s^*}$ and $L_{t^*}$ on the QAS3 stellarator,
    intersecting at the point $\v{\gamma^*}=\v\gamma(s^*,t^*)$. Here
    $s$ is the poloidal angle and $t$ is the toroidal angle.}
  \label{fig:stell_lines}
\end{figure}

In essence, the formula above means that we can compute an approximation to
the tensor product quadrature error by integrating the one-dimensional
error estimates that we have already derived. To expand this
statement, let us now focus on the first term of
\eqref{eq:proderr_approx} (the second one is treated identically). We
have that
\begin{align}
  \opint_s \opremt g_p
  %%\frac{f}{\norm{\v\gamma-\v x}^{2p}}
  =
  \int_{E_1}
  \left[
  \int_{E_2}
  \frac{f(s,t) \dif t}{\norm{\v\gamma(s,t)-\v x}^{2p}}
  -
  \sum_{l=1}^n 
  \frac{f(s,t_l) w_l}{\norm{\v\gamma(s,t_l)-\v x}^{2p}}
  \right] \dif s .
  \label{eq:IsEt}
\end{align}
The term in the brackets (i.e. $\opremt$) represents the quadrature
error on the line $L_s$ that for a given $s$ is defined as
\begin{align}
  L_s \coloneqq
  \left\{
  \v\gamma(s,t)
  \mid
  t \in E_2
  \right\}
  ,\quad
  s \in E_1
  .
\end{align}
See illustration in \cref{fig:stell_lines}. %
Estimating this error is precisely the problem that was treated in
\cref{sec:qerr1Dcurves}. As we have seen, the magnitude of the error
depends on the closest root to the squared distance function, here
defined as
\begin{align}
  R^2(s,t) \coloneqq \norm{\v\gamma(s,t)-\v x}^{2} 
  \label{eq:R2_def_surface} .
\end{align}
For a given $s$, we denote by $t_0(s)$ the complex root such that
\begin{align}
  R^2\left(s, t_0(s)\right) = 0.
\end{align}
In order to abbreviate our notation, let $\est(t_0,n,p)$ denote the
quadrature rule specific part of one of the estimates derived in
\cref{sec:qerr1Dcurves}, such that
\begin{align}
 | \opint_s \opremt g_p| \le \opint_s  |\opremt g_p| \approx 
  \int_{E_1} 
  \abs{f\pars{s,t_0(s)} G\pars{s,t_0(s) }^p}
  \est\pars{t_0(s),n,p} \dif s .
  \label{eq:estintegral}
\end{align}
Since we are considering problems in three dimensions we assume that
$p$ is a half-integer. Then we have from
\cref{est:trapz_p_halfint,est:gl_p_halfint}
(\cref{eq:trapz_p_halfint,eq:gl_p_halfint}) that
\begin{align}
  \est(t_0,n,p) =
  \frac{4\pi}{ \Gamma(p)}
  \begin{cases}
    n^{p-1}  e^{-n|\Im t_0|}
    & \quad \text{for trapezoidal,}
    \\
    \abs{ \frac{2n+1}{\sqrt{t_0^2-1}} }^{p-1}
    \rho(t_0)^{-(2n+1)}
    & \quad \text{for Gauss-Legendre.}
  \end{cases}
      \label{eq:est_def}
\end{align}
Our task now is to evaluate the integral \eqref{eq:estintegral} (and
the analogous integral to estimate $|\opint_t \oprems g_p|$).  The main
difficulty in doing this is that even though we have a closed form
expression for $\est(t_0,n,p)$, we do not have one for $\est(t_0(s),n,p)$,
since $t_0(s)$ is computed numerically using the technique outlined in
\ref{sec:curve-rootfinding}. We could still evaluate
\eqref{eq:estintegral} using quadrature, but that would require us to
repeat the numerical root finding procedure multiple times for a single
target point $\v x$, something which we deem would be too costly for
the purpose of error estimation. Instead, we use the following
semi-analytical approach.

%%\subsection{Numerical evaluation of estimates on surfaces}

\subsubsection{Best approximation}
In order to evaluate \eqref{eq:estintegral}, our first step is to find
the closest grid point on the surface, which we denote
$\v{\gamma^*} = \v\gamma(s^*,t^*)$. At this location, we assume that we
have access to the derivatives $\partial_t \v\gamma$
through $\partial_t^q \v\gamma$ (either analytically, or
computed numerically at the time of discretization). Then we can form
the univariate local approximation
\begin{align}
  \v{\tilde\gamma}(s^*,t) = 
  \sum_{j=0}^q   
  \frac{ \pars{t-t^*}^j }{j!}
  \dpd[j]{\v{\gamma}}{t} (s^*,t^*)
  .
  \label{eq:gamma_taylor_3d}
\end{align}
Alternatively, we can form a global $n_t$th order polynomial
approximation on the line $L_{s^*}$,
\begin{align}
  \v{\tilde\gamma}(s^*,t) = \poly_n
  \left[\v\gamma(s^*, \cdot)\right](t).
  \label{eq:gamma_poly_3d}
\end{align}
Typically we use the latter approximation for Gauss-Legendre, since it
is ``global'' only over a panel where $n$ is small, while we utilize the
local approximation for the trapezoidal rule. 
%\todo[inline]{Mention again that global works for GL because it is
%  global \emph{on the panel}.} %
We insert this (local or global) approximation into the squared
distance function \eqref{eq:R2_def_surface} and apply the techniques
discussed in \cref{sec:curve-rootfinding} to find the root
\begin{align}
  t_0^* \approx t_0(s^*).
\end{align}
This root represents our \emph{best approximation} of $t_0(s^*)${, where $s^*$ is the value of the parametrization in $s$ at the quadrature node that is closest to  $\bf{x}$.}

\subsubsection{Linear approximation}

\newcommand{\dgds}{\v\gamma^*_s}
\newcommand{\dgdt}{\v\gamma^*_t}

We also form the bivariate linear approximation
\begin{align}
  \v\gtilde_L(s, t) &= \v\gamma^*
                    + \partial_s\v\gamma(s^*,t^*)\Delta s
                    + \partial_t\v\gamma(s^*,t^*)\Delta t,
\end{align}
where $\Delta t=t-t^*$ and $\Delta s=s-s^*$. For brevity we write $\v r=\v\gamma^*-\v x$, $\dgds=\partial_s\v\gamma(s^*,t^*)$
and $\dgdt=\partial_t\v\gamma(s^*,t^*)$.
The squared distance function \eqref{eq:R2_def_surface} then takes the
form
\begin{align}
  R^2 \approx 
  \underbrace{
  \norm{\v r}^2 
  + 2 (\v r \cdot \dgds)\Delta s 
  + \norm{\dgds}^2 \Delta s^2 
  }_{a(\Delta s)}
  + 
  \underbrace{
  \pars{
  2 (\v r \cdot \dgdt)
  + 2 (\dgds\cdot\dgdt)\Delta s
  }
  }_{b(\Delta s)}
  \Delta t
  +
  \underbrace{
  \norm{\dgdt}^2
  }_{c(\Delta s)}
  \Delta t^2 .
  \label{eq:R2_linear}
\end{align}
Finding the roots of this by solving for $\Delta t$, we get
\begin{align}
  t_0^L(\Delta s) = t^* - \frac{b}{2c} 
  \pm i \frac{\sqrt{4ac-b^2}}{2c} .
  \label{eq:t0_linear}
\end{align}
This is our \emph{linear approximation} to $t_0(s)$.

\subsubsection{Combined approximation}
\label{sec:comb-appr}

Due to the exponential distance dependence of quadrature errors, we
expect $\est(t_0(s),n,p)$ to have a peak close to $s^*$, and then decay
exponentially with $|s-s^*|$. In order to capture the magnitude of
that peak as well as possible, 
while having a simple explicit dependence on $s$,
%while still observing exponential decay, 
we define the following \emph{combined
  approximation}:
\begin{align}
  \tilde t_0(s) = t_0^*  - t_0^L(0) + t_0^L(s-s^*).
  \label{eq:combined_approx}
\end{align}
Inserting this into \eqref{eq:estintegral}, and reasoning that
$\est(t_0(s),n,p)$ is the most rapidly varying factor (with a peak near
$s^*$),
\begin{align}
  \opint_s  |\opremt g_p| \approx 
  \abs{f\pars{s^*, t_0^*} G\pars{s^*, t_0^* }^p}
  \int_{E_1}  
  \est\pars{t_0^* - t_0^L(0) + t_0^L(s-s^*),n,p} \dif s .
  \label{eq:combined_integral}
\end{align}
We are now left with a definite integral of a closed-form
function. Since we only need to compute it to 1--2 digits of accuracy,
it can be rapidly evaluated using quadrature, the details of which
depend on which estimate we are integrating, as outlined below. This completes our method for estimating quadrature errors in 3D. In \ref{sec:summ-algor-surf} we summarize the required steps, and in \ref{sec:numer-exper} we demonstrate its performance.

\begin{figure}
% scripts/stellarator_line.m
  \centering
  \begin{subfigure}[t]{.45\textwidth}
    \includegraphics[width=\textwidth]{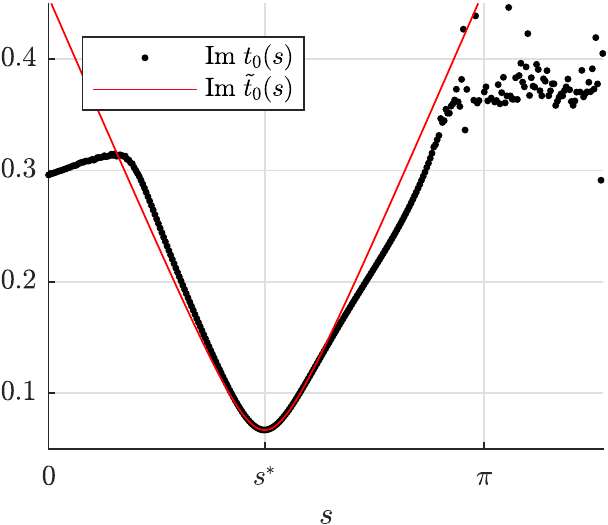}
    \caption{The root $t_0(s)$ is computed numerically on the line
      $L_s$, while $\tilde t_0(s)$ is computed using the combined
      approximation \eqref{eq:combined_approx}.}
    \label{fig:roots_of_s}
  \end{subfigure}
  \hfill
  \begin{subfigure}[t]{.45\textwidth}
    \includegraphics[width=\textwidth]{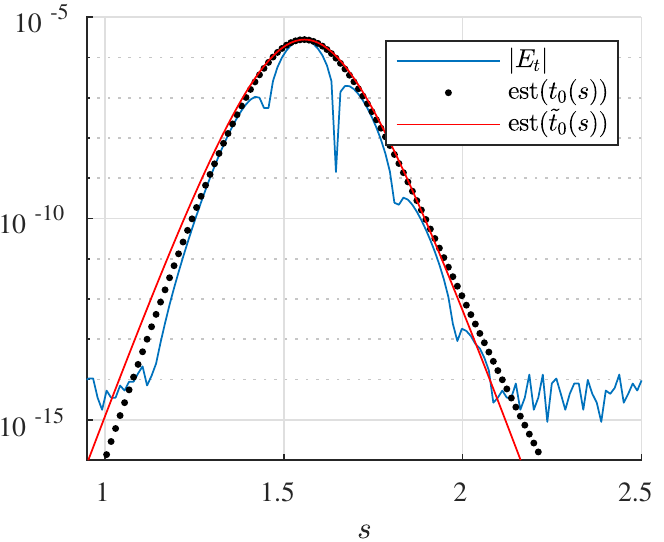}
    \caption{Actual trapezoidal rule error $|\opremt|$ on the line
      $L_s$, compared to estimates computed using the roots in
      \cref{fig:roots_of_s}.}
    \label{fig:ests_of_s}
  \end{subfigure}
  \caption{Quantities measured on the line marked $L_{t^*}$ in
    \cref{fig:stell_lines}. Kernel is single layer.}
  \label{fig:roots}
\end{figure}

\paragraph{Trapezoidal rule}%

In the case of the trapezoidal rule we are integrating the estimate
\eqref{eq:trapz_p_halfint} on the periodic interval $E_1 = [0,2\pi)$.
The linear approximation does not take the periodicity into account,
so it is reasonable to use $\Delta s \in [-\pi,\pi]$ as the interval
of integration. On this interval the estimate decays several orders of
magnitude, since it loops around the entire geometry in physical
space. Quantifying this decay, we have that \eqref{eq:trapz_p_halfint}
decays as $e^{-n~|\Im t_0|}$. For large $|s-s^*|$ the imaginary part
of our linear approximation grows as
\begin{align}
  \abs{\Im t_0(s)} \sim k\abs{s-s^*}, \quad k=\frac{\norm{\dgds}}{\norm{\dgdt}},
\end{align}
so asymptotically the estimate decays as (temporarily omitting $n,p$
in the argument of $\est(.)$),  
\begin{align}
  \est\pars{t_0(s)} \sim e^{-nk~|s-s^*|} .
  \label{eq:trapz_decay}
\end{align}
For our purposes (1--2 digits of accuracy) we can safely expand the
interval of integration in \eqref{eq:combined_integral} from
$[-\pi,\pi]$ to $[-\infty,\infty]$, as the added tails are negligible,
\begin{align}
  \int_{E_1}  
  \est\pars{\tilde t_0(s)} \dif s
  \approx
  \int_0^\infty
  \est\pars{\tilde t_0(s^*+\Delta s)} \dif\,(\Delta s)
    +
  \int_0^\infty
    \est\pars{\tilde t_0(s^*-\Delta s)} \dif\,(\Delta s) .
\end{align}
Then, Gauss-Laguerre quadrature is suitable, as it is a Gaussian
quadrature rule for integrals of the type
$\int_0^\infty g(x) e^{-x} \dif x$ \cite[\S3.5(v)]{NIST:DLMF}.
Substituting $x=nk\Delta s$,
\begin{align}
  \int_0^\infty
  \est\pars{\tilde t_0(s^*\pm\Delta s)} \dif\,(\Delta s)  
  &=
    \frac{1}{nk}
  \int_0^\infty
    \
    \underbrace{
    \left[
  \est\pars{\tilde t_0(s^*\pm x/nk)}
    e^{x}
    \right]
    }_{h^{\pm}(x)}
    e^{-x}
  \dif x .
\end{align}
We find that it is sufficient to apply 8-point Gauss-Laguerre
quadrature to~$h^{\pm}(x)$.

\paragraph{Gauss-Legendre rule}

The integration of in \eqref{eq:combined_integral} is more
straightforward in the case of the Gauss-Legendre estimate
\eqref{eq:gl_p_halfint}, where the interval $E_1=[-1,1]$ runs between
the edges of a panel in the neighborhood of the target point $\v x$.
Here we have found that it is sufficient to use Gauss-Legendre
quadrature with 8 points to evaluate
\eqref{eq:combined_integral}. Depending on the location of $s^*$, this
is done using either two 4-point rules or one 8-point rule,
\begin{align}
  \int_{-1}^1 = 
  \begin{cases}
    \int_{-1}^{s^*} + \int_{s^*}^{1} & \text{ if } |s^*| < 0.9,\\
    \int_{-1}^1 & \text{ otherwise }.
  \end{cases}
\end{align}
We remark that the second case above is applicable to the case also where
$s^*<-1$ and $>1$, which would typically occur when the target point
$\v x$ is closest to a neighboring panel. 

\subsection{Summary of algorithm for error estimation near surfaces}
\label{sec:summ-algor-surf}

We now summarize our algorithm for quadrature error estimation near
surfaces: Given a layer potential of the form \eqref{eq:3Dpot} with a
half-integer $p$, evaluated at a target point $\v x$, the quadrature
error due to a near singularity in the integrand can be accurately
estimated through the following steps:
\begin{enumerate}
\item Identify the grid point $\v\gamma^* = \v\gamma(s^*, t^*) $ on
  $S$ that is closest to $\v x$, {see eq. (\ref{eq:t_star})}.  
\item Form $\v{\tilde\gamma}(s^*,t)$ using either a local
  \eqref{eq:gamma_taylor_3d} or a global \eqref{eq:gamma_poly_3d}
  approximation.
  \label{item:secondstep}
\item Use Newton's method to find $t_0^*$ such that
  $\norm{\v{\tilde\gamma}(s^*,t_0^*) - \v x}^2=0$, as outlined in
  \cref{sec:curve-rootfinding}. Evaluate $G(s^*, t_0^*)$ (already
  computed in the Newton iterations), as this quantity is used in the
  evaluation of the error estimate.
\item Compute an approximation to $f(s^*, t_0^*)$, either as
  $\norm{f}_{\infty(S)}$, or through $\tilde f(s^*, t)$ formed using
  the same kind of approximation used for
  $\v{\tilde\gamma}$ in step \ref{item:secondstep}.
\item Form the combined root approximation $\tilde t_0(s)$ defined in
  \eqref{eq:combined_approx}, using \eqref{eq:R2_linear} and
  \eqref{eq:t0_linear} together with $t_0^*$ and $\dgds,\dgdt$.
\item Compute $\opint_s |\opremt g_p| $ through numerical integration of
  \eqref{eq:combined_integral} as outlined in \cref{sec:comb-appr},
  with $\est(t_0)$ given by \eqref{eq:est_def}, and using the quantities
  $G(s^*, t_0^*)$, $f(s^*, t_0^*)$, and $\tilde t_0(s)$ from previous
  steps.
\label{item:laststep}
\item Compute $\opint_t |\oprems g_p|$ by repeating steps
  \ref{item:secondstep}--\ref{item:laststep} with the roles of $s$ and
  $t$ interchanged. 
\item Estimate the quadrature error at $\v x$ as
  \begin{align}
    |\opremst g_p| \approx \opint_s  |\opremt g_p| + \opint_t
    |\oprems g_p| .
  \end{align}
\item (In the case of panel-based Gauss-Legendre quadrature, repeat
  the above steps for all panels $S$ close to $\v x$ and sum the
  contributions.)
\end{enumerate}

\section{Numerical experiments for a surface in $\reals^3$}
\label{sec:numer-exper}

We will now show that the method of \cref{sec:surface_estimates} can
be used for accurately estimating the nearly singular quadrature error
when evaluating a layer potential from a surface in three
dimensions. 

{In \cref{table_kernels} we identify the corresponding kernels and $p$ in the integral \eqref{eq:layerpot_S2} for different PDEs.
%We will start by rewriting $f$ in eq. (\ref{eq:3Dpot}) as
%\begin{equation}
%f(s,t)=\phi(s,t) \norm{\dpd{\v\gamma}{s} \times \dpd{\v\gamma}{t}}
%\label{eq:f_phi}
%\end{equation}
%\begin{align}{
%  u(\v x) = \iint_E 
%  \frac{k \pars{\v x, \v \gamma(s,t)}\sigma(\v\gamma(s,t)) }
%  {\norm{\v\gamma(s,t)-\v x}^{2p}} 
%  \norm{\dpd{\v\gamma}{s} \times \dpd{\v\gamma}{t}}
%  \dif t \dif s=  \iint_E \frac{\phi(s,t)}{\norm{\v\gamma(s,t)-\v x}^{2p}} \norm{\dpd{\v\gamma}{s} \times \dpd{\v\gamma}{t}} \dif t \dif s,}
%\label{eq:3Dpot_mod}
%\end{align}
%so that one can easily identify the corresponding kernel and value of $p$ for different PDEs, as listed in table \ref{table_kernels}. 
Here we consider only scalar kernels, but the method can be directly applied also to tensorial kernels, where the estimates are applied component by component to the vectorial output. From this, the function $f(s,t)$ introduced in \eqref{eq:3Dpot} and used in the error estimates can easily be identified.}

\begin{table}[h!] 
\begin{center}
{
\begin{tabular}{|| c | c| c ||} 
 \hline
%PDE & Single layer & Double layer \\ 
\multicolumn{3}{|c|}{}\\
\multicolumn{3}{|c|}{\large $ u(\v x) = \int_{\Surf} \frac{k \pars{\v x, \v y}\sigma(\v y)
  }{\norm{\v y - \v x}^{2p}} \dif S(\v y) $}  \\ 
  \multicolumn{3}{|c|}{}\\
  \hline
  & & \\
  \textbf{PDE} & \textbf{Single layer,} $p=1/2$ & \textbf{Double layer}, $p=3/2$ \\
   & & $k\pars{\v x, \v y}=\tilde{k}\pars{\v x, \v y} {\v n_y} \cdot  (\v y - \v x)$\\
 \hline\hline
 Harmonic &  & \\
 $\Delta u=0$ & $k \pars{\v x, \v y}=1$ & $\tilde{k} \pars{\v x, \v y}=1$\\
 \hline
  Helmholtz & & \\
  $(\Delta+\omega^2)u=0$&  $k \pars{\v x, \v y}=e^{i\omega \norm{\v y - \v x}} $ & $\tilde{k} \pars{\v x, \v y}=(i\omega\norm{\v y - \v x}-1)e^{i\omega\norm{\v y - \v x}} $\\
  \hline
   Mod. Helmholtz & & \\
  $(\Delta-\omega^2)u=0$&  $k \pars{\v x, \v y}=e^{-\omega \norm{\v y - \v x}} $ & $\tilde{k}  \pars{\v x, \v y}=(-\omega\norm{\v y - \v x}+1)e^{-\omega\norm{\v y - \v x}} $\\
  \hline
\end{tabular}}
\caption{{PDEs with corresponding single and double layer kernels, identifying the functions $k$ and the values of $p$ for integrals written in the form \eqref{eq:layerpot_S2}.}}
\label{table_kernels}
\end{center}
\end{table}

\begin{figure}[htbp]
  \centering
  \hfill
  \begin{subfigure}[b]{0.48\textwidth}
    \centering
\includegraphics[align=c,width=\textwidth]{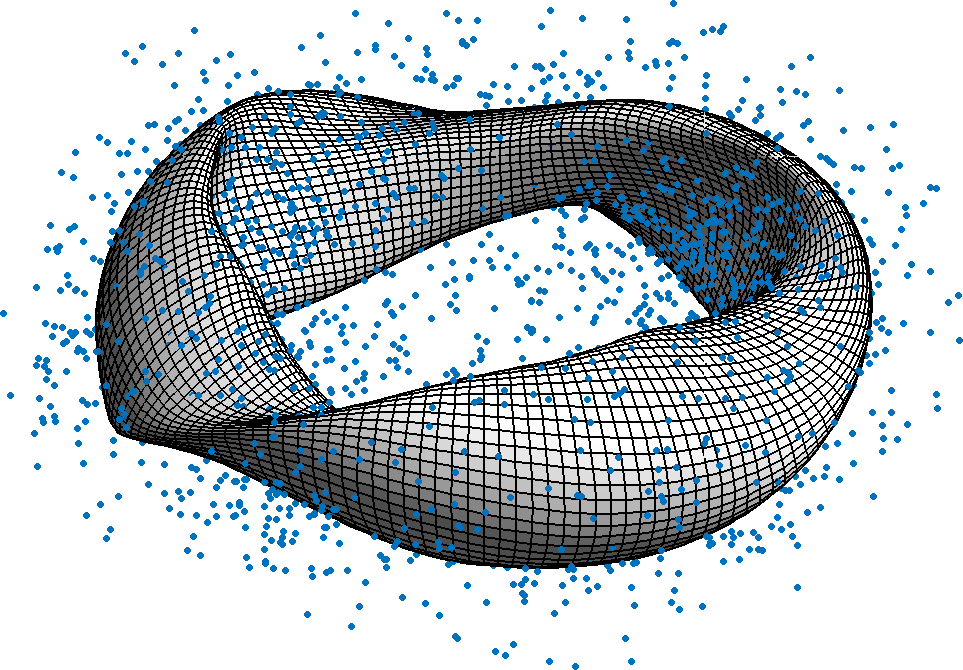}
\caption{Trapezoidal rule discretization with 3000 random test
  points.}
    \label{fig:stellerator_trapz_randpoints}
  \end{subfigure}
  \hfill
  \begin{subfigure}[b]{0.48\textwidth}
    \centering    \includegraphics[align=c,width=0.85\textwidth]{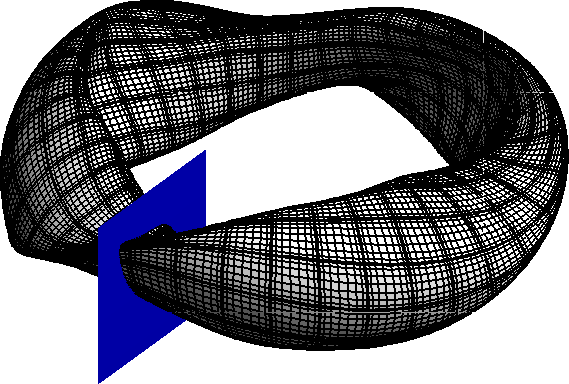}
    \vspace{1em}
    \caption{Gauss-Legendre discretization with test points on a plane
      cutting the surface.}
    \label{fig:stellerator_gl_cut}
  \end{subfigure}
  \caption{The QAS3 stellarator geometry used in our tests, showing two of the discretizations and test point distributions.}
  \label{fig:stell_grids}
\end{figure}

As our source geometry we will use the QAS 3 stellarator
\cite{Garabedian2002} that was used in \cref{sec:lines_in_3D}. We
discretize the surface in two ways: using the tensor product
trapezoidal rule with $50 \times 150$ points, and using a
$12\times 36$ grid of quadrilateral panels, each discretized with an
$8 \times 8$ tensor product Gauss-Legendre rule. These two
discretizations are shown in \cref{fig:stell_grids}. We choose to start to
evaluate the 3D harmonic double layer potential (in the form of
(\ref{eq:3Dpot}), with $p=3/2$)
\begin{align}
  D^{3D}_H[\sigma](\v x)=\int_{\Surf}\sigma(\v y) \frac{{\v n_y} \cdot
  (\v y - \v x)}{\norm{\v y - \v x}^3} \, \, \dif S(\v y),
 \label{eq:dl_harmonic_3D}
\end{align}
with the (arbitrarily chosen) density
\begin{align}
  \sigma(\v y)=\sigma(\v\gamma(s,t))=1+\cos(s)\sin(t) .
  \label{eq:sigma_experiments}
\end{align}
This density is illustrated in \cref{fig:toroidal_shell_geo}. 
We could instead have used a layer density that is the actual solution
to a corresponding integral equation, such that the layer potential
(\ref{eq:dl_harmonic_3D}) would produce the solution to a specific
Laplace boundary value problem. This is what was done in section
\ref{sec:results_complex}, for 2D results.
Our experience is however that the density does not influence the
nearly singular quadrature error much, as long as it is well resolved
by the discretization. 
In 2D, the integral equation for the layer density $\sigma$ can be
discretized by the regular Gauss-Legendre or trapezoidal rule.
In 3D however, the integrand has a singularity, and a special 
quadrature method is needed to obtain accurate results. By specifying
the density, we avoid any pollution from errors in the density as well
as building the infrastructure for the accurate solution of the
integral equation. 

In all our tests, we compute the layer potential error by comparing
with a potential computed using a grid with twice as many points in
each direction. This choice of reference is cost effective in 3D, and
sufficient for our purposes since we only need to know the error to
within 1--2 digits of accuracy. The error is then compared to the
estimate computed using the algorithm outlined in
\cref{sec:summ-algor-surf}, with $p=3/2$. Just as in
\cref{sec:lines_in_3D}, $\v{\tilde\gamma}$ is constructed using a 5th
order Taylor expansion in the case of the trapezoidal rule.

In order to illustrate how our estimates perform, we will below report
the results on several different sets of measurement points.

\subsection{Random test points}

As our first test, we compute the layer potential error at 3000 random
points located in both the interior and exterior of the stellarator
surface. The random points are generated by going a random distance in
the normal direction from a random point on the surface,
\begin{align}
  \v x = \gamma(s, t) + d \v n(s,t),
\end{align}
with $(s,t,d)$ uniformly distributed random variables,
\begin{align}
  (s,t,d) \sim ~\mathcal U(0,2\pi)
  \times \mathcal U(0,2\pi)
  \times \mathcal U(-h, h)
  \quad
  \mbox{ with } h=1.5~.
\end{align}
This results in the cloud of test points shown in \cref{fig:stellerator_trapz_randpoints}, although it must be noted that only the exterior points are visible in the figure.

\begin{figure}[htbp]
  \centering
  \begin{subfigure}{.48\textwidth}
    \includegraphics[width=\textwidth]{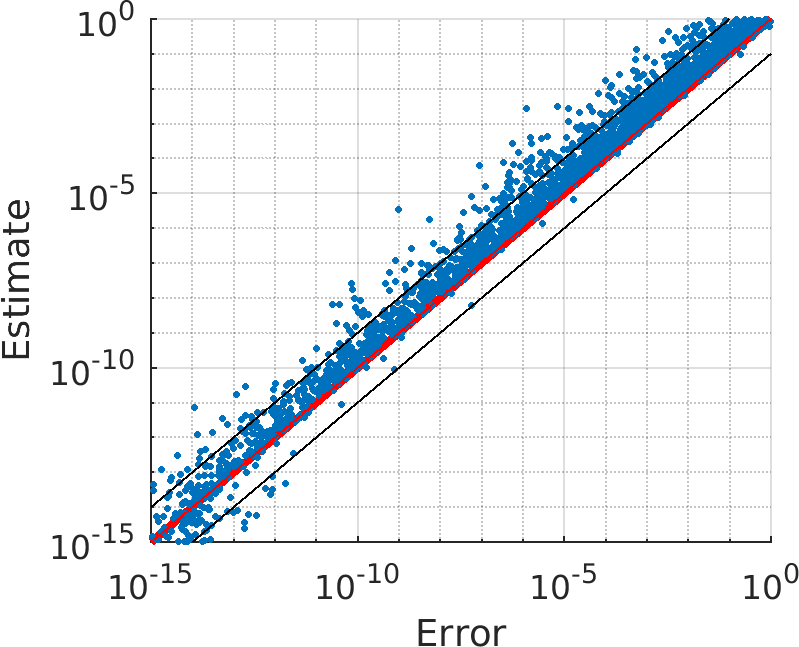}
    \caption{Trapezoidal}
  \end{subfigure}
  \hfill
  \begin{subfigure}{.48\textwidth}
    \includegraphics[width=\textwidth]{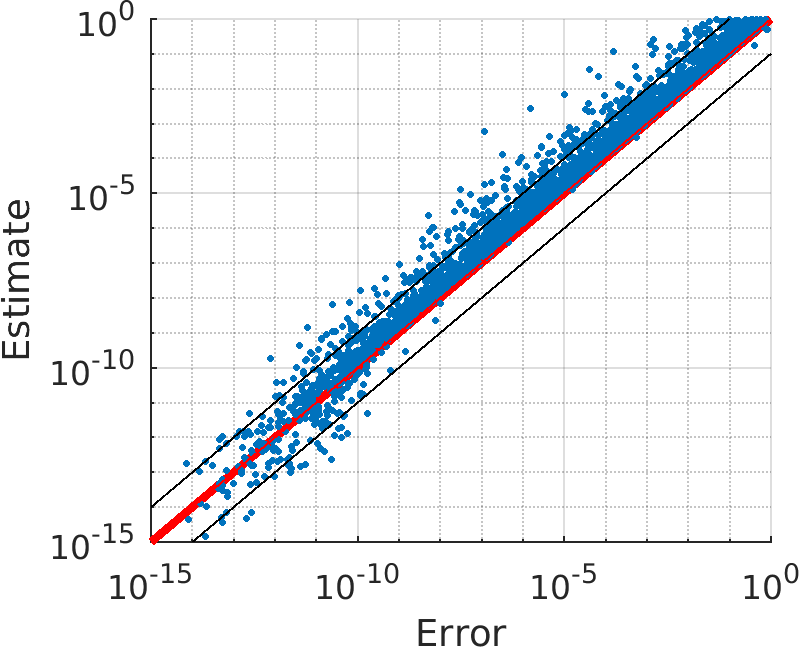}
    \caption{Gauss-Legendre}
  \end{subfigure}
  \caption{Scatter plot of error vs estimate on the random points shown in \cref{fig:stellerator_trapz_randpoints}. The red line indicates where error and estimate are equal, while the black lines indicate where they differ by factors $10$ and $1/10$, respectively.}
\label{fig:stellerator_randpoints_scatter}
\end{figure}

At each of the 3000 test points and for each of the two
discretizations, we compute the error in the layer potential and
compare it to our error estimate. This data is then used to generate
the scatter plots shown in \cref{fig:stellerator_randpoints_scatter}. 
These plots indicate that our estimates have the following important features:
\begin{itemize}
\item They are conservative most of the time (i.e. they overestimate
  the error). Only at a few points close to the surface where the
  errors are large do they slightly underestimate the error. 
\item They are within a factor 10 of the actual error most of the time.
\item Starting at errors as small as around $10^{-15}$ for trapezoidal and
  $10^{-10}$ for Gauss-Legendre, the estimates never underestimate the error by more than a factor 10.
%They never underestimate the error by more than a factor 10, as
 % long as the nearly singular error is larger than the resolution
 % error, which here is around $10^{-15}$ for trapezoidal and
 % $10^{-10}$ for Gauss-Legendre.
\end{itemize}

\subsection{Flat plane cutting the surface}
\label{sec:flat-plane-cutting}

\begin{figure}[htbp]
  \centering 
  \begin{subfigure}{0.44\textwidth}
    \includegraphics[width=\textwidth]{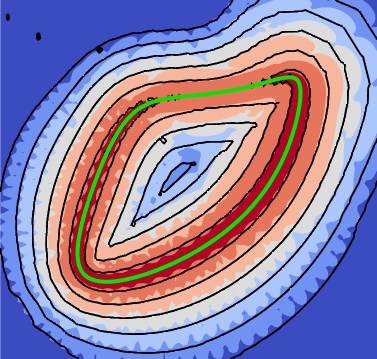}
    \caption{Trapezoidal}
    \label{fig:stellerator_trapz_cut_error}
  \end{subfigure}
  \begin{subfigure}{0.44\textwidth}
    \includegraphics[width=\textwidth]{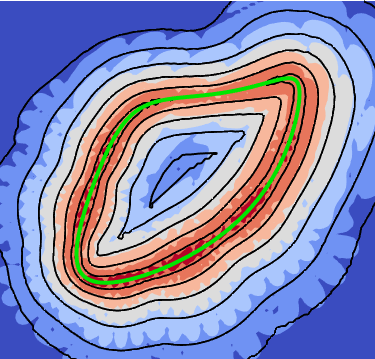}
    \caption{Gauss-Legendre}
    \label{fig:stellerator_gl_cut_error}
  \end{subfigure} \includegraphics[align=c,height=0.45\textwidth]{threedim_line_colorbar3}
  \caption{Results on the plane shown in \cref{fig:stellerator_gl_cut}
    for the two discretization, with the surface cross section marked
    green. The errors are shown as colored fields, with the contours
    of the estimates drawn in black for the levels
    $10^{\{-12,-10,-8,-6,-4,-2,0\}}$.}
  \label{fig:stellerator_cut}
\end{figure}

As our second test, we compare errors and estimate on a set of
$100 \times 100$ points covering the square plane shown in
\cref{fig:stellerator_gl_cut}. The results, displayed in
\cref{fig:stellerator_cut}, show that our estimates predict the error
levels very well close to the surface, which is where one would
normally want to use error estimates. However, the accuracy of the
Gauss-Legendre estimate (\cref{fig:stellerator_gl_cut_error})
deteriorates far away from the surface, meaning that the root finding
process has not converged to the correct root. This is likely because
the 8th order local polynomials used in the root finding get inaccurate
that far away (using higher order panels would most likely yield
better results). This is however far enough away from the surface that
error estimates would typically not be applied, it is close to the
surface that the error estimates are critical.  Note that there are also a few isolated points in
\cref{fig:stellerator_trapz_cut_error}, where the error is
overestimated. We can not explain why the root finding has failed in
these locations.
The results found here do support the conclusions made in the previous
subsection, when discussing the results in figure 
\ref{fig:stellerator_randpoints_scatter}.

\subsection{Toroidal shell}

\begin{figure}[htbp!]
  \centering
  \includegraphics[align=t,width=0.49\textwidth]{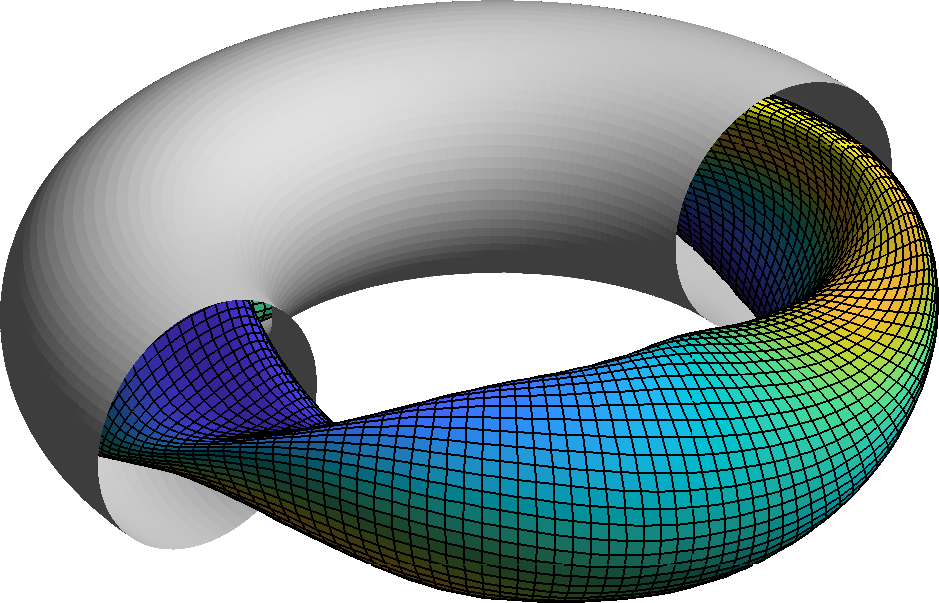}
  \includegraphics[align=t,width=0.49\textwidth]{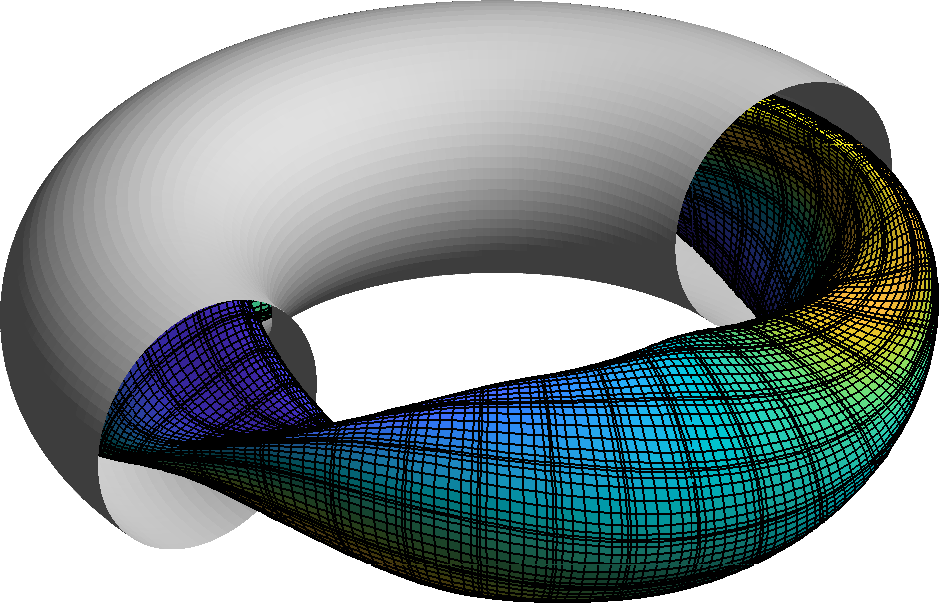}  
  \caption{Toroidal shell on which we compute quadrature errors and
    estimates (note that only half the shell is shown here). The
    coloring of the stellarator surface corresponds to the density
    $\sigma$ given by \eqref{eq:sigma_experiments}.}
  \label{fig:toroidal_shell_geo}
\end{figure}

As our final test, we let our tests points be a grid of
$200 \times 76$ points on the surface of a torus with major radius
$R=4.5$ and minor radius $r=1.7$, shown in
\cref{fig:toroidal_shell_geo}. This surface, which we denote the
\emph{toroidal shell}, encloses the stellarator from which the layer
potential is evaluated. 

\begin{figure}[htbp!]
  \centering
  \includegraphics[width=\textwidth]{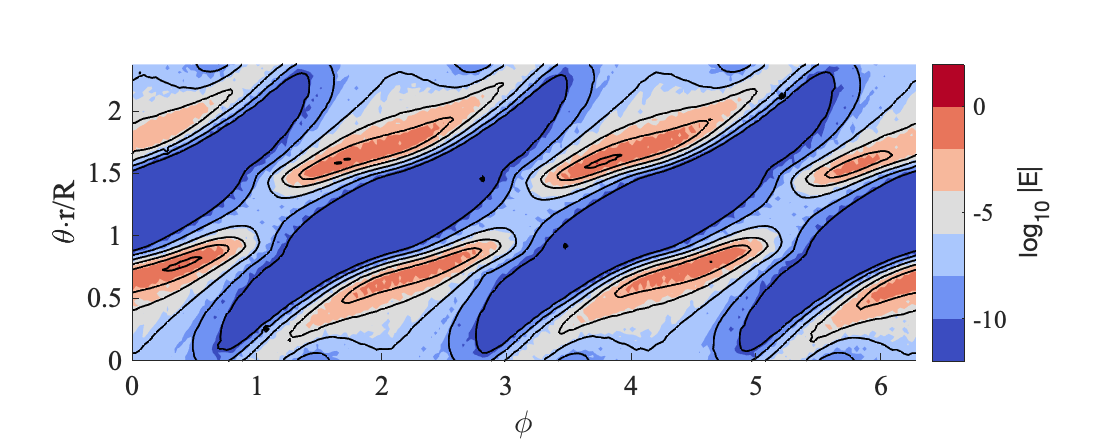}  
  \\
  \includegraphics[width=\textwidth]{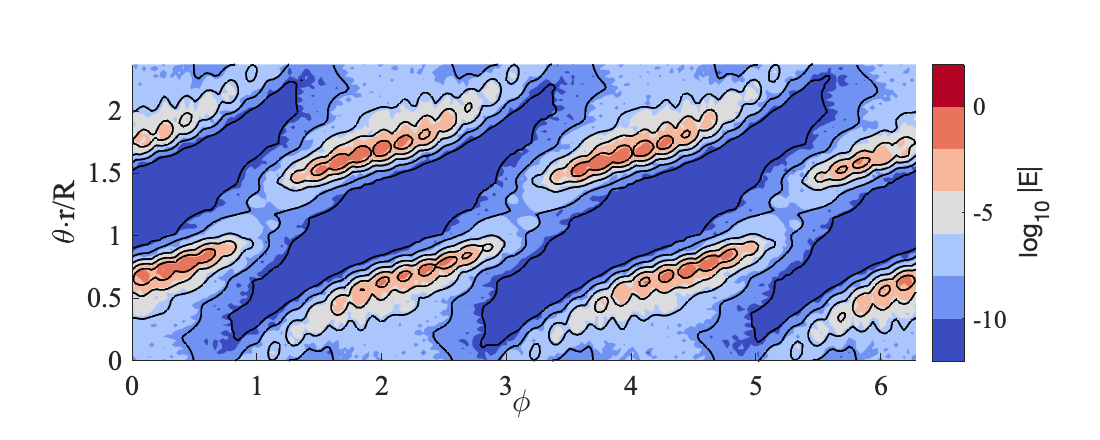}  
  \caption{Errors for the harmonic double layer potential as colored fields and estimates as black contours on
    the toroidal shell, for the trapezoidal rule (above) and the
    Gauss-Legendre rule (below). The contours are drawn for the error
    levels $10^{\{-10, -8,-6,-4,-2,0\}}$. Note that the poloidal
    angle (vertical axis) is scaled by the aspect ratio of the torus.}
  \label{fig:toroidal_shell_contours}
\end{figure}

\begin{figure}[htbp!]
  \centering
  \includegraphics[width=\textwidth]{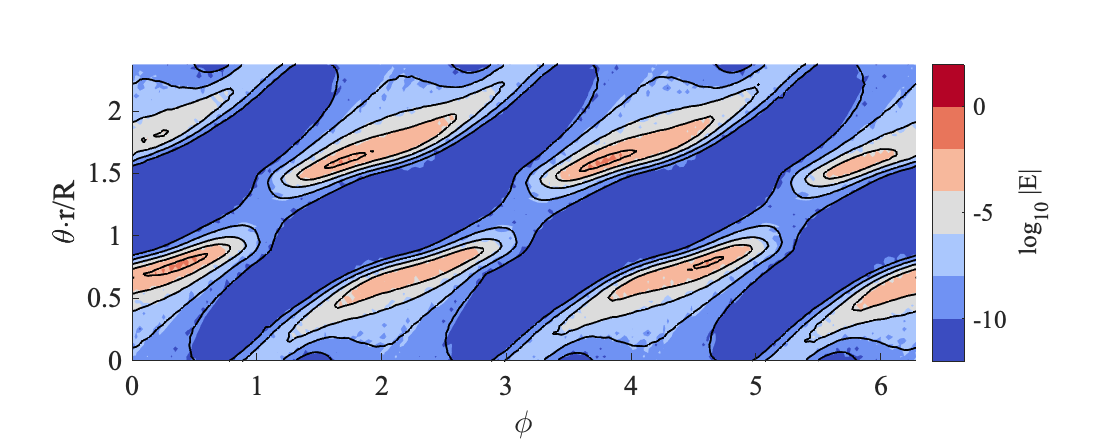}  
  \\
  \includegraphics[width=\textwidth]{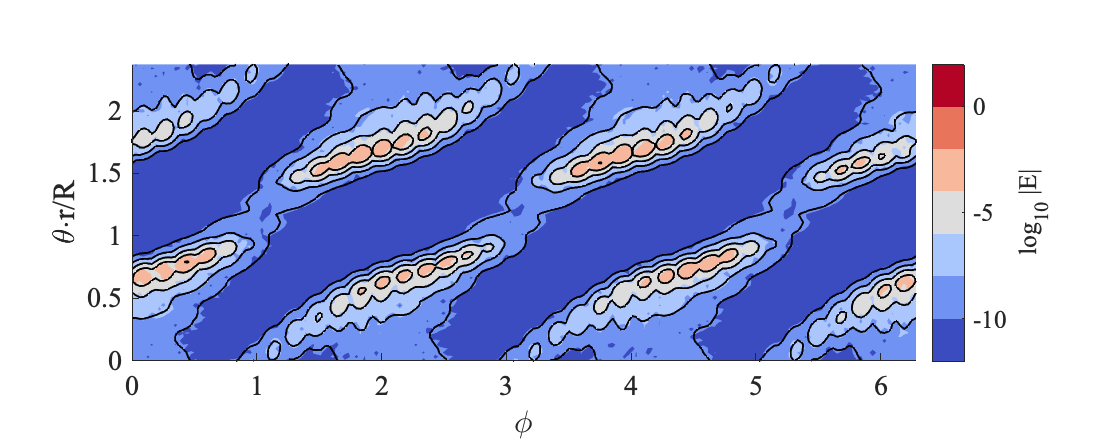}  
  \caption{{Errors for the single layer potential of the modified Helmholtz equation with $\omega=5$ as colored fields and estimates as black contours on
    the toroidal shell, for the trapezoidal rule (above) and the
    Gauss-Legendre rule (below). The contours are drawn for the error
    levels as in Figure \ref{fig:toroidal_shell_contours}.}}
  \label{fig:toroidal_shell_contours_MHsl}
\end{figure}

{We tested the performance of the estimates on different kernels (see \cref{table_kernels}): the behaviour of the error and of the corresponding estimates depends on the decay of the singularity, so the plots will look very similar for kernels with the same decay. In addition to the results for the harmonic double layer potential (\cref{fig:toroidal_shell_contours}), in \cref{fig:toroidal_shell_contours_MHsl} we show also the erros and error estimates for the modified Helmholtz single layer potential with $\omega=5$ and the same discretization and density as shown in \cref{fig:toroidal_shell_geo}.}

%The results, displayed in Fig. \ref{fig:toroidal_shell_contours}-{\ref{fig:toroidal_shell_contours_MHsl}}, are very similar to those in \cref{sec:flat-plane-cutting}. 

\section{Conclusions}

In this paper we have introduced a theoretical and computational
framework for estimating nearly singular quadrature errors in the
evaluation of layer potentials of the form \eqref{eq:layerpot_S2}, for
smooth source geometries that are either one-dimensional curves in
$\mathbb R^2$ or $\mathbb R^3$, or two-dimensional surfaces in
$\mathbb R^3$. This framework is defined for the trapezoidal and
composite Gauss-Legendre quadrature rules, which are two of the most
common choices in the integral equation field. However, generalization
to other quadrature rules is possible with the knowledge of the
remainder function $k_n$ (Elliott et al. derive an expression for
Clenshaw-Curtis quadrature in \cite{Elliott2008}). 

Our work on quadrature error estimates started in
\cite{AfKlinteberg2016quad}. It was extended and improved upon for
one-dimensional curves discretized using composite Gauss-Legendre
quadrature in \cite{AfKlinteberg2018}, e.g. introducing the
root-finding procedure needed for accurate estimation for curved
panels. In \cite{AfKlinteberg2020line},  a so-called singularity swap
quadrature method was introduced for curves in both $\mathbb R^2$ and
$\mathbb R^3$ again based on composite Gauss-Legendre
quadrature, introducing some key ideas for curves in $\mathbb R^3$
that we have explored in this work.

There are three major contributions of the current work:
(1) Showing how error estimation and root finding on one-dimensional
curves can be derived for and applied to the trapezoidal rule,
(2) extending the analysis for complex kernels on curves in $\C^2$ to
real valued kernels on curves in both $\mathbb R^2$ and $\mathbb R^3$,
%%for both the trapezoidal rule and the composite Gauss-Legendre quadrature rule,
(3) deriving the desired quadrature error estimates for two-dimensional
surfaces in $\mathbb R^3$ building on the results on one-dimensional
curves in $\mathbb R^3$. 

%showing how the results on one-dimensional curves in $\mathbb R^3$ can
%be extended in order to derive the desired quadrature error estimates
%for two-dimensional surfaces in $\mathbb R^3$.

%Our work builds on recent development in the contexts of adaptive
%quadrature by expansion \cite{AfKlinteberg2018} and singularity swap
%quadrature \cite{AfKlinteberg2020line}, both of which consider
%one-dimensional curves discretized using composite Gauss-Legendre
%quadrature. Here, we take this development further in two major ways:
%(1) By showing how rootfinding and error estimation on one-dimensional
%curves can be extended to the trapezoidal rule. (2) By showing how the
%results on one-dimensional curves can be extended in order to derive
%estimates for two-dimensional surfaces in $\mathbb R^3$.

As we have shown, our quadrature error estimates perform very well in
actual computations, consistently estimating the error to within one
order of magnitude of the actual value for layer potentials evaluated
over curved surfaces in $\reals^3$. For curves, the estimates are
remarkably precise already for moderate values of discretization
points $n$, even though they are asymptotic estimates.
The focus of this work has not been to
derive upper bounds of the error, even if such bounds would be
desirable. Error estimates without unknown coefficients are
more useful in actual simulations. Evaluation of our
estimates have a low per-point computational cost, since they only
require information from nearby surface grid points (if local
approximations are used). They can therefore be used on the fly in 2D
and 3D simulations to determine e.g. when the regular
quadrature ceases to be sufficiently accurate or what upsamplig rate
should be used.
{The error estimates can also be used to create adaptive
quadrature algorithms, such as was done in 2D \cite{AfKlinteberg2018}, especially needed in 3D applications with multiple particles interacting (e.g. drops \cite{Sorgentone2018,Sorgentone2021}, vesicles \cite{Rahimian2015}, etc.) to provide efficient quadrature methods with error control. }

%The error estimates can also be used to create adaptive
%quadrature algorithms, such as \cite{AfKlinteberg2018}, with the aim
%to provide quadrature methods with error control. 

\section*{Acknowledgements}

L.a.K. would like to thank the Knut and Alice Wallenberg Foundation
for their support under grant no.\ 2016.0410. C.S. acknowledges
support through The Dahlquist Research Fellowship financed by Comsol
AB. A.K.T acknowledges the support of the Swedish Research Council
under Grant No. 2015-04998. 

\appendix
\section{Two lemmas}
\label{app:lemmas}

In the derivation for half-integers in \cref{subsec:halfint}, we use a
result that we prove in \cref{lemma:main}. We start by providing
an intermediate result.

\begin{lemma}
Let $\pbar \ge 1$ be an integer. Then the following holds, 
\begin{equation}
\Gamma(1/2) \prod_{q=1}^{\pbar} \left(\pbar+1/2-q \right) =
%%\Gamma(1/2) \prod_{r=0}^{\pbar-1} \left(r+1/2 \right) =
 \Gamma(\pbar+1/2), 
\label{eqn:gammalemma}
\end{equation}
where $\Gamma(z)$ is the gamma function.
\label{lemma:Gammaproduct}
\end{lemma}

{\em Proof:} First make the simple substitution $r=\pbar-q$ to obtain
\begin{equation}
\Gamma(1/2) \prod_{q=1}^{\pbar} \left(\pbar+1/2-q \right) =
\Gamma(1/2) \prod_{r=0}^{\pbar-1} \left(r+1/2 \right).
\label{eqn:gammalemma_step1}
\end{equation}
We then make use of the relation $z \Gamma(z) = \Gamma(z+1)$,
that holds for all $z \in \C$ (see e.g. Eqn 1.2.1 in \cite{Lebedev1972}).
%\begin{equation}
%z \Gamma(z) = \Gamma(z+1).
%\end{equation}
Particularly, this yields
\begin{equation}
\frac{2\pbar+1}{2} \Gamma\left(\frac{2\pbar+1}{2}\right)=\Gamma\left(\frac{2\pbar+3}{2}\right),
\end{equation}
and so $\frac{1}{2} \, \Gamma(1/2)=\Gamma(3/2)$,  $\frac{3}{2} \, \Gamma(3/2)=\Gamma(5/2)$
etc. 
Starting from (\ref{eqn:gammalemma_step1}), we use this formula repeatedly,
\begin{align*}
\Gamma(1/2) \prod_{r=0}^{\pbar-1} \left(r+1/2 \right) & =
\Gamma(1/2) \cdot (1/2) \prod_{r=1}^{\pbar-1} \left(r+1/2 \right) =\Gamma(3/2) \prod_{r=1}^{\pbar-1} \left(r+1/2 \right) 
= \cdots = \\
&=\Gamma(\pbar-1+1/2) \cdot (\pbar-1+1/2)=\Gamma(\pbar+1/2),
\end{align*}
which yields the desired result.

\begin{lemma}
Let $\pbar \ge 1$ be an integer. Then the following holds, 
\begin{equation}
\frac{1}{\prod_{q=1}^{\pbar} \left(-\pbar-1/2+q \right)}=
%%\frac{(-1)^{\pbar} \sqrt{\pi}}{\Gamma(\pbar+1/2)}=
\frac{1}{\sqrt{\pi}} \Gamma(1-(\pbar+1/2)),
\label{eqn:partintfac}
\end{equation}
where $\Gamma(z)$ is the gamma function.
\label{lemma:main}
\end{lemma}

{\em Proof:} 
We can factor out a negative sign, and use \cref{lemma:Gammaproduct} together
with the fact that $\Gamma(1/2)=\sqrt{\pi}$, to get
\begin{equation}
\frac{1}{\prod_{q=1}^{\pbar} \left(-\pbar-1/2+q \right)} =\frac{(-1)^{\pbar} }{\prod_{q=1}^{\pbar} \left(\pbar+1/2-q \right)} 
=\frac{(-1)^{\pbar} \sqrt{\pi}}{\Gamma(\pbar+1/2)}.
\label{eqn:applylemma1}
\end{equation}

The so called Euler reflection formula  (see e.g. Eqn 1.2.2 in
\cite{Lebedev1972}) holds for all $z \in \C$,
\begin{equation}
\Gamma(z) \Gamma(1-z) = \frac{\pi}{\sin (\pi z)}.
\end{equation}
Specifically, this means that we have
\begin{equation}
\frac{1}{\Gamma(\pbar+1/2)}=\frac{(-1)^{\pbar}}{\pi}
\Gamma(1-(\pbar+1/2)).
\end{equation}
Combining this with (\ref{eqn:applylemma1}), we obtain the desired
result (\ref{eqn:partintfac}).

\section{Error estimates for cartesian and complex formulation of the harmonic double
  layer potential}
\label{app:DoubleLayerEst}

The harmonic double layer potential in two dimensions is given by
\begin{align}
  u(\v x) = \int_{\Gamma} 
\frac{\v \hat{\v n}_y \cdot(\v x-\v y)}{\norm{\v y - \v x}^{2}}  \sigma(\v y) 
\dif S(\v y),
\label{eq:doublelayer}
\end{align}
where $\hat{\v n}_y$ is the outward pointing normal at $\v y \in
\Gamma$.
With the curve parameterized by $\v \gamma(t): \reals \rightarrow \reals^2$, we have $\dif S(\v y)=\norm{\v\gamma'(t)} \dif t$.

Identify the vectors $\v x$, $\v y$ and $\hat{\v n}_y$ in $\reals^2$ by  $z$,
$\tau$ and  $n_{\tau}$ in $\C$. 
%Using also that for two vectors $\v v_1$ and $\v v_2$, identified with
%the two complex points $z_1$ and $z_2$, it holds $\v v_1 \cdot \v
%v_2=\Re (z_1 \bar{z}_2)$ and $|\v v_1-\v v_2|^2=(z_1-z_2)\overline{(z_1-z_2)}$.
We can then write
\begin{equation}
\frac{\v \hat{\v n}_y \cdot(\v x-\v y)}{\norm{\v y - \v x}^{2}} \rightarrow
\frac{\Re \left\{ n_{\tau}\overline{(z-\tau)} \right\} }{(z-\tau)
  \overline{(z-\tau)} } =
\Re \left\{ 
\frac{ n_{\tau}\overline{(z-\tau)} }{(z-\tau) \overline{(z-\tau)} } 
\right\} =
\Re \left\{ \frac{ n_{\tau}}{z-\tau} \right\}.
%\Re \left{n_{\tau}\overline{(z-\tau)} \right} }{a}
%%{ (z-\tau) \overline{(z-\tau)} }
\end{equation}
Now let $\tau(t): \reals \rightarrow \C$ be a complex
parameterization of $\Gamma$, a positively oriented curve. 
The outward pointing normal vector is then given by
$n_{\tau}=-i\tau'(t)/|\tau'(t)|$. Furthermore, the integration element $\dif S(\v
y)=\norm{\v\gamma'(t)} \dif t$ becomes $|\tau'(t)| \dif t=|\dif \tau|$.
This means that we can replace $\dif S(\v y)$ by $\frac{-i}{n_{\tau}}
\dif \tau =\Re \left\{ \frac{-i}{n_{\tau}} \dif \tau \right\}$.
Using also that $\sigma(\tau)$ is a real quantity, we get 
\begin{equation}
u(z)=\Re \left\{ \int \frac{ n_{\tau}}{z-\tau} 
\sigma(\tau) \frac{-i}{n_{\tau}} \dif \tau \right\} = 
\Re \left\{ -i  \int \frac{ \sigma(\tau)}{z-\tau} \dif
  \tau \right\}.
\end{equation}
Hence, we can also write
\begin{align}
u(z)= \Im \left\{ \int \frac{ \sigma(\tau)}{z-\tau} \dif
  \tau \right\} = 
\int  \sigma(\tau)  
\Im \left\{ 
\frac{\dif \tau } {z-\tau} 
\right\} 
\label{eq:doublelayerC}
\end{align}

Now, we want to consider the error estimates that have been derived
for the two different forms.
For the complex form in (\ref{eq:doublelayerC}) we can identify $g(t)$ in 
(\ref{eq:complex_psip}) with $g(t)=\sigma(t)$, if we ignore
taking the imaginary part. For $p=1$, the error estimate 
(\ref{eqn:En-complex-w-kn})  for approximating (\ref{eq:doublelayerC})
then simply reads
\begin{equation}
|E_n^{Complex}| \approx \abs{ \sigma(t_0)  k_n(t_0)}.
\label{eqn:En-complex-p1}
\end{equation} 
Here, $k_n(.)$ is specific for the quadrature rule used and can be
found in (\ref{eq:knTz}) and (\ref{eq:knGL}) for the trapezoidal rule and Gauss-Legendre
quadrature, respectively.

The cartesian form of the double layer potential in
(\ref{eq:doublelayer}) can be written in the form of 
(\ref{eq:base_integral}) with $p=1$ and 
\begin{equation}
f(t)=\v \hat{\v n}_{t} \cdot(\v x-\v \gamma(t)) \, \sigma(t) \norm{\v\gamma'(t)}.
\end{equation}
The error estimate (\ref{eq:base_est_integer}) reads
\begin{align}
  \abs{E_n^{Cartesian}} 
         & \approx
           2 \abs{f(t_0) G(t_0) k_n(t_0) } . 
           \label{eq:En_cart_p1}
\end{align}

Using that 
\begin{align}
 \hat{\v n}_{t}=\frac{(-\gamma'_2(t),\gamma'_1(t))}{ \norm{\v\gamma'(t)}},
\label{eq:nhat}
\end{align}
we can write the product 
\begin{equation}
  \begin{split}
    f(t) G(t) &= \frac{-(-\gamma'_2(t),\gamma'_1(t)) \cdot(\v x-\v
      \gamma(t)) }{2\v\gamma'(t) \cdot \pars{\v x - \v\gamma(t)}}\,
    \sigma(t) \\
    &= \frac{1}{2}\frac{
      \gamma'_2(t)(x_1-\gamma_1(t))-\gamma'_1(t)(x_2-\gamma_2(t)) } {
      \gamma'_1(t)(x_1-\gamma_1(t))+\gamma'_2(t)(x_2-\gamma_2(t)) }
    \sigma(t)
  \end{split}
\end{equation}
We have $R^2(t_0,\v x)=(\gamma_1(t_0)-x_1)^2+(\gamma_2(t_0)-x_2)^2=0$.
From this, we find
\begin{equation}
(\gamma_1(t_0)-x_1)= \pm i (\gamma_2(t_0)-x_2),
\end{equation} 
and at $t=t_0$ the expression above simplifies to
\begin{equation}
f(t_0) G(t_0)= \pm \frac{i}{2} \sigma(t_0).
\end{equation}
Hence, the error estimate  (\ref{eq:En_cart_p1}) simplifies to
\begin{align}
  \abs{E_n^{Cartesian}} 
         & \approx
           \abs{\sigma(t_0) k_n(t_0) },
           \label{eq:En_cart_p1_simpl}
\end{align}
which is identical to the error estimate for the complex kernel in
(\ref{eqn:En-complex-p1}).

In writing the estimate for the complex integral, we did not take into
account that we only consider the
imaginary part of the integral. If we do so, we can write the error
estimate for that integral as $ \abs{\Im \left\{\sigma(t_0) k_n(t_0)
  \right\}}$.

In the last step of deriving (\ref{eq:base_est_integer}), the real part is
skipped from the formula in (\ref{eq:base_est_int_w_Re}).
Using that estimate instead, we would get the error estimate, 
\begin{align}
\abs{ \Re \left\{ \mp i \sigma(t_0) k_n(t_0)  \right\}}=
\abs{ \Im \left\{ \sigma(t_0) k_n(t_0)  \right\}}, 
\end{align}
i.e. the same
error estimate as when we include the imaginary part for the complex
integral. 

The estimate including the imaginary part, include also the node oscillations of the error, and
alllows us to capture even these ``wiggles'' in the error
contours, while the estimates in (\ref{eq:En_cart_p1}) and
(\ref{eqn:En-complex-p1}) instead produce error curves that envelopes
the actual error. 
This is further discussed in section \ref{sec:results_complex}.
See figure \ref{fig:aqbx_plot} for the Gauss-Legende rule
and figure \ref{fig:trapz_complex} for the trapezoidal rule. 

\clearpage
\bibliography{library,library-local}
\bibliographystyle{abbrvnat_mod}

\end{document}